\documentclass{article}
\usepackage{latexsym}
\usepackage{epsfig}  

\begin{document}
%\fontsize{13}{13}\selectfont
\title{\huge Cosmic String, Harvey-Moore Conjecture
 and Family Seiberg-Witten Theory}
\author{Ai-Ko Liu\footnote{Current Address: 
Mathematics Department of U.C. Berkeley}\footnote{
 HomePage:math.berkeley.edu/$\sim$akliu}}

\maketitle

\newtheorem{conj}{Conjecture}[section]
\newtheorem{main}{Main Theorem}[section]
\newtheorem{theo}{Theorem}[section]
\newtheorem{lemm}{Lemma}[section]
\newtheorem{prop}{Proposition}[section]
\newtheorem{rem}{Remark}[section]
\newtheorem{cor}{Corollary}[section]
\newtheorem{mem}{Example}[section]
\newtheorem{defin}{Definition}[section]
\newtheorem{axiom}{Axiom}[section]
\newtheorem{obs}{Observation}
\newtheorem{assump}{Assumption}
\newtheorem{summ}{Summary}
\newtheorem{theo-corr}{Theorem-Corollary}[section]
\newtheorem{axioms}{Axiom:}[section]
\newtheorem{assum}{Assumption} 
\newtheorem{warn}{Warning}

\section{Preliminary}

In this paper we plan to discuss the curve enumeration of certain
algebraic three-folds with additional 
structures of $K3$ fibrations. Our motivation 
to study this question is to solve 
the conjecture of Harvey-Moore [HM1],[HM2] on the numbers of immersed
rational curves in the fibers of a $K3$ fibered Calabi-Yau three-fold.

 In a series of papers [HM1], [Moo], 
the counting of rational curves in a $K3$ fibered
 Calabi-Yau threefold has been tied to R. Borcherd's product formula [Bo1],
 [Bo2] of certain automorphic forms, hyperbolic Kac-Moody algebra [GN2] 
and arithmetic Mirror symmetry [D], [GN1]. It is one of the 
 major predictions of the
 so-called heterotic-type $II$ string-string duality [KV], [HM1], [HM2], [KLM], 
etc. 
 We solve this conjecture based on our theory of family
 Seiberg-Witten invariant [LL1], [Liu1], [Liu2], [Liu4], [Liu5], [Liu6]. 
One major difficulty
to resolve Harvey-Moore conjecture is that their prediction on the 
curve enumeration gives formulae on the
embedded rational curves instead of the stable maps. Thus we have to build a
new foundation to define Harvey and Moore 
``numbers of immersed rational curves'' mathematically instead
of using the standard Gromov-Witten invariants. 
Ideally they are related to the standard Gromov-Witten
 invariants by a multiple covering formula (Please consult the preprint form of
[HM1] page 45, equation 10.4 or section 3.4., equation (20), of [Moo ] 
for its relationship with the pre-potential of
 Gromov-Witten invariants). 
Up to now the algebraic
 definition of ``the number of 
immersed nodal curves'' for general Calabi-Yau three-folds has not been found
 yet, while a symplectic definition of the ${\bf Z}$-valued invariant of ``the
 number of immersed nodal curves''
 for Calabi-Yau three-folds has been announced by Parker-Ionel [IP] recently.

\medskip

  In this paper we interpret these numbers naturally 
as the virtual numbers counting
 nodal curves in the $K3$ fibers, adopting family Seiberg-Witten theory.
 The idea of ``families'' and the ``family invariant'' 
arise naturally as the Calabi-Yau
three-folds involved in the string-string duality are always $K3$ fibered and can
be viewed as one parameter families of K3 surfaces. Thus it is rather nature to 
resolve the problem by using the concept of \underline{family}
 Seiberg-Witten invariant.

 \medskip

 In this paper 
``the numbers of immersed rational curves'' are interpreted as the virtual 
numbers of nodal rational curves in algebraic geometry. Our approach enables
 us to answer, 

\medskip

\noindent (i). In what sense are these virtual numbers ``invariants''?

\medskip

\noindent (ii). Why is the generating function of these virtual numbers
 modular?

\medskip
 
 It will be extremely interesting to
 compare with the symplectic definition of nodal curve invariants by
 the symplectic geometers [IP].

\bigskip

 Given an algebraic surface $M$, the (closure of) the moduli spaces of nodal
curves are often non-smooth and are ill-behaved. It has been a difficult task
to define the invariant of nodal curves which are reduced to geometric
 countings of nodal curves whenever the numbers of nodal curves can be 
understood in the classical sense \footnote{I.e. by imposing geometric 
conditions, say requiring the curves to pass through a finite number
 of generic points, the appropriate 
moduli space of nodal curves can be cut down to a finite number of
 points, potentially with multiplicities. Then the weighted sum gives
 a geometric count of nodal curves.}. The Yau-Zaslow conjecture [YZ] and
 its Calabi-Yau 
three-fold generalization, Gopakurma-Vafa conjecture [GV], rely on a
 mathematical foundation of invariants which count immersed nodal curves
 on $K3$ or on Calabi-Yau three-folds.   
The theory developed in [Liu5], [Liu6] allows
us to give an algebraic definition of
 the ``virtual numbers of immersed nodal curves'' for
 algebraic surfaces. We show that the naive
definition of attaching a virtual fundamental cycle to the appropriate 
moduli space of nodal curves is indeed well defined and is independent to
the complex deformation of $M$, the ill-behavior of the
 geometric cycle of the moduli space of nodal curves, 
or the various choices involved in the
 definition.   
 In this paper, we extend the approach to the cases of 
$K3$ pencils and attach
 nodal curve counting invariants to them.

\medskip

 Let us state the main result of the paper.

\begin{main}\label{main; 1}
 Let $X$ be a simply connected algebraic threefold which has 
 a tamed\footnote{Refer to definition \ref{defin; tame} for its definition.}
 $K3$ fibration structure over ${\bf P}^1$ and let $H^{1, 1}(X, {\bf Z})_f$
 be the group of type $(1, 1)$ cohomology classes of $X$ modulo the
 the subgroup generated by the class of the fiber.

Then 

\medskip

\noindent (i). given a fiber-wise effective 
class $C\in H^{1, 1}(X, {\bf Z})_f$,
 the virtual number of $g$-node
 nodal curves along $X\mapsto {\bf P}^1$ in $C$ is well defined. I.e. the
 virtual number is independent to the complex deformation of the tamed 
$K3$ fibration and choices involved in the definition.

\medskip

\noindent (ii). When $H^{1, 1}(X, {\bf Z})_f\subset H^2(K3, {\bf Z})$ 
is an unimodular lattice, for all $g$
all the virtual numbers ($\not=0$) of ${C^2\over 2}-g+1$-node 
nodal curves in $C$ depend on $C$ through 
$C^2$ only.

\medskip

Once we realize that they depend on the self-intersection numbers $C^2$
 only, we may form the corresponding generating function in terms of
 the formal variable $q$, powered by one-halves of the self-intersection numbers.

\medskip

\noindent (iii). The generating function can be factorized into the
 product of ${q\over\Delta(q)}=\{{1\over \prod_{i\geq 0}(1-q^i)}\}^{24}$ and
 a power series in $q$.

 When $g=0$, the generating function of immersed rational curves is 
reduced to the product of ${q\over\Delta(q)}$ and a $SL_2({\bf Z})$ modular
form of weight $11-{{rank_{\bf Z} H^2(X, {\bf Z})_f}\over 2}$.
\end{main}

 When the $K3$ fibration 
$X$ is a Calabi-Yau three-fold, the main theorem answers affirmatively
 the original conjecture by [HM1] as its special case. 

\begin{cor}
 When \footnote{Consult section \ref{subsection; exam} 
for the definitions of $Z_0$ and
 $W_0$.} $X=Z_0$ (or $W_0$), the $K3$ fibrations considered by
 Harvey-Moore [HM1], then 
the generating function of nodal rational curves is
 $-2({q\over \Delta(q)})\cdot  E_6(q)$ (or 
$-2({q\over \Delta(q)})\cdot  E_4(q)E_6(q)$), where $E_{2k}(q)$ is the
 Eisenstein series of weight $k$.
\end{cor}

 Please compare with (A. 48) and (A. 46) on page 56-57 of the e-print [HM1].

\medskip

 We discuss the implications of our main theorem in the following series
 of remarks.

\begin{rem}\label{rem; a}
The main theorem has a few interesting implications.

(i). From the string theory point of view, the prediction of rational curves
 on $K3$ fibered Calabi-Yau three-folds relies heavily on the string-string
duality and the Calabi-Yau condition. Yet our theorem indicates that the
 $c_1(X)=0$ condition plays no essential role in the enumeration. The only
 factor depending on the Calabi-Yau condition is the $-2$ (see section
\ref{subsection; CY}) in front of the modular objects. 

\medskip

(ii). The above 
result relies heavily on the fact that $dim_{\bf C}X=3$. Our theory
 allows extensions to higher dimensional bases, yet suppose
 the $K3$ fibration is over a higher dimensional base, the result of the
 enumeration will be quite different from the $dim_{\bf C}X=3$ case.

\medskip

(iii). Because the appearance of the Yau-Zaslow [YZ] factor ${q\over
 \Delta(q)}$, people had speculated naturally that the curve counting within
the $K3$ pencils were related to curve counting on $K3$.
Yet the computation done in section \ref{subsection; CY} implies that
 for the Calabi-Yau $K3$ fibrations (Fermat hypersurfaces in weighted 
projective spaces) constructed by the physicists, the 
 curves contributions are completely from the singular fibers (with 
$p_g=0$) closely related to complex multiplications on $K3$ surfaces. 
This is quite opposite to the intuition that the curves 
 should come from $K3$ themselves\footnote{  
 Only after some perturbations of the $K3$ fibrations (consult proposition
 \ref{prop; kk} on page \pageref{prop; kk}), we can make the 
 curve enumerating invariants computable.}.

\medskip

(iv). The unimodularity of the fiberwise cohomology lattice $H^2(X, {\bf Z})_f$
 is the key for the whole expressions (e.g. $-2E_6(q){q\over \Delta(q)}$ or
 $-2E_4(q)\cdot E_6(q){q\over \Delta(q)}$) to be $SL_2({\bf Z})$ modular forms.
 These modular factors $E_6(q)$ or $E_4(q)\cdot E_6(q)$ 
 have encoded the intersection theory of the cosmic strings
 with special divisors in the moduli space of polarized algebraic $K3$.

\medskip

(v). In general we expect that the generating function of the 
Gromov-Witten invariants of (the fiberwise classes of) $X$ with
a large second betti number $b_2(X)$ to be a power series of 
$b_2(X)-1$ variables. What the theorem demonstrates is that 
for the special $K3$ fibrations considered by the string theorists, the
generating function of immersed nodal curves collapses to a power 
series of $q$ alone. 

\medskip

(vi). The ``virtual numbers of nodal curves'' for these
 $K3$ fibered Calabi-Yau threefolds are nothing but the
 Gopakumar-Vafa numbers [GV] of these $K3$ fibered Calabi-Yau threefolds. 
The argument in this paper along with
 the algebraic proof of universality theorem in [Liu5] give the
 direct algebraic geometric definition of Gopakumar-Vafa numbers in terms of
 intersection theory [F] for these tamed Calabi-Yau $K3$ fibrations. 
\end{rem}

 See remark \ref{rem; c} below for some more information when the
 lattice $H^2(X, {\bf Z})_f$ is non-unimodular.

\begin{rem}\label{rem; b}
 The corresponding predictions on the numbers of 
high genus curves of these Calabi-Yau $K3$ fibrations are given by 
\footnote{Modulo certain algebraic manipulations on the modular forms.} [MM].
 Our theorem shows that the generating functions of higher genera nodal curve
 invariants of $Z_0$ and $W_0$ are also factorized into the products of
 $g$ independent cosmic string factors 
$-2E_6(q)$ and $-2E_4(q)E_6(q)$ and some $g$ dependent factors. 

 Our argument in section \ref{section; nodal} indicates that 
these $g$ independent cosmic string
 factors depend only on the $K3$ pencils through lattice theory and 
Howe duality from $Sp(1)$ to $SO(p, 2)$. This provides a beautiful mathematical
 explanation of type $II-A$ heterotic duality\footnote{The reader should notice
 that Borcherds' work on product formula [Bo1], [Ko] 
and its interaction with the string theory of type $IIA$-heterotic duality 
should be also understood in a similar context.}.
 
 Surprisingly the $g$ dependent 
factors encode the enumerative information on a single algebraic $K3$ and
 our theorem asserts that they are independent to the specified $K3$ pencils. 

\medskip

In fact, Kawai [Ka] has computed the $g=1$ predictions for 
 several Calabi-Yau $K3$ fibrations \footnote{Several of these
 $K3$ fibrations have non-unimodular fiberwise cohomology lattices.}. It is
rather non-trivial to observe from his equations (24), (26), (32), (34),
 (55), (57), (80)-(82), (85), (87), (90), etc. that the 
ratios ${\tilde{H}_{\diamondsuit}(\tau)
\over H_{\diamondsuit}(\tau)}$ for $\diamondsuit=A$, $B$, $C$, $\cdots$, 
are all equal to $E_2(q)$ and are therefore pencil independent.

Based on the string theory prediction [MM], [Ka], [CCLM] and Gopakumar-Vafa
 conjecture [GV] and [HST],
 such $g$ dependent terms can be determined explicitly and are
 (quasi) modular forms, expressed
 as polynomials of $E_2(q)$, $E_4(q)$ and $E_6(q)$.
\end{rem}

 It is a very interesting mathematical question to determine the
 universal ($K3$-pencil-independent) $g$-dependent factors without using
string duality, as these factors generate the Gopakumar-Vafa numbers of $K3$.   
 Once this is achieved, it provides a mathematical way to determine the
 Gopakumar-Vafa numbers of all these tamed $K3$ fibered Calabi-Yau threefolds.
 We hope to come back to this subject in the future.

\begin{rem}\label{rem; c}
When we drop 
the assumption (ii). of the main theorem \ref{main; 1} on the unimodularity
 of the fiberwise cohomology $H^2(X, {\bf Z})_f$, our
 technique is still applicable. Our argument implies that the virtual
 number of nodal curves in $C$ depends not only on $C^2$ but also on
 $C$'s pairings with a finite number of cohomology classes. In such a
 situation, the generating function of such virtual numbers has the
 pattern of (a portion of) the multi-variables 
'q expansion' of a Siegel modular form on some 
 Siegel upper half space (or equivalently a Jacobi form).
 On the other hand, the known prediction [Ka], [Ka2] of
 some explicit example from physicists indicates that the answer can be
also coded by classical modular forms \footnote{These modular forms
 may have $q$ expansions with fractional powers.} 
of suitable congruent subgroups of
 $SL_2({\bf Z})$. Potentially this may indicate that there is some interesting
 lifting from classical modular forms to Siegel modular forms which
 relates the string theory prediction and the calculation from the 
mathematical side.  
\end{rem}

\begin{rem}\label{rem; d}
The Yau-Zaslow [YZ] expression ${q\over \Delta(q)}$ 
appearing in the main theorem has been identified [Liu1] in the context
 of the generating function of nodal curves on a single $K3$,
 by using the ${\cal C}^{\infty}$ method (i.e. applying 
the argument of Taubes' ``SW=Gr'' to
 some Kahler families). A purely algebraic identification of the
 Yau-Zaslow factor [YZ] is yet to be found.
\end{rem}

\subsection{The Layout of the Current Paper}\label{subsection; layout}

\bigskip

 The layout of the current paper is as the following.

 In the next subsection \ref{subsection; survey}, we survey the
 Harvey-Moore's conjecture briefly and give a few important references from
 string theory.

\medskip

 In section \ref{section; K3f}, we recall some basic knowledge about the
 cohomologies of $K3$ fibrations. We also give a few examples of $K3$
fibrations, including the original examples studied in Harvey-Moore [HM1].
 In section \ref{subsection; KV}, we review the Kawamata-Viehweg covering
 trick briefly and use it to construct $K3$ fiber bundles from tamed $K3$
fibrations. 
 In section \ref{subsection; iso}, 
we relate these examples with the concept of complex multiplications of
 $K3$ and point out that these families are indeed iso-trivial $K3$ families.

 In section \ref{section; cosmic}, 
we introduce the concept of cosmic string
 (brane) and study its relationship with the family Seiberg-Witten invariants
 of $K3$ fibrations. In section \ref{subsection; general}, we point out that
 the family invariants of the $K3$ fibrations is closely related to the 
 Weil-Peterson Volume of the cosmic brane map. In section \ref{subsection; 
pathetic},
we point out in the $dim_{\bf C}B=1$ case why multiple-coverings of $-2$ curves 
 which potentially damage the curve counting do no harm to the enumeration result
 at all. The vanishing result we derive
 can be viewed as a warm up of the similar vanishing result for type $II$
 curves in section \ref{subsection; rel}.

 In section \ref{subsection; CY}, we study how does the Calabi-Yau 
 condition affect the family invariant. We also 
derive a simple defect relationship
 for the family invariant. This formula enables to understand why the
 family invariant of iso-trivial $K3$ fibrations comes solely 
from the singular fibers and relate this with a fractional bubbling off
 phenomenon of cosmic string maps.

\medskip

 In section \ref{section; proof}, we apply the general machineries of [Liu5] to
the relative setting and study the virtual numbers of nodal curves along
 the $K3$ fibrations. Even though our technique is also applicable to
 $K3$ fibrations with higher dimensional bases, 
we concentrate on the case of one dimensional base
 to simplify the discussion.  As the so-called ``universality theorem''
 has been interpreted in [Liu7] as a non-linear version of enumerative 
 Riemann-Roch formula, the discussion in section \ref{section; proof} can
 be viewed as the ``family version'' of our enumerative Riemann-Roch formula.

\medskip

 This section has been the backbone of the paper. Because much of the focus
  in section \ref{subsection; rel}
 is parallel to the absolute $\tilde{B}=pt$ case, we will refer to the
 long paper [Liu5] frequently and emphasize the parts which need to be
 adjusted. The reader who is interested in getting to the technical details
 should go back to [Liu5] and [Liu7].

\medskip

In section \ref{section; nodal}, we combine the various 
results from section \ref{section; cosmic}, section \ref{section; proof}
 to prove the main result of the paper. By combining the ideas of
 ``cosmic string'' [GSYV], family Seiberg-Witten theory, and the Howe duality,
 we derive the main theorem assuming that the fiberwise cohomology lattice
 is unimodular. It turns out the intersection theory of the cosmic string
 with the special divisorial cycles in the moduli space of lattice 
polarized $K3$s plays a vital role here.

To reduce the length of the paper, the corresponding 
results for non-unimodular lattices will be discussed elsewhere. 

\medskip

 The author would like to express gratitude to 
S.T. Yau, C.H. Taubes, H. H. Wu for their interest in the work. He also 
wants to thank S. Givental, Kefeng Liu, J. K. Yu, M. Marino, A. Klemm, 
E. Zaslow for helpful 
 discussions.

\subsection{A Short Survey on Harvey-Moore Conjecture}\label{subsection; survey}

\bigskip

 In the section, we give a brief introduction to the conjecture of Harvey-Moore
 regarding the counting of rational curves in a Calabi-Yau $K3$ fibrations.

The concept of type $IIA-Heterotic$ String-String duality was first
 raised by C. Vafa and S. Karchu in [KV]. It was then discovered that every
 Calabi-Yau threefold which had a ``heterotic dual'' was $K3$ fibered. The
 reader may consult [KLM] for examples of $K3$-fibered Calabi-Yau threefold.

\medskip

 By using the concept of type $IIA-Heterotic$ duality [KV], Harvey and Moore 
[HM1], [HM2] were
able to relate some super-conformal theory on a
$K3$ fibered Calabi Yau threefold to certain super-conformal
 theory constructed on $K3\times {\bf T}^2$, known as heterotic string theory.
 It turns out the super-conformal theory depends on some continuous moduli
 depending on the complex structures of $K3$ and $T^2$.

  It is well known that the moduli of algebraic $K3$ surfaces forms
  some arithmetic quotient of a type $IV$ bounded symmetric domain.

\medskip

 Borcherds [Bo1], [Bo2] constructs certain automorphic form on such symmetric
 domain which allows a product formula similar to the Weyl-Kac character 
formula of a affine Kac-Moody algebra [GN2]. In certain cases, the product formula
 of such automorphic form can be realized as the Weyl-Kac-Borcherds 
character formula of certain hyperbolic generalized Kac-Moody algebra.
 The exponents of the product formula corresponds to dimension of the
 root space ${\bf g}_{\alpha}$.

$$e^{2\pi\sqrt{-1}\rho\cdot y}\prod_{r>0}(1-e^{2\pi\sqrt{-1}r\hat{\cdot}y})^{
 c({r^2\over 2})},$$

   One of Borcherds' major discoveries is that the generating function 
$\sum_mc(m)q^m$ of the dimensions of the
 root spaces is nothing but the $q$ expansion of
 some classical modular form on the upper half plane. 

  It turns out that the same type of 
product formula appears naturally in the framework of
 Harvey and Moore [HM1].

 As a part of string theory datum,
 the super-conformal theory also depends on the maximal
 torus of ${\bf E}_8\times {\bf E}_8$.   Based on the calculation of
 Elliptic genera of $K3$ on the heterotic side, they were able to compute
 such automorphic form and the exponents of the product factors explicitly.
 Based on the conjecture that Heterotic String compactified on $K3\times T^2$ 
is dual to the type $IIA$ string
 on certain $K3$ fibered Calabi-Yau manifold, 
they identified the hypothetical $g=0$ invariants of 
embedded rational curves \footnote{Closed related to the $g=0$ Gromov-Witten 
invariants.} of those $K3$ fibered Calabi-Yau threefolds as the exponents
 of the product expansion of the 
specific automorphic form on the complex moduli of the $K3$ surface.

 As a consequence, based on type $IIA$-heterotic duality, Harvey and Moore
 predict that the generating function of 
the number of embedded rational curves on such Calabi-Yau $K3$ fibered 
threefold must be modular.

 The following conjecture is extracted from the predictions of 
Harvey-Moore to $K3$ fibered Calabi-Yau threefolds $Z_0$ and $W_0$.

\begin{conj}\label{conj; HM}
 Let $X$ be a tamed $K3$ fibered Calabi-Yau threefold and let $n_C$ denote
the virtual number of rational curves in the fiberwise 
class $C\in H_2(X, {\bf Z})_f$. Suppose that the fiberwise intersection
 lattice ${\bf M}\cong H^2(X, {\bf Z})_f$ is unimodular,
then $n_C$ depends on $C$ through the
 self-intersection number $C^2=2\delta-2$ 
along the fiber and the generating function
 of $n_\delta=n_C, C^2=2\delta-2$  

$$\sum_{g\geq 0} n_{\delta} q^{\delta}$$
  is the $q$ expansion of an explicitly constructed modular form. 
\end{conj}

 If $X=Z_0$, then the modular form is equal to $-2{q\over \Delta(q)}E_6(q)$.
 If $X=W_0$, then the modular form is equal to $-2{q\over \Delta(q)}E_4(q)E_6(q)$.

 The modular forms presented here differ from these from [HM1] up to a minus sign.
 This overall sign difference is due to the different conventions on the 
orientations of the moduli space of curves. Please consult the surveys by 
Moore [Moo] and Kontsevich [Ko] about this conjecture.
  
 The above modularity conjecture of Harvey-Moore is a deep phenomenon closely 
related to the arithmetic mirror conjecture [GN1], [GN2], [D] and Borcherds' 
work on product formula. The high genera extension of the conjecture has
 been formulated by Marino and Moore [MM], and is refined by [Ka], [Ka2], 
 [KY], [HST], etc.

 In the current paper, we only deal with tamed $K3$ fibrations. It is 
very interesting to extend the above study to the cases with
 normal-crossing singular fibers or other singular fibers. 

  They also predicts the existence of some hyperbolic Kac-Moody algebra acting
 on the moduli space of ``BPS'' states such that the cited product formula
 is indeed its character formula [HM2].

\section{The Topology of K3 Fibrations and the Examples}\label{section; K3f}

\bigskip

In this section, we study the topology of $K3$ fibrations and 
introduce a few concepts.

Recall the following definition on $K3$ surfaces,

\begin{defin}\label{defin; K3}
 An algebraic surface $M$ over ${\bf C}$
is an algebraic $K3$ surface if
\noindent (i). it is simply connected.
\noindent (ii). its first Chern class vanishes.
\end{defin}

If $M$ is a simply connected Kahler surface with $c_1(M)=0$, then the
 corresponding complex surface may be non-algebraic.

\begin{defin} \label{defin; fibr}
 Let $X_0$ be a smooth algebraic manifold. $X_0$ is said to carry a
 $K3$ fibration structure if there exists a surjective algebraic morphism 
$\pi_0:X_0\mapsto B_0$ such that the generic fibers are smooth $K3$ surfaces.
\end{defin}

 In this paper, we assume additionally that the base $B_0$ is
smooth.

\medskip

 While the regular fibers are $K3$, the degenerated fibers may contain
 singularities.  For enumerative purpose, we impose additional restrictions on
 the types of
 singular fibers in the given $K3$ fibration.

\begin{defin}\label{defin; tame}
  A $K3$ fibration is said to be tamed
if it satisfies the following conditions.

\medskip

\noindent (i). The singular fibers are all irreducible and reduced.

\medskip

\noindent (ii). In a given singular fiber, there are at most a finite number of
 isolated singularities.

\medskip

\noindent (iii). The classical monodromy operator around each singular fiber is of 
finite order. 
\end{defin}

 The primary goal of our paper is to study the curve counting of 
such tamed $K3$
fibrations.

Our main focus is the $dim_{\bf C}X_0=3$ case. It will be shown in our
discussion that $dim_{\bf C}X_0$ plays a very crucial role in proving
 the main theorem. We also assume that
$h^{2,0}(X_0)=0$ as in the case of usual Calabi-Yau three-folds.

\medskip

 Let us recall some basic facts about the cohomology of $K3$ surfaces. 

Let $M$ be a $K3$ surface. It is well known that $c_1(M)=0$ and $\int_M
 c_2(M)=24$.
 Because simply connectedness (i.e. $\pi_1(M)=\{1\}$)
 and Poincare duality, it immediately
follows that $H^1(M, {\bf Z})=H^3(M, {\bf Z})=\{0\}$. On the other hand,
 the middle cohomology $H^2(M, {\bf Z})$ is an even unimodular lattice
 of rank $22$ which can be decomposed into $3{\bf H}\oplus 2(-{\bf 
E}_8)$, usually called the $K3$ lattice.

 The symbol ${\bf H}$ denotes the standard rank two 
even hyperbolic lattice while ${\bf E}_8$ 
denotes
 the $E_8$ lattice, the root lattice of the Lie algebra of $E_8$. 
 In this paper, we denote the $K3$ lattice by ${\bf L}$. 

\bigskip

 Because $X_0$ is algebraic and therefore Kahler, 
the cohomologies of $X_0$ is decomposed 
into the
bi-graded vector spaces $H^{p, q}(X_0, {\bf C})$, $p, q\in {\bf N}\cup \{0\}$. 
We have,

 $$H^1(X_0, {\bf C})\cong H^{1,0}(X_0, {\bf C})\oplus H^{0,1}(X_0, {\bf 
C}),$$

 $$H^2(X_0, {\bf C})\cong H^{2,0}(X_0, {\bf C})\oplus H^{0,2}(X_0, {\bf 
C})\oplus H^{1,1}(X_0, {\bf C}),$$

 respectively.

Because $X_0$ is algebraic, $H^{1,1}(X_0, {\bf C})$ is at least one 
dimensional. It contains the subgroup generated by
 the ample polarization of $X_0$.

 Pick an arbitrary regular value $z\in B_0$ of the fibration, $i_z: 
\pi^{-1}(z)\subset X_0$ induces an
inclusion of the algebraic $K3$ fiber over $z$ into $X_0$.

The induced morphism $$i^{\ast}_z: H^2(X_0, {\bf C})\mapsto 
H^2(\pi^{-1}(z), {\bf C})$$
preserves the gradation of Hodge decomposition and induces morphisms

$$H^{1,1}(X_0, {\bf C})\mapsto H^{1,1}(\pi^{-1}(z), {\bf C})$$
 and $$H^{2,0}(X_0, {\bf C})\mapsto H^{2,0}(\pi^{-1}(z), {\bf C})$$ on 
their $(1,1)$ and $(2, 0)$ components.

 Recall the Hodge decomposition of the middle cohomologies of $K3=M$,
$$H^2(M, {\bf C})=H^{2,0}(M, {\bf C})\oplus H^{1,1}(M,
{\bf C})\oplus H^{0,2}(M, {\bf C}))$$
 with hodge numbers $h^{1,1}=20, h^{2,0}=h^{0,2}=1$.

Given an algebraic $K3$ surface $M$, both $H^{1,1}(M, {\bf C})$ and
 $H^2(M, {\bf Z})$ can be viewed as subgroups of $H^2(M, {\bf C})$.
 Then the Picard lattice 
$H^2(M, {\bf Z})\cap H^{1,1}(M, {\bf C})$ is a sub-lattice of
 $H^2(M, {\bf Z})$, which may not be unimodular. It consists of the
 integral cohomology classes of $M$ which are of type $(1, 1)$.

Because $M$ is algebraic, $\rho_M=
rank_{\bf Z} H^2(M, {\bf Z})\cap H^{1,1}(M, {\bf C})$
is nonzero. In general, $\rho_M$ is bounded within the range 
$1\leq \rho_M\leq 20$.

The assumption $h^{2,0}(X_0)=0$ implies the map
 $$H^{2,0}(X_0, {\bf C})\mapsto H^{2,0}(\pi_0^{-1}(z), {\bf C})$$  to be 
null.

In general, the map
$$H^{1,1}(X_0, {\bf C})\mapsto H^{1,1}(\pi_0^{-1}(z), {\bf 
C})$$ 

is neither injective or surjective.

On the other hand, $h^{2,0}(X_0)=0$ implies that
 $H^2(X_0, {\bf Z})$ is a sub-lattice of $H^{1,1}(X_0, {\bf C})$.

The image of the map 
$H^{1,1}(X_0, {\bf C})\mapsto H^{1,1}(\pi^{-1}(z), {\bf C})$
is then generated by the integral
divisor classes in $\pi_0^{-1}(z)$ which are the restriction of
 divisor classes on $X_0$ to the fiber $\pi_0^{-1}(z)$.

The kernel is generated by the divisors in $X_0$ which
 restricts trivially to $\pi_0^{-1}(z)$.  When $dim_{\bf C}X_0=3$,
 $\pi_0^{-1}(z')$, for all 
$z'\in B_0$ determines a unique divisor class $[F]$ which
satisfies $[F]\cdot [F]=0$. In particular, it
lies in the kernel of the restriction map.

In the following, we prove that the kernel is exactly ${\bf C}[F]$.

\begin{prop}\label{prop; kernel}
Assuming that $dim_{\bf C}X_0=3$, $h^{2,0}(X_0)=0$, 
then the kernel of the restriction map
$$H^{1,1}(X_0, {\bf C})\mapsto H^{1,1}(\pi_0^{-1}(z), {\bf C})$$
 is generated by the fiber class $[F]$.
\end{prop}\label{prop}

\noindent Proof: Take $p_0: B_0\mapsto pt$ to be the constant map 
and consider the Leray 
spectral sequence 
associated with the composition
 $p_0\circ \pi_0$.  The $E_2$ term of the spectral sequence
 is 

$${\cal R}^i(p_0)_{\ast}({\cal R}^j(\pi_0)_{\ast}{\bf Z}), $$
which converges to ${\cal R}^{i+j}(p_0\circ \pi_0)_{\ast}{\bf Z}$.

 We pay our attention to the $i+j=2$ piece above. 
The graded two piece of $E_2$ spectral sequence
 consists of three components,
 ${\cal R}^2(p_0)_{\ast}({\cal R}^0(\pi_0)_{\ast}{\bf Z})$, 
 ${\cal R}^1(p_0)_{\ast}({\cal R}^1(\pi_0)_{\ast}{\bf Z})$,
 and ${\cal R}^0(p_0)_{\ast}({\cal R}^2(\pi_0)_{\ast}{\bf Z})$.
 Because the regular fibers are simply connected, the sheaf 
 ${\cal R}^1\pi_{\ast}{\bf Z}$ is supported upon the singular values
 of $\pi_0:X_0\mapsto B_0$, which is a finite set of points in $B_0$.
Therefore ${\cal R}^1 p_{0\ast}\bigl({\cal R}^1\pi_{\ast}{\bf Z}\bigr)=0$.

 In particular, this implies that all but one dimension of the cohomology 
$H^2(X_0, {\bf Z})={\cal R}^2(p_0\circ \pi_0)_{\ast}({\bf Z})$
comes from the monodromy invariant cohomology of the smooth 
fibers. On the other hand, the extra generator
 can be identified with $\pi_0^{\ast}([B_0])=F$.  This proves the
 proposition.   $\Box$

\bigskip

  By our assumption, $H^2(X_0, {\bf C})$ is of purely type $(1,1)$. Thus, all
 the integral second cohomology classes can be represented by
 divisors.  On the other hand, by Poincare duality $h^{2,0}(X_0)=0$ also implies 
 that $H^4(X_0,{\bf Z})$ is of type $(2, 2)$.

The following lemma confirms that all the classes in $H^4(X_0, {\bf Z})_{free}$ 
can be represented by holomorphic curves in $X_0$.

\begin{lemm} \label{lemm; effective}
Let $c\in H^4(X_0, {\bf Z})$ be a non-torsion
 integral degree four cohomology class, then
there exists an positive integer $k$, such that $kc$ is represented by
an integral combination of holomorphic curves in $X_0$. 
\end{lemm}

\noindent Proof of lemma \ref{lemm; effective}: 
Because $H^2(X_0, {\bf Z})_{free}\subset H^{1, 1}(X_0, {\bf C})$, any
 element in $H^2(X_0, {\bf Z})_{free}$ can be represented by an integral
linear combination of effective divisors in $X_0$.

Firstly
 take the ample polarization 
$[\omega]$ of $X_0$. Then the Hard Lefschetz theorem implies that

$$\cup[\omega]:H^2(X_0, {\bf C})\mapsto H^4(X_0, {\bf C})$$ induces
 an isomorphism. The morphism $\cup [\omega]$ is of type $(1, 1)$.
 This implies that all the elements in
 the free module $H_2(X_0, {\bf Z})_{free}$ can be represented by
the intersection of a
divisor class with the ample polarization class.
Recall that any divisor can be written as the formal difference of
 effective divisors. Moreover, by Bertini theorem there exists a $k\in {\bf N}$ 
such that the representatives from high multiples ($\gg k$) of ample divisor
 classes
 can be made to have proper intersections with any given effective
divisor.

In particular, this implies that a high multiple ($\gg k$) 
of any element in
 $H_2(X_0, {\bf Z})_{free}$ can be expressed as a difference of
 effective curve classes.

\medskip

 If the element $c$ is torsion, a high multiple of $c$ vanishes. The
 statement in the lemma holds trivially.
$\Box$

\bigskip

 Define ${\cal C}_{X_0}\in H_2(X_0, {\bf Z})_{free}$ 
to be the curve cone of all the integral 
classes representable by effective combinations of holomorphic curves in 
 $X_0$. Then the elements in ${\cal C}_{X_0}$ generates $H_2(X_0, {\bf R})$.

 The image of ${\pi_0}_{\ast}:{\cal C}_{X_0}\mapsto 
 H_2(B_0, {\bf Z})\cong {\bf Z}$ 
generates a monoid. Because of lemma \ref{lemm; effective},
 we know that this monoid is not trivial.
If the image is equal to ${\bf N}[B_0]\subset H_2(B_0, {\bf Z})$, then
 we take a class $c\in {\cal C}_{X_0}$ which maps onto $[B_0]$.

 By assumption, there exists a holomorphic curve $\Sigma$ 
representing $c$. We 
calculate the intersection pairing $\Sigma \cdot \pi_0^{-1}(z), z\in B_0$, by
 the projection formula, 
 $(c, F)_{X_0}=((\pi_0)_{\ast}c, [B_0]_{coh})_{B_0}$. Here $[B_0]_{coh}$
denote the oriented generator of $H^2(B_0, {\bf Z})$.

 It implies that $\Sigma
\cap \pi_0^{-1}(z)$ is one point, for all $z\in B_0$. In other words, $\Sigma$
 realizes a holomorphic
 section $B_0\mapsto X_0$.

In the general situation, the image of ${\cal C}_{X_0}$ is generated by
 $m[B_0]$, for some $m\in {\bf N}$. Then
a multi-section of $B_0\mapsto X_0$ can be found
 by a class $c\in {\cal C}_{X_0}, (\pi_0)_{\ast}(c)=m[B_0]$.

\bigskip

 Recall the Yukawa cubic intersection pairing on threefolds,
 $$H^2(X_0, {\bf Z})\otimes H^2(X_0, {\bf Z})\otimes H^2(X_0, {\bf Z})\mapsto
 {\bf Z},$$ 

 by sending $a\otimes b\otimes c$ to $\int_{X_0}a\cup b\cup c$.

 By fixing one of the arguments to be $F$, it induces a 
 quadratic pairing on $H^2(X_0, {\bf Z})$,

$$Q:H^2(X_0, {\bf Z})\otimes H^2(X_0, {\bf Z})\mapsto {\bf Z}.$$ 

 The pairing induces an intersection form on the ${\bf Z}$ module
 $H_2(X_0, {\bf Z})_{free}$, which degenerates on ${\bf Z}F$.  
 Given a regular fiber
 $\pi_0^{-1}(z)$, the pairing is compatible with Poincare duality in the
 following sense:

\medskip

 Given any divisor class $[D]$ in $X_0$, its restriction on 
 $\pi_0^{-1}(z)$ induces a divisor class on this fiber. As the fiber
is smooth of complex dimension two, it is a curve class in 
 $\pi_0^{-1}(z)$.

 Then the quadratic pairing 

 $$Q([D]_1, [D_2])=([D]_1, [D]_2, F)_{X_0}=([D]_1|_{\pi^{-1}_0(z)},
 [D]_2|_{\pi^{-1}_0(z)})_{ \pi^{-1}_0(z)}$$ induces the
 intersection pairing on the kernel of  

$${\pi_0}_{\ast}:H_2(X_0, {\bf Z})_{free}\mapsto H_2(B_0, {\bf Z})\cong
 {\bf Z}$$

 through Poincare duality on $\pi_0^{-1}(z)$.

 The pull-back
 map of the inclusion $i_z:\pi^{-1}_0(z)\mapsto X_0$
 induces an embedding of the lattice
$H^2(X_0, {\bf Z})_{free}/Ker(i_z^{\ast})$ into the $K3$ lattice ${\bf L}$.
 The spectral sequence computation in proposition \ref{prop; kernel} implies that
 $F$ is primitive in $H^2(X_0, {\bf Z})$, thus 
 $Ker(i_z{\ast})={\bf Z}F$.

 The Hodge index theorem implies that the positive eigenspace of 
 $Q$ is at most one dimensional. On the other hand, the ample polarization 
$[\omega]$
 restricts to an ample polarization over each fiber $\pi^{-1}_0(z)$ and
 its self-intersection pairing is positive definite.
In particular, the lattice $H^2(X_0, {\bf Z})_{free}/{\bf Z}F$ is
of signature $(1, b^2(X_0)-1)$.
 As a consequence, one gets an upper bound on $b_2(X_0)$.

\begin{lemm}\label{lemm; bound}
The second Betti number of $X_0$ does not exceed $20$.
\end{lemm}

\noindent Proof: It is because the lattice $H^2(X_0, {\bf Z})_{free}/{\bf Z}F$ 
embeds into ${\bf L}$,  while the
 negative eigenspace of ${\bf L}$ is of $19$ dimension. $\Box$

 When $z$ varies over the regular values of $\pi_0$ in $B_0$, the fibers
  $\pi^{-1}_0(z)$ form a one parameter family of
 algebraic $K3$s with varying complex structures.

 In general, the Picard lattice of $\pi^{-1}_0(z)$ changes with respect
 to $z$ and their sizes depend on how $H^2(\pi^{-1}_0(z), {\bf Z})$ intersect
 with the type $(1, 1)$ subspaces $H^{1,1}(\pi^{-1}_0(z), {\bf C})
\subset H^2(\pi^{-1}_0(z), {\bf C})$.  In general, it is a very 
delicate Diophantine phenomenon and fiberwise Picard lattices within the family
 can jump randomly. 
 On the other hand, the inclusion
 of $H^2(X_0, {\bf Z})_{free}/{\bf Z}F$ into
 $H^2(\pi^{-1}_0(z), {\bf Z})$ also gives an lower bound on
 these Picard lattices as all classes in $H^2(X_0, {\bf Z})_{free}/{\bf Z}F$
 are of type $(1, 1)$. 

 Recall the following well known lemma in lattice theory,

\begin{lemm}\label{lemm; uni} 
Let ${\bf M}\subset {\bf L}$ be an inclusion of ${\bf M}$ into an unimodular
 lattice ${\bf L}$. 
 The following conditions are equivalent,

\medskip

\noindent (i). The lattice ${\bf M}$ is unimodular.

\medskip

\noindent (ii). ${\bf M}^{\bot}$ is unimodular, where ${\bf M}^{\bot}$ is formed by
 all elements $\in {\bf L}$ which pair trivially with ${\bf M}$.

\medskip

\noindent (iii). ${\bf M}\oplus {\bf M}^{\bot}={\bf L}$.
\end{lemm}

  In our situation, the $K3$ lattice ${\bf L}$ is even-unimodular.
  If ${\bf M}=H^2(X_0, {\bf Z})_{free}/{\bf Z}F$ is also unimodular,
 it is an indefinite even unimodular lattice of signature $(1, b_2(X_0)-1)$.
 Such lattices can be classified easily by Hesse-Mankowski 
 theorem (see e.g. [Se]).

 Then ${\bf M}$ must be isomorphic to either 
${\bf H}$, ${\bf H}\oplus (-{\bf E}_8)$
 or ${\bf H}\oplus (-2{\bf E}_8)$.  Some of the interesting examples
 from $K3$ fibrations in 
 string theory give intersection lattices in the list. However,
 one does not expect ${\bf M}$ for an arbitrary $K3$ fibration 
to be always unimodular. We will give a few
explicit examples later in section \ref{subsection; exam}.

\bigskip

Because $F^2=0$, consider the $F$-chain complex
formed by

$$0\mapsto H^2(X_0, {\bf Z})\stackrel{\cup F}{\longrightarrow}
 H^4(X_0, {\bf Z})\stackrel{\cup F}{\longrightarrow}
 H^6(X_0, {\bf Z})\mapsto 0.$$

The $F$-cohomology is defined to be $Ker(\cup F)/Im(\cup F)$. We
 explain briefly what does this cohomology group measure.

The space 
$Ker(\cup F)\subset H^4(X_0, {\bf Z})$ 
measures those curve classes (through Poincare duality) which
are induced from the fibers. On the other hand, the image of $F$ from
 $H^2(X_0, {\bf Z})$ are those curve classes which can be constructed by
 intersecting an divisor class $\in H^2(X_0, {\bf Z})$ with $F$. 

 The $F$cohomology is of pure torsion. If the above 
$F$-cohomology is non-trivial, it measures the discrepancy of those fiberwise
 curve classes which are not induced from an intersection of divisors.
 The following remark is a simple consequence of the previous discussion.

\begin{rem}\label{rem; cohomology}
 The following conditions are equivalent,

\medskip

\noindent (i). The lattice $H^2(X_0, {\bf Z})_{free}/{\bf Z}F$ is unimodular. 

\medskip

\noindent(ii). The $F$-cohomology is trivial

\medskip

\noindent (iii). All the holomorphic curves in the fibers of $X_0\mapsto B_0$ are
induced from the intersection of an effective divisor of $X_0$ with some fiber
representing $F$.
\end{rem} 

\bigskip

\subsection{\bf The Construction of Tamed $K3$ fibrations} 
\label{subsection; exam}

In this subsection, we review several examples of $K3$ fibrations which
 satisfy the conditions introduced earlier.  

 The first example is the well-known Lefschetz pencils in ${\bf P}^3$.

\medskip

  Consider ${\bf P}^3$ and the degree four irreducible hypersurfaces 
$\in H^0({\bf P}^3, {\cal O}(4))$ in ${\bf P}^3$. 

 Following Lefschetz,
 one can choose a generic pencil in ${\bf P}(H^0({\bf P}^3, {\cal O}(4))$
such that the quartics within this pencil contain at most ordinary double points. 
 Blowing up the base locus of the pencil of quartics in ${\bf P}^3$, one
 gets a rational algebraic 
manifold $Y_0$, which has a structure of $K3$ fibration over ${\bf P}^1$.

The following lemma identifies the cohomologies of $Y_0$.

\begin{lemm}\label{lemm; Y_0}
Let $Y_0$ be the algebraic threefold constructed above, then $Y_0$ is
 simply connected and $b_2(Y_0)=2$. The hodge numbers of $Y_0$ are given by
 $h^{0,0}(Y_0)=1, h^{1,0}(Y_0)=h^{0,1}(Y_0)=h^{2,0}(Y_0)=h^{0,2}(Y_0)=0$,
 $h^{1,1}=2, h^{3,0}(Y_0)=h^{0,3}(Y_0)=0, h^{2,1}(Y_0)=h^{1,2}(Y_0)=33$. The
other hodge numbers can be calculated by hodge duality 
$h^{3-p, 3-q}(Y_0)=h^{p, q}(Y_0)$.
\end{lemm}

\noindent Proof of lemma \ref{lemm; Y_0}: 
The computation is elementary. Nevertheless for completeness, we give a
 simple argument.  Firstly, the
 exceptional locus of $Y_0\mapsto {\bf P}^3$ has a ${\bf P}^1$ bundle
 structure over an smooth algebraic curve.  By adjunction equality,
 $2g-2=\int_{4[H]}[4H]|_{[4H]}\cup [4H]_{[4H]}=\int_{Y_0}
[4H]\cup [4H]\cup [4H]=4\cdot 4\cdot 4=64$. 
 Thus $g=33$. The blowing up procedure 
replaces the curve (the blowing up center) by a ${\bf P}^1$ bundle over 
the curve, which
 changes the Euler number by $-66$. On the other hand, 
$h^{i,0}(Y_0)=h^{i,0}({\bf P}^3)=0$, $i\geq 0$, because the birational morphism
$Y_0\mapsto {\bf P}^3$ does not change the pluri-genera.  Finally, the
 exceptional divisor contributes $1$ to the second Betti number 
$b_2(Y_0)=h^{1, 1}(Y_0)$. Then all the other hodge numbers can be
 derived easily. $\Box$

 One can determine the intersection matrix of $H^2(Y_0, {\bf Z})$
 as following. Firstly, $b_2(Y_0)=2$ implies that $[H]$ and $[E]$ generate
 $H^2(Y_0, {\bf Z})$.
 The hyperplane class $[H]$ induces a relative
 ample polarization, $[H]\cdot [H]\cdot F=4$.  On the other hand,
 $F=4[H]-[E]$. From this, we derive the intersection pairing
 between $[H]$ and $[E]$ by using $F^2=0$.  It follows that
 the fiber curve classes which are
 the restriction of the divisor classes of $Y_0$ form a one dimensional
 lattice generated by $[H]\cdot F$. The lattice is isomorphic to
 $(4)$.
 
 In this example, the lattice is not unimodular.

\medskip

 By the construction of the Lefschetz fibration, a simple calculation concludes
that,

\begin{prop}
 There are $110$ singular fibers in the $K3$ fibration $Y_0\mapsto {\bf P}^1$,
 and each singular fiber is a singular orbifold $K3$ surface with an isolated
 ordinary double point.
\end{prop} 

\noindent Proof: 
 The conclusion about ordinary double points is straightforward from the
 construction of Lefschetz fibrations.
 The number of singular fibers is determined by the Euler number of
 $Y_0$ and the Euler numbers of the singular fibers. The Euler number of
 a $K3$ fiber bundle
 over ${\bf P}^1$ is $48$. Yet the space $Y_0$ has Euler number 
$4+(-64)=-62$. This implies that the discrepancy of Euler numbers is
 $110$. This implies that there are $110$ singular fibers each with
 Euler number $23$. $\Box$

\bigskip

 The second example and the third examples are constructed by string theorists 
[HM1], [KLM].

  Consider the four dimensional weighted projective spaces
${\bf P}_{24}(1, 1, 2, 8, 12)$,  ${\bf P}_{84}(1, 1, 12, 28, 42)$  
and the Fermat hyper-surfaces in the anti-canonical linear systems

   $$f=x_0^{24}+x_1^{24}+x_2^{12}+x_3^3+x_4^2=0,$$
and

 $$g=x_0^{84}+x_1^{84}+x_2^7+x_3^3+x_4^2=0.$$

These hypersurfaces define singular Calabi-Yau threefolds.  The
 projections to $(x_0, x_1)$ give birationally defined
 rational maps to ${\bf P}^1$ which
 are well-defined away from the curve defined by 
$x_0=x_1=0,  x_2^{12}+x_3^3+x_4^2=0$
 in ${\bf P}_{24}(1, 1, 2, 8, 12)$ or $x_0=x_1=0, x_2^7+x_3^3+x_4^2=0$ in
 ${\bf P}_{84}(1, 1, 12, 28, 42)$. 

The weighted projective subspace ${\bf P}_{12}(1, 4, 6)$ defined by
 $x_0=x_1=0$ in ${\bf P}_{24}(1, 1, 2, 8, 12)$ is the singular set of
 ${\bf P}_{24}(1, 1, 2, 8, 12)$.

  By blowing up along the subspace in ${\bf P}_{24}(1, 1, 2, 8, 12)$, one
 establishes a weighted projective 
 ${\bf P}_{12}(1, 2, 4, 6)$ fiber bundle structure over ${\bf P}^1$. The proper
 transformation under the toric blowing up of the weighted hyper-surface $f$
 is denoted by $Z_0$.  Similarly by blowing up along the subspace
 ${\bf P}_{42}(6, 14, 21)$, the resulting threefold constructed from
 resolving $g=x_0^{84}+x_1^{84}+x_2^7+x_3^3+x_4^2=0$ is denoted by $W_0$.

\begin{prop}\label{prop; smooth}
The algebraic threefold $Z_0, W_0$ are smooth and the restrictions of the
 toric fibrations to $Z_0, W_0$ induce $K3$ fibration structures
over ${\bf P}^1$.
The $K3$ fibration structure on $Z_0$ and (on $W_0$)
have $12$ and (respectively $42$) singular fibers, each of them isomorphic to
 a reduced singular $p_g=0$ surface with a unique hyperbolic singularity.
\end{prop} 

\noindent Proof: Both the weighted projective spaces 
 ${\bf P}_{84}(1, 1, 12, 28, 42)$ and ${\bf P}_{24}(1, 1, 2, 8, 12)$ are singular.
Because they are toric, they can be described by toric datum of fans. Denote
 $e_1=(1, 0, 0, 0)$, $e_2=(0, 1, 0, 0)$, $e_3=(0, 0, 1, 0)$ and
 $e_4=(0, 0, 0, 1)$. Denote $e_5=(-1, -12, -28, -42)$ and $\hat{e}_5=
(-1, -2, -8, -12)$. 

 Then 

$$e_1+12e_2+28e_3+42e_4+e_5=0, e_1+2e_2+8e_3+12e_4+\hat{e}_5=0$$

Then ${\bf P}_{84}(1, 1, 12, 28, 42)$ is
 determined by the fan formed by the $1$-edges generated by
$e_1, e_2, e_3, e_4$ and $e_5$. Likewise ${\bf P}_{24}(1, 1, 2, 8, 12)$ is
determined by the fan formed by the $1$ edges generated by 
$e_1, e_2, e_3, e_4$ and $\hat e_5$.

Recall the well known criterion of smoothness of the affine pieces,

\begin{lemm}
 The $n$ dimensional 
 affine toric variety determined by an integral cone is smooth if and only if

\noindent  (i). The cone is simplicial. i.e.\ it is generated by $n$ edges.

\noindent (ii). 
The primitive integral elements of the $n$ edges form a ${\bf Z}$ basis of
 ${\bf Z}^n$.  

\end{lemm}

 One can check without much difficulty that 
the singularity sets are the weighted projective subspaces
 ${\bf P}_{42}(6, 14, 21)$ and ${\bf P}_{12}(1, 4, 6)$ with the normal cones
 modeled on ${\bf Z}_2$ surface quotient \footnote{The ${\bf Z}_2$
 orbifold singularities appear here because the greatest common divisor of the 
tuple $(12, 28, 42)$ (or $(2, 8, 12)$) is $2$.} orbifold singularities.

 Since the hypersurfaces are both of Fermat types, it is
easy to check that the hypersurfaces are transversal in the weighted
 projective spaces. Their singularities are induced from the
 toric singularities of the ambient spaces. 

 Recall that a single ${\bf Z}_2$ orbifold quotient singularity can be
 resolved by a single blowing up at the singular point and the exceptional
 curve is a ${\bf P}^1$.
 Likewise the toric singularities are resolved by single blowing ups along
the codimension two singular loci in the ambient toric varieties, and the
exceptional divisors have the structures of ${\bf P}^1$ bundles.
 The resolution also give $Z_0$ and $W_0$
 $K3$ fibration structures over ${\bf P}^1$. We denote the projection maps
 to ${\bf P}^1$ by $\pi$.

 The generic fibers $\pi^{-1}(\lambda), \lambda \in {\bf C}\subset 
 {\bf P}^1$, $\lambda^{24}\not=-1$, or $\lambda^{84}\not=-1$ can be identified with
 the $\lambda-$ dependent 
hypersurfaces $$\lambda^{24}x_1^{24}+x_1^{24}+x_2^{12}+x_3^3+x_4^2=0,$$
 and $$\lambda^{84}x_1^{84}+x_1^{84}+x_2^7+x_3^3+x_4^2=0,$$

 which are isomorphic to
 $$(\lambda^{24}+1)y_1^{12}+x_2^{12}+x_3^3+x_4^2=0,$$

 and 
$$(\lambda^{84}+1)y_1^{42}+x_2^7+x_3^3+x_4^2=0$$
in  ${\bf P}_{12}(1, 1, 4, 6)$ and ${\bf P}_{42}(1, 6, 14, 21)$,
 after the change of variable $y_1=x_1^2$ and the corresponding weight.

For generic $\lambda$, 
these hypersurfaces are smooth Fermat $K3$ surfaces. The singularities
appear in the fiber $\pi^{-1}(\lambda)$ 
when $\lambda^{24}=-1$ in the first case 
or $\lambda^{84}=-1$ in the second case, and the hypersurfaces
 degenerate to $x_2^7+x_3^3+x_4^2=0$ and $x_2^{12}+x_3^3+x_4^2=0$ with
 weighted coordinates $(x_1, x_2, x_3, x_4)$. 
These hypersurfaces carry isolated singularities at
 $(1, 0, 0, 0)$ in the weighted projective
 spaces.

 According to Arnold's classification of isolated singularities, these
 singularities are hyperbolic singularities. These singularities are
 also closed related to the $K3$ singularities in the literature. Please 
consult the proof of the assertion (3). of proposition \ref{prop; kk}
 for their relationship with the $K3$ singularities.

\medskip

 The Milnor numbers $\mu$ can be calculated easily, and they are
 $12$ and $22$, respectively. One checks that 
different smooth fibers of each $K3$ fibration are bi-holomorphic to each other.
 The bi-holomorphic isomorphism is constructed by identifying
 $\pi^{-1}(\lambda)$ with $\pi^{-1}(0)$ through
$(y_1, x_2, x_3, x_4)\mapsto (({1\over 1+\lambda^{k}})^{2\over k}y_1, x_2, 
x_3, x_4)$, $k=24$ or $84$ and by composing above 
the identifications for different
 $\lambda$.

In particular, the one parameter family of bi-holomorphisms 
 induced by $\lambda_0+\epsilon e^{i\theta}$, varying $\theta$ in the range 
$0\leq \theta <2\pi$, induces 
a nontrivial monodromy action on the middle cohomology of fiber $K3$, where the
 classical monodromy operator is a root of unity.
$\Box$ 

 Such bi-holomorphic monodromy diffeomorphisms are
 $K3$ complex multiplication in the literature, which
do not preserve the holomorphic two forms on the smooth fibers.
 This can be shown by a direct computation on the action of holomorphic two
forms.

\medskip

   Knowing that the classical monodromy diffeomorphism
 is a complex multiplication,
 then the singular fiber must be of $p_g=0$. Namely, they must
 be singular rational surfaces.

On the other hand, the following proposition allows one to
 perturb the complex structure of the total spaces $Z_0$ or $W_0$ into
 tamed $K3$ fibrations such that
 their singular fibers carry only ordinary double points.

\begin{prop}\label{prop; kk} 
 Let $Z_0$ and $W_0$ be the smooth Fermat type hypersurfaces embedded in the 
toric fiber bundles defined as above. Then there exist arbitrarily small
 deformations of the defining equations of $Z_0$ and $W_0$
such that the zero loci of the defining equations satisfy the
 following conditions.

\noindent (1). The deformed hypersurface is smooth.

\medskip

\noindent (2). It still carry a $K3$ fibration structure induced by the
 toric fiber bundle structure of the ambient space.

\medskip

\noindent (3). All the singular fibers carry single 
isolated singularities which are ordinary double points.
\end{prop}

\noindent Proof: The assertion 
(1). follows from the fact that the smoothness of the hypersurface 
is an open condition
 within the space of polynomial perturbations of $Z_0$ or $W_0$.

The condition (2). follows from the fact that the total spaces of the toric
 varieties project naturally to ${\bf P}^1$. The induced fibration is
still a $K3$ fibration since the generic fibers of the deformed fibration
 are perturbations of the original $K3$ fibers.

The assertion (3). deserves some additional attention.  

 Let $f(x_2, x_3, x_4)=x_2^{12}+x_3^3+x_4^2$ or 
$x_2^7+x_3^3+x_4^2$. By
 setting $\epsilon=(1+\lambda^{k\over 2})$, with $k=24$ or $84$, then 
 the pencil of $K3$ can be expressed as 

 $$\epsilon y_1^{k\over 2}+f(x_2, x_3, x_4)=0,$$
 
 and the local $K3$ fibration structure is given by 
 $(y_1, x_2, x_3, x_4)\mapsto \epsilon$.

  To understand how to improve the singularity of the central fiber, we
study the geometric structure of the fibration and its relationship with
 the product family.

 Consider the $K3$, $M$,
 defined by the equation $y_1^{k\over 2}+f(x_2, x_3, x_4)=0$
 in the weighted projective space. Take $M\times {\bf C}$ with the
affine coordinate $a$ on ${\bf C}$. For $a\in {\bf C}-\{0\}$, 
 $(y_1, x_2, x_3, x_4)\mapsto (ay_1, x_2, x_3, x_4)$ is a one parameter
family of automorphisms. After substituting $y_1$ by $ay_1$, the
  equation becomes $a^{k\over 2}y_1^{k\over 2}+f(x_2, x_3, x_4)=0$. 

 For $a=0$, the central fiber algebraic surface 
defined by $f(x_2, x_3, x_4)=0$ is 
isomorphic to the projectified cone over the projective curve $\Sigma$ 
defined by 
$f(x_2, x_3, x_4)=0$ in the weighted projective subspace $y_1=0$. The curve is also
the fixed locus of the one parameter family of automorphisms.

 To get this new family \footnote{The total space of the new family
is singular.} from $M\times {\bf C}$, we blow up $M\times {\bf C}$
 along $\Sigma\times \{0\}$. The central fiber consists of two components,
 a ruled surface (which is the exceptional divisor of the blowing up) and
 $M$. After blowing down along $M$ in the central fiber, it creates a so-called
 $K3$ singularity and it collapses the ruled surface into a projectified cone
 over $\Sigma$.

 The $K3$ singularity is nothing but the singularity at the origin of the
 projectified cone over the curve $\Sigma$.
 
 To recover the smooth model we start with, we consider the ${\bf Z}_{k\over 2}$
 action generated by $(y_1, x_2, x_3, x_4, a)\mapsto (\eta y_1, x_2, x_3, x_4, 
\eta a)$, where $\eta$, $\eta^{k\over 2}=1$ is a primitive $k\over 2$-th
 root of unity\footnote{This is consistent with the action of the
 monodromy.}. The quotient variety by the cyclic ${\bf Z}_{k\over 2}$ action
 is nothing but the original family after we perform the substitution 
$a^{k\over 2}\mapsto \epsilon$. After the quotient, the isolated 
singularity at the center fiber becomes a hypersurface singularity 
in the orbifold ${\bf C}^3/{\bf Z}_{k\over 2}$.

Next recall the basic fact about ``morsification'' in singularity
theory [Arnold]. 

Let ${\bf f}^{-1}(0)$, with
${\bf f}:{\bf C}^n\mapsto {\bf C}$, be a germ \footnote{In general, the
 domain of ${\bf f}$ should be a open neighborhood of ${\bf 0}$ instead of
 ${\bf C}^n$. We use ${\bf C}^n$ here because the hyperbolic singularities
 of the singular $K3$ fibers are ${\bf f}^{-1}(0)$ for some
 ${\bf f}:{\bf C}^3\mapsto {\bf C}$.}
of an isolated singularity
 at ${\bf 0}\in {\bf C}^n$.  Then it is a well known fact
 in singularity theory that a small perturbation of ${\bf f}$ by
generic linear functionals in ${\bf C}^n$ have only non-degenerated 
singularities. Such non-degenerated singularities are the ordinary
 double points.

 Let us formulate it as a lemma,

\begin{lemm}\label{lemm; morse}
 Let ${\bf 0}\in {\bf f}^{-1}(0)$, with ${\bf f}:{\bf C}^n\mapsto {\bf C}$,
 be an arbitrary 
germ of an isolated hypersurface singularity
 at ${\bf 0}\in {\bf C}^n$ and let $g_i:{\bf C}^n\mapsto {\bf C}, 
1\leq i\leq n$ be the $i-$th coordinate functions. Then for generic
 choices of $\underline{\lambda}=(\lambda_1, \lambda_2, \cdots, \lambda_n)$,
 $${\bf f}_{\underline{\lambda}}={\bf f}+\underline{\lambda}{\bf g}=
  {\bf f}+\sum_i{\lambda}_i{ g}_i$$
induces a holomorphic map ${\bf f}_{\underline{\lambda}}:{\bf C}^n\mapsto
 {\bf C}$ whose critical points are non-degenerated.
 The $\underline{\lambda}$ tuple can be chosen such that each singular fiber
contains exactly one non-degenerated critical point.
\end{lemm}

 A proof of the lemma can be found in [Arnold] page 30-31.

In our situations, we consider the family 
$a^{k\over 2}y_1^{k\over 2}+f(x_2, x_3, x_4)=0$, before performing the
 ${\bf Z}_{k\over 2}$ quotient. 
we take $n=3$ and ${\bf f}(\hat{x}_2, \hat{x}_3, \hat{x}_3)
=\hat{x}_2^7+\hat{x}_3^3+\hat{x}_4^2$ or 
$\hat{x}_2^{12}+\hat{x}_3^3+\hat{x}_4^2$. The
 isolated singularity at the origin of ${\bf C}^3$ 
determined by ${\bf f}=0$ is the $K3$ singularity mentioned
 above.

 Then observe that by substituting $\hat{x}_i={x_i\over y_1^{w_i}}$, 
the linear perturbation
 by $\lambda_i g_i$ can be prolonged to a polynomial perturbation
 $\lambda_i y_1^{{k\over 2}-w_i}x_i$, with $k=24$ or $84$ for the equations 

$$a^{12}y_1^{12}+x_2^{12}+x_3^3+x_4^2=0,$$

 and 

$$a^{42}y_1^{42}+x_2^7+x_3^3+x_4^2=0.$$

Over here the symbol
$w_i$ means the weight of the variable $x_i$.

 We extend the ${\bf Z}_{k\over 2}$ action to the tuple $\underline{\lambda}$
 as well by sending $\lambda_i\mapsto \eta^{w_i}\lambda_i$. It is easy
 to see that the deformed equation is ${\bf Z}_{k\over 2}$-equivariant and
 can be descended into the ${\bf Z}_{k\over 2}$-quotient.

\bigskip

 Once the tuple $\underline{\lambda}$ is chosen to be small enough, all the
 ordinary double points morsified from the original
 singularity correspond to the values $a\not=0$. $\Box$

\medskip

\begin{rem}
 The above proposition can be viewed as an analogue of the
 well known fact that any regular elliptic surface without multiple-fibers 
can be deformed into
 a ``generic'' elliptic fibration with only type $I_1$ singular fibers (i.e.
 unique $A_1$ nodal singularities in the singular fibers).
\end{rem}

\medskip
 
  There are a whole class of Calabi-Yau $K3$ fibrations over ${\bf P}^1$ 
found by the string theorists and algebraic geometers [KV], [KLM], [B], etc.
  Start from the basic building blocks, one may construct new $K3$ 
fibrations by the fiber sum construction. They usually fall out of
 the Calabi-Yau category, but they are still interesting examples to 
study, as our theory indicates that the Calabi-Yau condition plays 
 a minor role in the enumerative theory of curves on $K3$ fibrations.

\subsection{The Kawamata-Viehweg Covering Trick and the $K3$ Fiber Bundle} 
\label{subsection; KV} 

 In studying the curve counting on the original $K3$ fibration, 
$X_0\mapsto B_0$, we recast it into a better model
 such that the new space has a structure of
 a relative algebraic $K3$ fiber bundle over its base.

Firstly, let us recall the well known trick in algebraic geometry (see e.g. lemma 
5 of [Kaw]),

\begin{lemm}(Covering trick) (Kawamata-Viehweg)\label{lemm; KV} 
 Let $N$ be a non-singular projective algebraic variety
 of dimension $k$. Let $\cup_{i\in I} D_i$ be a simple normal
 crossing divisor with smooth irreducible components. Suppose that
 $m_i$, $i\in I$, are a sequence of positive integers attached to $D_i$, $i\in I$,
 then there exists a non-singular projective algebraic variety
 $N'$ constructed as a branched
 covering of $N$ such that the covering is ramified at each of  $D_i$
with ramification multiplicity $m_i$.
\end{lemm} 

  The usage of the covering trick to the case of complex 
curves is particularly simple and its proof is straightforward.
 Let $B_0$ be the base curve of $X_0\mapsto B_0$ and let
 $p_i, 1\leq i\leq n$, be the singular values of ${\cal \pi}_0: X_0\mapsto B_0$. 
 We pick $N=B_0$ and $D_i=p_i, 1\leq i\leq n$. Mumford's semi-stable reduction
 theorem on algebraic surfaces implies the following:
Let $z=p_i$ be a singular value of $\pi_0:{\cal X}_0\mapsto B_0$.
  There exists a finite ramified covering of
a neighborhood of $z$, denoted by $g_z: 
\tilde{\cal N}_i\mapsto {\cal N}_i$, 
 such that $ \tilde{\cal N}_i\times_{{\cal N}_i} {\cal X}_0\mapsto {\cal N}_i$ is
 birational to a new local $K3$ fibration over 
 $\tilde{N}_i$ whose central singular fiber consists of 
simple normal crossing smooth divisors.

According to the birational 
classification result (e.g. [Mo]) of the semi-stable fibers of
 normal-crossing $K3$, a priori, there are three distinct cases to consider.

\medskip

\noindent (I). The type $I$ special fiber is a smooth K3 surface. 

\medskip

\noindent (II). The type $II$ special fiber. The irreducible components of
 the normal crossing divisor are all rational. The dual complex of the
 normal crossing divisor gives a triangulation of $S^2$.

\medskip

\noindent (III). The type $III$ special fiber. All the irreducible components are
 rational or ruled.

\medskip

When the global $K3$ fibration is tamed, we take the integers $m_i$ to
be the order of the classical local monodromy operator around $\pi_0^{-1}(p_i)$.

 Then the
 covering trick implies the existence of a finite branched covering
$f: B\mapsto B_0$ such that
 there exists a semi-stable $K3$ fibration over $B$ which is birational to
 the pulled back fibration ${\cal X}_0\times_{B_0}B\mapsto B$. 
By the choices of the
 multiplicities $m_i$, the local monodromy around the singular fibers
 are all trivial.

 Because the local monodromy operators around the type $II$ and type $III$
 fibers are known to be of infinite orders, 
 the semi-stable central fibers must be of type $I$, i.e. smooth $K3$ surfaces. 
In other words,
 one can blow down the exceptional divisors to make the new $K3$ fibration
 an algebraic $K3$ fiber bundle.

\subsection{The Iso-trivial Families of $K3$}\label{subsection; iso}

\bigskip

Recall that a holomorphic 
fibration $X_0\mapsto B_0$ is said to be an iso-trivial
 family if there exists a finite branched covering $B_0'\mapsto B_0$
 such that $B_0'\times_{B_0}{\cal X}_0$ is birational to 
the trivial product $B_0'\times X$.

 We have the following proposition characterizing the $K3$ fibrations
 appearing in section \ref{subsection; exam}. 

\begin{prop} \label{prop; iso}
 The Calabi-Yau $K3$ fibrations $Z_0$ and $W_0$ in section \ref{subsection; exam}
 are both iso-trivial $K3$ fibrations.
\end{prop} 

\noindent Proof: Because the isolated singularities of the singular fibers in
 $Z_0$ and $W_0$ are both Fermat type quasi-homogeneous 
hyperbolic singularities, their
 classical monodromies are finite \footnote{See section \ref{subsection; exam}
 for the
 identification of their monodromies.}. Therefore, both $Z_0\mapsto B_0$ and 
$W_0\mapsto B_0$ are tamed $K3$ fibrations.  By applying the covering trick, 
there exists finite ramified coverings such
 that the pulled-back fibrations are birational to relatively 
algebraic $K3$ fiber bundles.

 Because the regular fibers of $Z_0$ (and $W_0$) are all bi-holomorphic to each
 other and the set of regular values are open and dense in ${\bf P}^1$, the
 fibers of the algebraic $K3$ fiber bundle constructed from 
 $Z_0$ or $W_0$ are all bi-holomorphic to a fixed
  smooth fiber. These bi-holomorphisms enable us to 
 construct the required isomorphism from the product variety
  to the algebraic $K3$ fiber bundle. So the $K3$ fiber bundle is
 a trivial product. $\Box$

\bigskip

  The proposition enables us to re-construct $Z_0$ and $W_0$
 ${\bf birationally}$ by the following quotient process. 

\begin{defin}\label{defin; CM}
 Let $X$ be an algebraic $K3$ surface and let $\Omega_X$ be the canonical bundle.
 An automorphism $\sigma\in Aut(X)$ is said to be a complex multiplication
 if $\sigma$ induces a non-trivial character on the one dimensional vector 
space of holomorphic two forms $H^0(X, \Omega_X)$.
\end{defin}

 The concept of complex multiplications have played an important role
 in the theory of elliptic curves (one dimensional Calabi-Yau space).
 The reader may consult [Kon], [N1] for more details about complex
 multiplications on $K3$ surfaces.

\bigskip

 Let $\mu$ denote the order of $\sigma$ and let $\mu_0$ denote the
 order of the image of $\sigma$ in 
 $End(H^0(X, \omega_X))$. We consider the $\mu<\infty$ case only.

 Take $B_0={\bf P}^1$ and take $p_i, 1\leq i\leq n$, to be a finite number of
 points in $B_0$.  By applying lemma \ref{lemm; KV}, one can find a smooth
 algebraic curve $\Sigma\mapsto B_0$ such that the projection map ramifies
 along $p_i$ with multiplicity $\mu_0$. Consider $\sigma'$ to be
 the covering automorphism of $\Sigma\mapsto B_0$ with fixed points. 

 Take $X\times \Sigma$ and consider the product action of the cyclic group
 generated by $\sigma\times \sigma'$.  Then the
 quotient $X\times \Sigma/\langle \sigma \times \sigma'
\rangle$ is an algebraic threefold
 which allows a $K3$ fibration structure. 
 By construction, the regular fibers of the
$K3$ fibrations are bi-holomorphic to $X$. 

 The reader may consult Borcea's paper [B] for examples of Calabi-Yau K3 fibrations
 constructed from complex involutions on algebraic K3s. 

\medskip

\section{The Cosmic String and Family Seiberg-Witten Invariant}
\label{section; cosmic} 

\bigskip

In this section, we plan to discuss the relationship between the idea of
 cosmic string [GSYV] and its special 
role in the family Seiberg-Witten theory of K3 surfaces.

 In subsection \ref{subsection; general}, we deal with the general 
$K3$ fibrations whose singular fibers carry $A-D-E$ singularities.

\medskip

 In subsection \ref{subsection; pathetic}, we discuss briefly the problem of
 multiple-coverings of $-2$ curves along the $K3$ pencil.

\medskip

 In subsection \ref{subsection; CY}, we discuss how does the Calabi-Yau
 condition constrain the degree of the cosmic string. At the end of 
subsection \ref{subsection; CY}, we prove a defect relationship and 
explain how does the family invariant get ``sucked'' into the singular
 fibers of the iso-trivial $K3$ fibration when we degenerate from a
 ``generic'' $K3$ fibration to an iso-trivial one. This can be 
interpreted as a fractional bubbling-off phenomenon of the cosmic string map
 $\Phi_{\cal Y}:\tilde{B}\mapsto \underline{\cal M}_{\bf M}$.

\medskip 

 The concept of cosmic string will play an essential role in proving
 the main theorem of the paper.

\subsection{The Family Seiberg-Witten Invariants of $K3$ Fibrations}
\label{subsection; general}

In this subsection, we consider the $K3$ fibrations whose singular
 fibers contain only $A-D-E$ singularities. Our goal is to relate the
 Weil-Peterson 
symplectic volume (degree) of the cosmic brane and family Seiberg-Witten
 invariant of a given $K3$ fibration.

\begin{prop} \label{prop; wall} 
 Let $\pi:{\cal X}\mapsto B$ be a tamed 
$K3$ fibration from a smooth total space
 ${\cal X}$, $h^{(2, 0)}({\cal X})=0$, to a smooth 
base space $B$. Let $\omega_{{\cal X}/B}$ denote the relative
 polarization.
Suppose that the singular fibers (if there is any)
 contain only rational double points ($A-D-E$ surface singularities), then
 the algebraic family Seiberg-Witten invariants of a fiberwise class 
 $C\in H^2({\cal X}, {\bf Z})_f$, with $\int_{{\cal X}/B}
\omega_{{\cal X}/B}\cup C>0$, is well-defined and it can be
 expressed as 

 $$(-1)^{dim_{\bf C}B-1}
\int_B c_1^{dim_{\bf C}B}({\cal R}^2\pi_{\ast}{\cal O}_{{\cal X}})+....,$$
where the correction term $....$ in the formula represents a pairing
 $\int_B  c_1({\cal R}^2\pi_{\ast}{\cal O}_{{\cal X}})\cup U(C)$. The class 
$U(C)\in H^{2dim_{\bf C}B-2}(B, {\bf Z})$
 represents a universal polynomial expression (in terms of the cup product) of 
 the push-forward expressions of $\pi_{\ast}(Todd_{{\cal X}/B}C^k)$, $k\in
{\bf N}$ and $c_1({\cal R}^2\pi_{\ast}{\cal O}_{{\cal X}})$. When the fibration
${\cal X}\mapsto B$ is smooth, it is reduced to a polynomial expression of 
$$\pi_{\ast}(c_1^a({\bf T}({\cal X}/B))c_2^b({\bf T}({\cal X}/B))C^c), c>0$$
 and $c_1({\cal R}^2\pi_{\ast}{\cal O}_{{\cal X}})$.
\end{prop}

\noindent Proof of proposition \ref{prop; wall}: 

\noindent Step 1:
 We begin by recalling the algebraic family Kuranishi model
 construction of $C\in H^2({\cal X}, {\bf Z})_f$. For simplicity we assume
 that the first betti number of ${\cal X}$ is zero. By the 
assumption $h^{2, 0}({\cal X})=0$, the class $C$ determines 
 a holomorphic line bundle (invertible sheaf) on ${\cal X}$, denoted by 
${\cal E}_C$.

  By considering a very ample effective divisor $D$ which
 restricts to very ample divisors on the fibers and
 a large enough $n\in {\bf N}$, 
we consider the following derived long exact sequence,

$$\hskip -.3in
0\mapsto {\cal R}^0\pi_{\ast}\bigl({\cal E}_C\bigr)\mapsto {\cal R}^0\pi_{\ast}
\bigl({\cal O}(nD)\otimes {\cal E}_C\bigr)\mapsto 
{\cal R}^0\pi_{\ast}
\bigl({\cal O}_{nD}(nD)\otimes {\cal E}_C\bigr)\mapsto 
  {\cal R}^1\pi_{\ast}\bigl({\cal E}_C\bigr)\mapsto 0.$$

The two terms in the middle 
are locally free for $n\gg 0$ by Serre vanishing theorem
(see e.g. theorem 5.2 on page 228 of [Ha]).

The difference of their ranks can be calculated by surface Riemann-Roch formula
and it is equal to 

$$1-q+p_g+{C^2-c_1(K3)\cdot C\over 2}=1+1+{C^2\over 2}.$$

The condition $\int_{{\cal X}/B}\omega_{{\cal X}/B}\cup C>0$ on the
 relative degree and the fact that the singular fibers are $K3$ orbifolds, and 
the relative Serre duality implies that 
${\cal R}^2\pi_{\ast}\bigl({\cal E}_C\bigr)=0$.

 On the other hand, for non-iso-trivial fibrations the
 invertible sheaf ${\cal R}^2\pi_{\ast}\bigl(
{\cal O}_{\cal X}\bigr)$ is non-trivial. According to [Liu2], the formal
 base dimension $febd(C, {\cal X}/B)=0$ and the expected dimension of the 
algebraic family Seiberg-Witten invariant is ${C^2\over 2}+dim_{\bf C}B$.

  So we consider ${\cal V}={\cal R}^0\pi_{\ast}
\bigl({\cal O}(nD)\otimes {\cal E}_C\bigr)$ and
 ${\cal W}= {\cal R}^0\pi_{\ast}
\bigl({\cal O}_{nD}(nD)\otimes {\cal E}_C\bigr)$. Then
 the sheaf morphism 
${\cal V}\mapsto {\cal W}$ induces a bundle map $\Phi_{{\bf V}{\bf W}}:
{\bf V}\mapsto {\bf W}$ and
its projectified kernel cone can be identified with the algebraic family
 moduli space of $C$, ${\cal M}_C$. It can be viewed as the zero locus of the
 canonical section of $\pi_{{\bf P}({\bf V})}^{\ast}{\bf W}\otimes {\bf H}$
 over ${\bf P}({\bf V})$ determined by $\Phi_{{\bf V}{\bf W}}$..

  Because the obstruction sheaf of the universal curve 
 ${\cal C}\mapsto {\cal M}_C$,  
 ${\cal R}^1\pi_{\ast}\bigl({\cal O}_{\cal C}({\cal C})\bigr)$ fits into
the following exact sequence,

$${\cal R}^1\pi_{\ast}\bigl({\cal O}_{{\cal X}\times_B{\cal M}_C}
({\cal C})\bigr)\otimes {\cal H}\mapsto 
{\cal R}^1\pi_{\ast}\bigl({\cal O}_{\cal C}({\cal C})\bigr)\otimes {\cal H}\mapsto 
{\cal R}^2\pi_{\ast}\bigl({\cal O}_{\cal X}\bigr)\mapsto 0,$$

 and because the surjectivity of

$$\pi_{{\bf P}({\bf V})}^{\ast}{\cal W}\otimes {\cal H}\mapsto 
{\cal R}^1\pi_{\ast}\bigl({\cal O}_{{\cal X}\times_B{\cal M}_C}
({\cal C})\bigr)\otimes {\cal H},$$
 we may define ${\cal AFSW}_{{\cal X}\mapsto B}(1, C)$ for 
 the fibration ${\cal X}\mapsto B$ to be
 $$\int_{{\bf P}({\bf V})}
c_{top}(\pi_{{\bf P}({\bf V})}^{\ast}{\bf W}\otimes {\bf H}\oplus 
 {\bf R}^2\pi_{\ast}({\cal O}_{\cal X}))\cup
 c_1^{{C^2\over 2}+dim_{\bf C}B}({\bf H}).$$

\medskip

 By a standard computation using the definition of Segre classes, the above
expression can be reduced to 
 $$\int_Bc_1({\bf R}^2\pi_{\ast}\bigl({\cal O}_{\cal X}\bigr))\cup
c_{total}({\bf W})\cup s_{total}({\bf V})=
\int_Bc_1({\bf R}^2\pi_{\ast}\bigl({\cal O}_{\cal X}\bigr))\cup 
c_{dim_{\bf C}B-1}({\bf W}-{\bf V}).$$

\bigskip

\noindent Step 2: We evaluate the above expression by realizing
 that 
$${\cal V}-{\cal W}={\cal R}^0\pi_{\ast}\bigl({\cal E}_C\bigr)-
{\cal R}^1\pi_{\ast}\bigl({\cal E}_C\bigr)=\pi_{\ast}{\cal E}_C$$
 in the K group of $B$ and it can be determined by the Grothendieck Riemann-Roch
 formula (family index formula in differential topology),

$$ch(\pi_{\ast}{\cal E}_C)=\int_{{\cal X}/B}Todd_{{\cal X}/B}ch({\cal E}_C).$$

We have $ch({\cal E}_C)=e^C$.
\medskip

 We separate the push-forward on the right hand side
 into two groups, one part containing push-forwards of powers $C^c, c>0$, 
the other part collecting all the terms which are $C$ independent.

As the most relevant case is when $\pi:{\cal X}\mapsto B$ is smooth, in the 
following we deal with the case when $\pi:{\cal X}\mapsto B$ is smooth. 
It is straight-forward to deal with the more general case.

 We would like to prove that the class 
$U(C)\in H^{2dim_{\bf C}B-2}(B, {\bf Z})$ is a 
universal polynomial expression involving the variables 
 $x_{a, b, c}=\pi_{\ast}(c_1^a({\bf T}({\cal X}/B))c_2^b({\bf T}({\cal X}/B))C^c)$,
 and $c_1({\cal R}^1\pi_{\ast}{\cal O}_{\cal X})$ in $H^{\ast}(B, {\bf Z})$. 

\medskip

 Firstly, we consider the relative Todd class of the vertical tangent bundle
 ${\bf T}({\cal X}/B)$ in $H^{\ast}({\cal X}, {\bf Z})$.

 $$Todd_{{\cal X}/B}=\sum_{i\geq 0} Todd_i(c_1, c_2),$$ 
where $Todd_i(c_1, c_2)$ is the degree $2i$ term of the relative Todd class 
and is a universal quasi-homogeneous polynomial of 
the variables $c_1({\bf T}({\cal X}/B)), c_2({\bf T}({\cal X}/B))$ with 
weight $1$ and $2$, respectively.

  These universal polynomials are determined implicitly by using the
 generating function of Bernoulli numbers 

   $${1\over e^t-1}=\sum_{i\geq 0} B_i {t^i\over i!},$$

 and the recursive formula

  $$x^{i+1}+y^{i+1}=(x+y)(x^i+y^i)-xy(x^{i-1}+y^{i-1})$$

with the convention $x+y\longrightarrow c_1$ and $xy\longrightarrow c_2$.

  By applying the Grothendieck Riemann-Roch formula 

$$ch(\pi_{\ast}{\cal E}_C)=\pi_{\ast}\bigl(Todd_{{\cal X}/B}ch({\cal E}_C)
\bigr),$$ we find that for 
$j\in {\bf N}$, the $j-$th term of the chern character of the direct
 image is equal to

  $$\sum_{m+n=j+2} {1\over n!}\pi_{\ast}\bigl(Todd_m(c_1, c_2)
c_1^n({\cal E}_C)\bigr)\in H^{2j}(B, {\bf Z}),$$

which (after using $c_1({\cal E}_C)=C$) 
is a linear combination of the push-forward, 

$$\pi_{\ast}\bigl( c_1^a({\bf T}({\cal X}/B))c_2^b({\bf T}({\cal X}/B))C^c 
\bigr), a+2b+c=j+2.$$

 Let us recall the proposition outlined on page 805 of [LL1], 
regarding the calculation of the total Chern (Segre) classes from the 
Chern characters of the index virtual bundle.

\begin{prop} \label{prop; chern-segre}
 Let $S(t)=\sum c_i t^i$ denote the Segre polynomial of the 
index virtual bundle\footnote{the direct image sheaf in the algebraic category.},
 then it can be expressed by the Chern characters through the following formula,

 $$S(t)=Exp\bigl(\sum_i (-1)^i (i-1)!ch_i\bigr)=
Exp\bigl(-ch_1 t+ch_2t^2-2!ch_3t^3+\cdots\bigr).$$
\end{prop}

Because all the $ch_i\in H^{even}(B, {\bf Z})$ are commutative to each other, the
 relative 
orders of different $ch_i$ do not matter. According to our discussion, each $ch_i$
 can be written as sums of $\pi_{\ast}\bigl( c_1^a({\bf T}({\cal X}/B))
c_2^b({\bf T}({\cal X}/B))C^c \bigr)$.

We may separate the sum
 $\sum_i(-1)^i (i-1)! ch_i$ into two different 
sums, the former collects all the terms 
 $\pi_{\ast}\bigl(c_1^ac_2^b\bigr)$ without
 $C^c$, and latter collects all the terms of the form
 $\pi_{\ast}(c_1^ac_2^bC^c), c>0$.
 
$$\sum_i(-1)^i (i-1)!ch_i=\sum_i(-1)^i(i-1)!ch_i'+\sum_i(-1)^i(i-1)!ch_i''.$$

By proposition \ref{prop; chern-segre}, the Segre polynomial can be expressed as
 the product

    $$S(t)=Exp\bigl(\sum_i(-1)^i(i-1)! ch_i'\bigr)
\cdot Exp\bigl(\sum_i(-1)^i(i-1)! ch_i'' t^i \bigr)$$

  According to the calculation in step $1$, the first factor on the right hand side is the $C$ independent term, which is 
nothing but the Segre polynomial of $\pi_{\ast}{\cal O}_{\cal X}$, equal to 
the Segre polynomial of ${\cal R}^2\pi_{\ast}\bigl({\cal O}_{\cal X}\bigr)$.

 Thus the Segre polynomial of the direct image is equal to

$$\bigl(\sum s_i({\cal R}^2\pi_{\ast}\bigl({\cal O}_{\cal X}\bigr))t^i\bigr)\cdot 
Exp\bigl(\sum_{i>0}(-1)^i (i-1)!ch_i''t^i\bigr).$$

   By expanding
 out the degree $dim_{\bf C}B-1$ term of $S(t)$, the family invariant is equal to 

$$(-1)^{dim_{\bf C}B-1}
\int_{B}c_1^{dim_{\bf C}B}({\cal R}^2\pi_{\ast}\bigl({\cal O}_{\cal X}\bigr))+ corrections,$$
 
and the 'correction terms' is a universal polynomial of 
$c_1({\cal R}^2\pi_{\ast}\bigl({\cal O}_{\cal X}\bigr))$ and the various $ch_i^{''}, i>0$, which can be expressed
 as universal polynomials of formal variables $\pi_{\ast}
\bigl(c_1^ac_2^bC^c\bigr)$.

 This finishes the proof of the proposition. $\Box$

\subsubsection{The Hodge Bundle over the Moduli Space of Marked $K3$ Surfaces}
\label{subsubsection; hodge} 

 As before let ${\bf M}\subset {\bf L}$ be a sub-lattice of the $K3$ lattice
 ${\bf L}$.
   Consider the moduli space of ${\bf M}-$marked algebraic
 K3 surfaces ${\cal M}_{\bf M}$. It is known [D]
that
 ${\cal M}_{\bf M}$ is a quotient of a complex hyperbolic
 space ${\bf P}$. The space ${\bf P}$ can be viewed as a
 subspace of the period domain and parametrizes the
 deformations of Hodge structures of $H^{2, 0}(K3, {\bf C})\subset
 H^2(K3, {\bf C})$.
The universal line bundle ${\bf U}$ on ${\bf P}$ of the lines
 $H^{2, 0}(K3, {\bf C})$ descends to a line
 bundle over ${\cal M}_{\bf M}$. Given the family of
 complex deformations of a fixed underlying smooth manifold
 $M$ diffeomorphic to $K3$, the local family of holomorphic
 two forms ${\bf \Omega}$ on $M$ defines a local holomorphic
 trivialization of ${\bf U}$.
 
 The first Chern form $\partial\bar{\partial}ln(\int_M
{\bf \Omega}\wedge \bar{\bf \Omega})$ of ${\bf U}$ is invariant
 under ${\bf \Omega}\mapsto e^f{\bf \Omega}$ for holomorphic functions 
 $f$ on ${\bf P}$ and defines
 a Kahler form on ${\bf P}$. Its descend onto the quotient ${\cal M}_{\bf M}$
 is known to be the Weil-Peterson form $\varpi_{wp}$ [T], [To].

\medskip

  From our discussion we have $c_1({\bf U})=[\varpi_{wp}]$.
 On the one hand, a $K3$ fiber bundle $\pi:{\cal Y}\mapsto 
B$ determines a cosmic brane map $\Phi_{\cal Y}:B\mapsto 
{\cal M}_{\bf M}$. One may pull back the line bundle
 ${\bf U}$ and the Weil-Peterson form $\varpi_{wp}$ to
 $B$.

 On the other hand, the second derived image bundle
 ${\bf R}^2\pi_{\ast}\bigl({\cal O}_{\cal Y}\bigr)$
 is isomorphic to ${\bf R}^0\pi_{\ast}\bigl(\Omega^2_{{\cal Y}/B}\bigr)
\cong \Phi_{\cal Y}^{\ast}{\bf U}$ by relative Serre duality.

 Therefore $c_1({\bf R}^2\pi_{\ast}\bigl({\cal O}_{\cal Y}\bigr))$
 is equal to $-c_1({\bf U})=-\Phi_{\cal Y}^{\ast}[\varpi_{wp}]$.
 
 \medskip

 By combining this observation with the conclusion of 
proposition \ref{prop; wall}, we find that,

\begin{lemm}\label{lemm; brane}
 Let ${\cal Y}\mapsto \tilde{B}$ be 
an algebraic fiber bundle of algebraic $K3$ surfaces
 and $C\in H^2({\cal Y}, {\bf Z})_f$, ${C^2\over 2}+dim_{\bf C}\tilde{B}\geq 0$. 
Then the algebraic family invariant ${\cal AFSW}_{
{\cal Y}\mapsto \tilde{B}}(1, C)$
  is equal to $-\int_{\tilde{B}} \Phi_{\cal Y}^{\ast}\varpi_{wp}^{top}+
corrections$. 
\end{lemm}

 Notice that up to a sign the leading term is nothing but the
symplectic volume (degree) of the cosmic brane image.

\begin{rem}
When $dim_{\bf C}\tilde{B}=1$, it is easy to observe from the
 proof of proposition \ref{prop; wall} that family invariant
 ${\cal AFSW}_{{\cal Y}\mapsto \tilde{B}}(1, C)$ is reduced to
 $-\int_{\tilde{B}}\Phi_{\cal Y}^{\ast}\varpi_{wp}$.

\medskip
 
 The fact that when $dim_{\bf C}\tilde{B}$ 
the above answer is independent to $C$ is crucial for the
 simplicity of the conclusion of the main theorem.
\end{rem}

\bigskip

\subsection{The Pathetic Symptom of Multiple Coverings of $-2$ Curves}
\label{subsection; pathetic}
 
\medskip

 In this subsection, we reflect why the algebraic enumeration of curves in the
 family ${\cal X}\mapsto B$ has to be interpreted in the sense of 
virtual numbers. For simplicity, we may assume $dim_{\bf C}B=1$.

   One difficulty in dealing with curve enumeration of ${\cal X}\mapsto B$
 is related to the fact that $dim_{\bf C}H^0({\cal X}_b, {\cal E}_C|_b)$ may 
jump when the point $b\in B$ moves along the base $B$.
 Their dimensions are usually different from the expected dimension 
${C\cdot C\over 2}+1$.

\medskip

\noindent 
(i). There may exist infinite many $-2$ classes in $H^2({\cal X}, {\bf Z})_f$ 
which are effective on the different fibers
 ${\cal X}_b, b\in B$.

\medskip

\noindent 
(ii). As there is no a priori control on the Picard numbers of the fibers over 
$b\in B$, ${\cal X}_b$, it may happen that some
 $-2$ class $\in H^2(X_b, {\bf Z})$ suddenly becomes effective when $b$ 
is specialized to special points in $B$.

\medskip

Because of the appearance of the $-2$ classes by either $(i)$ and $(ii)$, 
some effective curve dual to $C$ in a given fiber ${\cal X}_b$ 
may break off multiples of different $-2$ curves. In symplectic topology, this phenomenon is known to be
 the "bubbling off".  These "bubbling off" occur at the "boundary points" of
 the moduli space of pseudo-holomorphic maps.

  The decomposition of the curves into different irreducible components 
yields a corresponding decomposition of the class $C$ as $C=C'+\sum m_i D_i$,
 where $D_i\cdot D_i=-2$, $m_i\geq 1$  are the irreducible smooth
 $-2$ classes and $C'$ is dual to the sum of the classes of the 
components (counted 
with multiplicities) which are not $-2$ curve classes. By utilizing the
fact that relative canonical bundle is trivial on the classes $C$, $C'$ and
 $D_i$, etc., the expected family dimension of the class $C$ and
 $C'$ are related by

 $${C\cdot C\over 2}+dim_{\bf C}B={C'\cdot C'\over 2}+dim_{\bf C}B
 -\sum_i m_i^2+\sum_{i<j} m_i\cdot m_j D_i\cdot D_j.$$

The $-$ sign in front of the term $\sum_i m_i^2$ is due to the negativity
 of the self intersection number $D_i\cdot D_i$ and can frequently cause the
 expect family dimension of $C'$ to exceed that of $C$.
\medskip

\noindent (A). At a given generic 
fiber ${\cal X}_b, b\in B$, the generic members of 
 the the linear system $|{\cal E}_C|$ may not be smooth irreducible. 
 If $C'\cdot D_i\not=0$, this is can interpreted as a bubbling off
 phenomenon.
Thus, generic members in $|{\cal E}_C|$
 may bubble off different multiplicities of multiple coverings of $-2$ curves.
\medskip

\noindent (B). At a special point $b$, some new generically non-effective curve 
classes (which may not be monodromy invariant) may become effective in the
 fiber $X_b$, which gives additional $-2$ curve
classes that a given monodromy invariant $C$ may bubble off. 
\medskip

  If the phenomenon in (A). is fiber-independent and might be hopeful to be
 analyzed by a case by case study, the spontaneous appearance of 
 the new $-2$ curve classes is hard to analyze directly.

\medskip

  Therefore it is vital to realize that the curves involving 
these $-2$ curves do no contribute
 to the family invariant at all when $dim_{\bf C}B=1$. We apply the 
 analysis of type $II$ exceptional classes in [Liu6] to the current context. 
 When non-iso-trivial fibration with 
$dim_{\bf C}B=1$, the only exceptional classes with non-negative
 family dimension are $-2$ classes.
The excess contribution \footnote{As $febd(C, {\cal X}/B)=0$, there are no
 $c_1({\cal R}^2\pi_{\ast}{\cal X})$ factors.}
from a collection of exceptional classes
 $D_1, \cdots, D_k$ is ${\cal AFSW}_{{\cal X}\mapsto B}(
\times_B^{i\leq k}\pi_{i\ast}
[\times_B^{i\leq k}{\cal M}_{D_i}]_{vir}\cap \tau, C-\sum_{i\leq k}D_i)$, for some
 suitable $\tau$ class.

 Given a $-2$ class, the expected dimension of its family moduli space
 is $0$. So $\pi_{i\ast}
[{\cal M}_{D_i}]_{vir}$ determines a zero degree cycle class
 in ${\cal A}_0(B)$. When $k>1$, the 
$\times_B^{i\leq k}\pi_{i\ast}
[\times_B^{1\leq i\leq k}{\cal M}_{D_i}]_{vir}$ is in 
${\cal A}_{<0}(B)=\emptyset$.
    In any case the above mixed invariant is zero. When $k=1$, it is 
proportional to
 $\int_{{\bf P}({\bf V})_b}c_{top}(\pi_{{\bf P}({\bf V})}^{\ast}{\bf W}\otimes
 {\bf H})\cup c_1({\cal R}^2\pi_{\ast}{\cal O}_{\cal X})\cup 
 c_1^{{C'^2\over 2}+1}({\bf H})=0$, because 
$c_1({\cal R}^2\pi_{\ast}{\cal O}_{\cal X})|_b=0$.

 So the curves involving (multiple coverings of) $-2$ curves do not contribute
to the family invariant of $C$ at all and the algebraic family invariant 
${\cal AFSW}_{{\cal X}\mapsto B}(1, C)$ 
counts smooth curves dual to $C$.

\medskip

\begin{rem}
In section \ref{section; proof}, we will deal with the virtual numbers of
 nodal curves in the $K3$ fibrations. We adopt a similar argument
 based on the theory of excess contributions of 
type $II$ exceptional classes [Liu6] to derive the
 vanishing result of excess contributions from 
curves involving type $II$ exceptional curves. 
\end{rem}

\subsection{The Calabi-Yau condition and the Family Invariant}
\label{subsection; CY}

\bigskip

Earlier we have identified
 the algebraic family Seiberg-Witten invariants of complex three 
dimensional $K3$ fiber bundles. In this subsection, we address the 
situation that ${\cal X}_0\mapsto B_0$ is a Calabi-Yau K3
 fibration possibly with singular fibers.

 The following proposition is the main result of the subsection.

\begin{prop}\label{prop; minus 2}
 Let $X_0\mapsto B_0$ be a Calabi-Yau algebraic $K3$ fibration such that all
the singular fibers have at most $A, D, E$ type orbifold singularities.
 Then $$\int_{B_0}\varpi_{wp}^{dim_{\bf C}B_0}=c_1(B_0)^{dim_{\bf C}B_0}[B_0].$$
\end{prop}

\medskip

\noindent Proof of proposition \ref{prop; minus 2}: 
By taking the proper push-forward of the structure sheaf
${\cal O}_{X_0}$ along $\pi_0:X_0\mapsto B_0$, 
its direct image is 

$$\pi_{0\ast}{\cal O}_{X_0}={\cal R}^0\pi_{0\ast}({\cal O}_{X_0})+
{\cal R}^2\pi_{0\ast}({\cal O}_{X_0}).$$

 The first term is isomorphic to ${\cal O}_{B_0}$ and the second
 term is invertible, by the base change theorem and the fact 
that $H^2(\pi_0^{-1}(b'), {\cal O}_{\pi_0^{-1}(b')})\cong {\bf C}$ for 
 fiber (orbifold) $K3$.  

\medskip

Grothendieck Riemann-Roch theorem enables us to calculate its Chern character and
 to identify it with $$\int_{X_0/B_0}Todd_{X_0/B_0}
=\int_{X_0/B_0}
 Todd_{X_0}/Todd_{B_0}$$
 $$=\int_{X_0/B_0}(1+{c_1(X_0)
\over 2}+{c_1(X_0)^2+c_2(X_0)\over 
 12}+\cdots)/(1-{c_1(B_0)\over 2}+\cdots).$$

 Because the sheaf 
${\cal R}^2(\pi_0)_{\ast}({\cal O}_{X_0})$ 
is invertible, its first Chern class determines its Chern character
 completely,

$$c_1({\cal R}^2\pi_{0\ast}({\cal O}_{X_0}))
=-{1\over 24}\int_{X_0/B_0}c_2(X_0)
\cup (c_1(B_0)-c_1(X_0))$$

$$=-c_1(B_0)+c_1(X_0)=-c_1(B_0),$$

 by using the Calabi-Yau condition $c_1(X_0)=0$.

 When the $K3$ fibration contains only $A$, $D$, $E$ type singular
fibers, there is a well defined cosmic string map 
$\Phi_{X_0}$ from $B_0$ to the appropriated moduli space
 of marked $K3$ surfaces.
 Then $c_1({\cal R}^2({\bf \pi}_0)_{\ast}({\cal O}_{X_0}))$ 
is nothing but  
$c_1(\Phi_{X_0}^{\ast}{\bf U}^{\ast})=-[\varpi_{wp}]$. 

 Thus $$\int_{B_0}\varpi_{wp}^{dim_{\bf C}B_0}
=c_1(B_0)^{dim_{\bf C}B_0}[B_0].$$  $\Box$ 
\medskip

\begin{cor}  \label{corr; P}
 If $B_0={\bf P}^1$, then the integral pairing $\int_{B_0}\varpi_{wp}$
 is reduced to $2$.
\end{cor} 

  In fact, this is the only numerical restriction of the Calabi-Yau condition upon
the invariant computation of $X_0\mapsto B_0$.  If the singular fibers of 
$X_0\mapsto B_0$ contain singularities other than rational
 double points, then the above 
assertion does not hold.  Under some additional assumption, we derive 
a defect formula to capture the discrepancy, 
which is applicable to the examples $Z_0$, $W_0$ in
 subsection \ref{subsection; exam}. 

\medskip

For simplicity, let us assume $dim_{\bf C}B_0=1$ without going 
into much technical details. 

 Let $\cup_i p_i$ be the singular values of the Calabi-Yau $K3$ 
fibration $\pi_0:X_0\mapsto B_0$.

  We assume the existence of a branched covering 
$h:\tilde{B}_0\mapsto B_0$ branched along the singular values 
$\cup_i p_i$, such that the
 pulled-back fibration $X_0\times_{B_0}\tilde{B}_0$ is birational
 to a smooth $K3$ fiber bundle $\tilde{X}_0\mapsto \tilde{B}_0$.

 Then there exists a well defined "cosmic string" map
 $\Phi_{\tilde{X}_0}$ from $\tilde{B}_0$ to the appropriated 
 ${\bf M}-$marked moduli space of K3 surfaces.

\medskip 

 Each of the closed loop $\gamma(i)$ surrounding $p_i$ induces a
 finite order automorphism $\in Aut(\tilde{B}_0)$ with fixed points
 and a monodromy action on the middle cohomology of the fiber $K3$. 
Its induced actions on the complex line of holomorphic two forms 
$\cong H^0(\tilde{\pi}_0^{-1}(b'), \Omega^2)$ ($b'$ fixed by the
 automorphism of $\tilde{B}_0$) must
be a complex root of unity $\lambda_i$.
The number $\lambda_i=1$ if the singular fiber 
of $X_0\mapsto B_0$ above $p_i$ contain only rational
 double points.

\begin{defin}\label{defin; defect}
 Define the local defect by 
 $$\delta(p_i)={arg(\lambda_i)\over 2\pi}, 0\leq arg(\lambda_i)<2\pi.$$
\end{defin}

 If the singularities of the singular fiber above $p_i$ 
are all isolated rational double points, then $\delta(p_i)=0$.

 We have the following proposition regarding the local defects,

\begin{prop}\label{prop; defect}
Let $\pi_0:X_0\mapsto B_0$
 be a Calabi-Yau $K3$ fibration over a smooth curve $B_0$.
 Suppose that 

\medskip

(i). The fibration $X_0\mapsto B_0$ can be deformed into a $K3$ fibration whose
 singular fibers contain only $A$, $D$, or $E$ types of singularities,

\medskip

(ii). there exists a smooth
 $\tilde{\pi}_0: \tilde{X}_0\mapsto \tilde{B}_0$ birational to the
 $K3$ fibration pulled back by the branched covering $h:\tilde{B}_0\mapsto B_0$,

then the local defects defined in definition \ref{defin; defect}, 
Weil-Peterson degree (symplectic volume) of 
$B_0$ and the first Chern class of
 $B_0$ are related by,

$$\int_{B_0}\varpi_{wp}+\sum_i \delta(p_i)=c_1(B_0)[B_0].$$
\end{prop}

The above formula imposes a constraint on the integrability of 
the sum of rational numbers $\sum_i \delta(p_i)$.

\medskip

\noindent Proof of proposition \ref{prop; defect}: 
The symbols $\int_{B_0}\varpi_{wp}$ deserves some explanation.

Firstly we
explain why the Weil-Peterson volume on $B_0$
 is still finite.
The Weil-Peterson form on $\tilde{B}_0$ is 
pulled back by the cosmic string map $\Phi_{\tilde{X}_0}$
induced by the $K3$ fiber bundle structure
 $\tilde{X}_0\mapsto \tilde{B}_0$. 
The Weil-Peterson form $\varpi_{wp}$ 
on the punctured curve $B_0-\cup_i p_i$ is pulled back
 by a (punctured) cosmic string map induced by the
 $K3$ fiber bundle $X_0\times_{B_0}(B_0-\cup_i p_i)\mapsto B_0-\cup_i p_i$.
 Its $h^{\ast}-$pull-back to $\tilde{B}_0-h^{-1}(\cup_i p_i)$ matches up with the
Weil-Peterson form on $\tilde{B}_0$.
 Because the Weil-Peterson volume of $\tilde{B}_0$ is
 finite, so is $\int_{B_0-\cup_i p_i}\varpi_{wp}$. This
 is what we mean by $\int_{B_0}\varpi_{wp}$.

\medskip

Put a Riemannian metric on $\tilde{B}_0$ invariant under the
 Decke transformations of the ramified covering $h:\tilde{B}_0\mapsto B_0$ 
and consider a small disc $D$ around
 a given pre-image point $p\in h^{-1}(p_i)$. The image $h(D-h^{-1}(p_i))$ 
is a punctured
 disc in $B_0$. Let $z_i$ denote a local holomorphic coordinate around
 $p_i\in B_0$.
Because $\tilde{X}_0\mapsto \tilde{B}_0$ is a $K3$ fiber bundle, the
restriction $\tilde{X}_0\times_{\tilde{B}_0}D\mapsto D$ can be trivialized
 to a product fiber bundle $M\times D\mapsto D$ differentially (where
 $M$ is diffeomorphic to $K3$ four-manifold).
 Under such a trivialization, the holomorphic 
family of fiberwise holomorphic two forms along
 $\tilde{X}_0\times_{\tilde{B}_0}D\mapsto D$ can be descended to
 an $h(D)$-family of fiberwise holomorphic two forms, denoted by
 ${\bf \Omega}(z_i)$.

 On the other hand, $X_0\times_{B_0}h(D-p)\mapsto h(D-h^{-1}(p_i))$ is a family
 of smooth $K3$s over the punctured disc $h(D-h^{-1}(p_i))=h(D)-p_i$. 
 The fiber bundle $X_0-\cup_i\pi_0^{-1}(p_i)\mapsto B_0-\cup_i p_i$ induces
 the (punctured)
 cosmic string map from the punctured curve $B_0-\cup_i p_i$. Because
 of the existence of the birational model, 
$\tilde{X}_0\times_{\tilde{B}_0}D\mapsto D$, the cosmic string map can
 be extended locally from $h(D)-p_i$ onto $h(D)$, and therefore from 
$B_0-\cup_i p_i$ onto $B_0$. We denote the resulting map $\Phi_{X_0}$.

In terms of
 ${\bf \Omega}(z_i)$, the monodromy non-invariant 
$h(D)-p_i$-family of fiberwise 
holomorphic two forms must be
 of the form $z_i^{\delta(p_i)}e^{f(z_i)}{\bf \Omega}(z_i)$ for some
 holomorphic function $f(z_i)$ of $z_i$, with
 the correct monodromy behavior around the 
 small loop $\gamma(i)$ surrounding $p_i$.

\medskip

 This implies that the distribution-valued first Chern form of the line bundle
$\Phi_{X_0}^{\ast}{\bf U}$ is 

 $$\partial\bar{\partial}ln(|z_i|^{2\delta(p_i)}|e^{f(z_i)}|^2
\int_M{\bf \Omega(z_i)}
\cup \overline{\bf \Omega}(z_i))=\varpi_{wp}+\lambda(p_i)\delta(p_i),$$

 over $h(D)$. Apply the same argument to all $p_i$.
Thus the distribution valued form
 is equal to $\varpi_{wp}+\sum_i\lambda(p_i)\delta(p_i)$ on $B_0$.

\medskip

 By the assumption (i) of our proposition,
 the deformation invariance of the first Chern 
number, and the proposition \ref{prop; minus 2}, $\int_{B_0}c_1({\bf U})
=c_1(B_0)[B_0]$ for $K3$ fibrations with rational double points along the
singular fibers, then

$$\int_{B_0}(\varpi_{wp}+\sum_i \lambda(p_i)\delta(p_i))=
\int_{B_0}\varpi_{wp}+\sum_i\lambda(p_i)=c_1(B_0)[B_0].$$ 

 So the proposition is proved. $\Box$

\medskip

 The above identity has some novel implication on the 
``bubbling off'' phenomenon of cosmic strings.
 The singular fibers of the Calabi-Yau $K3$ fibration $Z_0$ or $W_0$
in section \ref{subsection; exam} has isolated hyperbolic singularities.
 As these $K3$ fibrations are iso-trivial, $\int_{B_0}\varpi_{wp}=0$. Yet
 $c_1(B_0)[B_0]=2$. Thus, the above equality implies
 that the contribution merely comes from the summation of the local defects.

On the other hand, by proposition \ref{prop; kk} 
these $K3$ fibrations can be degenerated from
some tamed $K3$ fibrations with only ordinary double points among their 
singular fibers. By proposition \ref{prop; minus 2}, $\int_{B_0}\varpi_{wp}=2$
 for these ``generic'' $K3$ fibrations.

 Thus during the degeneration process from a ``generic''
 $K3$ fibration with $A_1$ singularities along the singular fibers 
to an iso-trivial $K3$ fibration,
 the harmonic energy of 
the cosmic string suddenly drops to zero and
 concentrates to the singular values of $\pi_0:X_0\mapsto B_0$, 
similar to the bubbling off phenomenon in the context of harmonic maps or
 $J$-holomorphic maps.

 The major difference from the usual bubbling off 
is that the energy lost at each $p_i$
 is generally fractional $\in {\bf Q}$.

 Consider the examples of iso-trivial families of $K3$ surfaces 
$Z_0$ and $W_0$ in section \ref{subsection; exam},
 the monodromy actions are complex multiplications of the 
$K3$ surfaces.

 The orders of the monodromies on the holomorphic two forms are of order
$12$ and $42$ respectively. 

 On the one hand, an easy calculation shows 
that there are $24$ and
$84$ identical singular fibers (which correspond to the number of
 the roots of $\lambda^k+1=0$, with
$k=24$ or $84$). Thus the total 
contribution of the local defects are $24\times {1\over 12}$ 
or $84\times {1\over 42}$, so they are both equal to $2$. 
Because the fibration is iso-trivial, the Weil-Peterson degree vanishes.

 In general, for non-iso-trivial $K3$ fibrations, 
the harmonic energy $\int_{B_0}\varpi_{wp}$ may be partially
 exhausted by hyperbolic singularities in the singular fibers.

\bigskip

\section{A Relative Version of the Universality Theorem} 
\label{section; proof}

\medskip

  The goal of this section is to extend the algebraic 
proof of universality theorem [Liu5] to a relative setting.

 The section \ref{subsection; rel}
 contains the most technical part of the paper in which we outline
a relative extension of the universality theorem. 

\medskip

 The basic object we will study is the modified family invariant
 ${\cal AFSW}_{({\cal Y}/\tilde{B})_{l+1}\mapsto ({\cal Y}/\tilde{B})_l}^{\ast}(1, 
C-2\sum_{i\leq l}E_i)$, which (up to a constant factor) encodes the
 virtual number of $l-$node nodal curves in $C$ along the
 family ${\cal Y}\mapsto \tilde{B}$. They key idea behind our approach is that
 the family blowup formula [Liu2] allows us to express it in terms of
 ${\cal AFSW}_{{\cal Y}\mapsto \tilde{B}}(1, C)$, which encodes the
 virtual number of embedded smooth curves in $C$. I.e. the family blowup
 formula allows us to relate the virtual numbers of immersed 
nodal curves with different genera. This observation will play a crucial
role in the proof of the main theorem in section \ref{section; nodal}.
 
\medskip

 In section \ref{subsection; extension}, we discuss how to extend
 the discussion of the family invariants to include the monodromy
 non-invariant classes as well.

 In section \ref{subsection; genera}, we discuss the correct mixed 
family Seiberg-Witten invariant (and the corresponding 
choice of $\eta_g$) to use in
 defining the virtual number of higher genera nodal curve invariants.
 The class $\eta_g$ is closely related to the Hodge bundle in Gromov-Witten
 theory.

\bigskip

\subsection{The Relative Universal Spaces and the Family Invariants}
\label{subsection; rel}

 In this 
subsection, we review the construction of the relative universal spaces 
[Liu5] and discuss the modified algebraic family Seiberg-Witten
 invariants on them.

 Let ${\cal X}\mapsto B$ be a fibration of
 algebraic surfaces\footnote{The general construction does
 not require the smoothness of the map $\pi$.} over a 
smooth base $B$.  
We construct the relative universal space $({\cal X}/B)_l$ by an 
inductive procedure, parallel to the construction in section 2 of [Liu5].

 For $l=0$, we define $({\cal X}/B)_0$ to be $B$ itself. 
If $({\cal X}/B)_l$ has
 been defined, then define $({\cal X}/B)_{l+1}$
 to be the
 blowing up of $({\cal X}/B)_l\times_{({\cal X}/B)_{l-1}}
 ({\cal X}/B)_l$ along
 the relative
 diagonal $\Delta_l: ({\cal X}/B)_l\mapsto ({\cal X}/B)_l
\times_{({\cal X}/B)_{l-1}}({\cal X}/B)_l$.
 It is easy to see that there is an induced map ${\bf f}_l:
({\cal X}/B)_{l+1}\mapsto 
({\cal X}/B)_l$ commuting with the projection maps
 to the base $B$.  

\bigskip

 The fibers of the relative universal spaces $({\cal X}/B)_l$ 
are exactly the $l-$th universal spaces of
 the fibers of ${\cal X}\mapsto B$.

 On the other hand, the covering morphism $\tilde{B}\mapsto B$ of the
 base spaces induces a pull-back fibration ${\cal X}\times_B
\tilde{B}\mapsto \tilde{B}$. It is not hard to check that
 the construction of the relative universal spaces commutes with
 base change, i.e. we have

$$({\cal X}/B)_l\times_B \tilde{B}=({\cal X}\times_B
\tilde{B}/\tilde{B})_l.$$

\medskip

 Let ${\cal X}\mapsto B$ be a three dimensional 
$K3$ fibration over $B$, i.e. $dim_{\bf C}B=1$, with
 only A-D-E rational double point singularities along the fibers. 
 We apply the covering trick in section \ref{subsection; KV} 
to find a ramified covering
 $\tilde{B}\mapsto B$ ramified over the singular values
 of ${\cal X}\mapsto B$ and by our discussion in the same section
 the pull-back fibration is
 birational to a smooth K3 fiber bundle 
$\pi:{\cal Y}\mapsto \tilde{B}$.

 Then there is a composite morphism ${\cal Y}\mapsto 
{\cal X}\times_B\tilde{B}\mapsto {\cal X}$ of the birational map and
 the ramified covering map. The exceptional locus of the
 birational morphism ${\cal Y}\mapsto {\cal X}\times_B\tilde{B}$
 is a disjoint union of trees of $-2$ rational curves, one for each
rational double point in the fiber of ${\cal X}\mapsto B$.
 This observation 
allows us to transform the nodal curve counting problem
 along ${\cal X}\mapsto B$ to a nodal curve counting
 problem along the enhanced fibration $\pi:{\cal Y}\mapsto \tilde{B}$.

\medskip

 In [Liu1] and [Liu5], we have demonstrated that the counting of 
$l-$node nodal curves in
 an algebraic surface $M$ is closely related to the enumeration of the
 family Seiberg-Witten 
invariant of $C-2\sum_{i\leq l}E_i$ along $M_{l+1}\mapsto M_l$.

 In the following, we discuss the analogue for the relative setting.
 Following the same philosophy, we may want to start with\footnote{The genus
 $g$ and $l$ are related by $g+l={C^2\over 2}+1$. The class $\eta_0=1$ for
 nodal rational curves. For general $g$, the right choice
 of the cycle class $\eta_g$ is
determined in subsection \ref{subsection; genera}.} 
 ${1\over l!}{\cal AFSW}_{({\cal X}/B)_{l+1}\mapsto 
({\cal X}/B)_l}(\eta_g, C-2\sum_{i\leq l}E_i)$ in enumerating the
 $l-$node nodal curves along the family ${\cal X}\mapsto B$.

\medskip

 There are a few questions that we need to address:

\medskip

\noindent (i). How do the family invariants change under the
 covering ${\cal X}\times_B\tilde{B}\mapsto {\cal X}$.

\medskip

\noindent (ii). Are the family invariants unchanged under the 
birational morphism ${\cal Y}/\tilde{B}\mapsto {\cal X}\times_B\tilde{B}/
\tilde{B}$?

\bigskip

\noindent Response to (i).:
 By using the construction of algebraic family Kuranishi models, it is easy
to see that the family invariant along ${\cal X}\times_B\tilde{B}\mapsto
 \tilde{B}$ is
 $[\tilde{B}:B]$-multiple of the family invariant along ${\cal X}\mapsto B$.

\begin{lemm}\label{lemm; cov}
The morphism $({\cal X}\times_B\tilde{B}/\tilde{B})_l\mapsto ({\cal X}/B)_l$ is
 a ramified covering of degree $[\tilde{B}:B]$.
\end{lemm}

\noindent Proof: By the definition of the relative universal spaces, it 
is easy to see that $({\cal X}\times_B\tilde{B}/\tilde{B})_l$
 is the $\tilde{B}\mapsto B$-pull-back of $({\cal X}/B)_l$. So the above
map must be a ramified covering induced by $\tilde{B}\mapsto B$.
$\Box$

 By the above lemma, we can derive that the family invariant of
 $C-2\sum_{i\leq l}E_i$ along
 $({\cal X}\times_B\tilde{B})_{l+1}\mapsto 
({\cal X}\times_B\tilde{B})_l$ is the $[\tilde{B}:B]$-multiple
of the family invariant along $({\cal X}/B)_{l+1}\mapsto 
({\cal X}/B)_l$.

\medskip

\noindent Response to (ii).: The space ${\cal X}\times_B\tilde{B}$ has a 
finite number (assuming $dim_{\bf C}B=1$)
 of isolated singularities. The map ${\cal Y}\mapsto 
 {\cal X}\times_B\tilde{B}$ is the minimal resolution which resolves the
 isolated singularities.

 We have the following commutative diagram of relative universal spaces,

\[
\begin{array}{ccc}
({\cal Y}/\tilde{B})_{l+1} & \mapsto & 
({\cal X}\times_B\tilde{B}/\tilde{B})_{l+1}\\
\Big\downarrow & & \Big\downarrow\\
({\cal Y}/\tilde{B})_l & \stackrel{h_l}{\longrightarrow}
& ({\cal X}\times_B\tilde{B}/\tilde{B})_l
\end{array}
\]

, where the horizontal maps are birational.

 By using this commutative diagram, we can check that
 $$\hskip -.6in
{\cal AFSW}_{({\cal Y}/\tilde{B})_{l+1}\mapsto ({\cal Y}/\tilde{B})_l}(\eta,
 C-2\sum_{i\leq l}E_i)={\cal AFSW}_{({\cal X}\times_B\tilde{B}/\tilde{B})_{l+1}
\mapsto ({\cal X}\times_B\tilde{B}/\tilde{B})_l}(h_{l\ast}\eta, C-2\sum_{i\leq l}
E_i).$$

  Therefore
 we will discuss the nodal curve counting along the smooth $\pi:{\cal Y}\mapsto
 \tilde{B}$.

To count the virtual number of nodal curves in 
$C\in H^{1, 1}({\cal Y}, {\bf Z})_f$ with $l$ nodes, 
 we consider the family $({\cal Y}/\tilde{B})_{l+1}\mapsto
({\cal Y}/\tilde{B})_l$ and start with the following family invariant,
 ${1\over [\tilde{B}:B]}{1\over l!}
{\cal AFSW}_{({\cal Y}/\tilde{B})_{l+1}\mapsto 
({\cal Y}/\tilde{B})_l}(\eta_g, C-\sum_{i\leq l} 2E_i)$ for some suitable  
 $\eta_g\in {\cal A}_{\cdot}(({\cal Y}/\tilde{B})_l)$,
 $g={C^2\over 2}+1-l$. The right choice of $\eta_g$ is addressed in detail in
 subsection \ref{subsection; genera}.

\medskip

The cycle class $\eta_g$ we determine in section \ref{subsection; genera} is 
 a polynomial combination of 
${\bf\pi}_i^{\ast}\hat{C}\in {\cal A}_{\cdot}({\cal Y}_l)$
 and $E_{i; j}\in {\cal A}_{\cdot}(({\cal Y}/B)_l), i\leq j\leq l$
 and the Chern classes of the relative tangent bundles 
$c_{total}({\bf T}({\cal Y}/
\tilde{B})_{k+1}/({\cal Y}/\tilde{B})_k)$, etc.   

 Parallel to the $B=pt$ case, the family invariant contains excess
 contributions of type $I$ exceptional classes. In [Liu1] and [Liu5], we
 had introduced the concept of modified family invariant which 
 captured the residual contribution removing all the excess contributions from
 the various type $I$ exceptional classes. Please consult [Liu5] section
 5.2. for
 its definitions.

 The $\tilde{B}-$relative version 
of the type $I$ modified family
 invariant in [Liu5] defines the modified
 invariant ${\cal AFSW}_{({\cal Y}/\tilde{B})_{l+1}
\mapsto ({\cal Y}/\tilde{B})_l}^{\ast}(\eta_g, C-\sum_{i\leq l} 2E_i)$ of 
the class $C$.

\medskip

Let us state the main result as a theorem,

\begin{theo}\label{theo; main}
Let ${\cal Y}\mapsto \tilde{B}$ be a smooth K3 fiber bundle with
 $dim_{\bf C}\tilde{B}=1$. 
Let $\eta_g$ be a polynomial combination of
 ${\bf \pi}_i^{\ast}C$ 
and $E_{i; j}\in {\cal A}_{\cdot}({\cal Y}_l), i\leq j\leq l$
 and the Chern classes of 
${\bf T}({\cal Y}/\tilde{B})_{k+1}/({\cal Y}/\tilde{B})_k$, 
 then the modified family invariant (by translating the
 definition 13 and 14 of [Liu5] to the relative setting) 

 $${\cal AFSW}_{({\cal Y}/\tilde{B})_{l+1}
\mapsto ({\cal Y}/\tilde{B})_l}^{\ast}(\eta_g, C-\sum_{i\leq l} 2E_i)$$

 can be identified with $l!$ times the virtual number of $l-$node nodal curves
 in the class $C$ and
 can be simplified (by using the family blowup formula
 [Liu2]) to be the product of the algebraic 
family Seiberg-Witten invariant of $C$,
 ${\cal AFSW}_{{\cal Y}\mapsto \tilde{B}}(1, C)$, and the
 degree $l$ universal polynomial of $C^2$ and $c_2(K3)$. 
\end{theo}

\noindent {\bf Sketch} of the proof of the theorem: 
 The proof of the theorem is a slight modification
 of the algebraic proof [Liu5] of the following universality theorem , 
applied to the special case of $M=K3$.

\begin{theo}
Let $l\in {\bf N}$ denote the number of nodal singularities.
Let $L$ be a $5l-1$ very-ample line bundle on an algebraic surface $M$,
 then the number of $l-$node nodal singular curves in a generic 
 $\delta$ dimensional linear sub-system of $|L|$ can be expressed as 
 a universal polynomial (independent to $M$) 
of $c_1(L)^2$, $c_1(L)\cdot c_1(M)$, $c_1(M)^2$, 
$c_2(M)$ of degree $l$.
\end{theo}
 
\medskip

 As the original proof has been lengthy, we do not 
 go through the full details again here. Instead
we emphasize on
 the difference which needs some modification. 

\medskip

 In our discussion, we consider $M=K3$ only, so $c_1(M)=0$.
 In the original set up, the complex family dimension of 
$C-\sum_{i\leq l} 2E_i$ is equal to 
${C^2-C\cdot c_1({\bf K}_M)\over 2}-l={C^2\over 2}-l$. 
In the current relative version,
it is replaced by ${C^2\over 2}-l+dim_{\bf C}\tilde{B}$. 
 In the relative setting, the
 addition of the family dimensions of a class 
$C_0$ and $p$ distinct type $I$ exceptional classes $e_{k_i}$
 obeys the following fiber product axiom, 

\begin{axiom}\label{axiom; fiber} 
 The expected dimension of the family moduli space of co-existence of 
$(C_0, e_{i_1}, e_{i_2}, \cdots, e_{i_p})$ is given by the formula,

$$\hskip -.3in
\{
{C_0^2-c_1({\bf K}_{{\cal Y}/B})\cdot C_0\over 2}+dim_{\bf C}\tilde{B}\}+
\sum_{1\leq k\leq p}\{{e_{i_k}^2-c_1({\bf K}_{{\cal Y}/B})\cdot e_{i_k}\over 2}+
dim_{\bf C}\tilde{B}\}-p\cdot dim_{\bf C}\tilde{B}$$
$$=dim_{\bf C}\tilde{B}+
{C_0^2-c_1({\bf K}_{{\cal Y}/B})\cdot C_0\over 2}+
\sum_{1\leq k\leq p}{e_{i_k}^2-c_1({\bf K}_{{\cal Y}/B})\cdot e_{i_k}\over 2}.$$
\end{axiom}

\medskip

 Let ${\cal M}_{C_0}\mapsto \tilde{B}$ and 
${\cal M}_{e_{i_k}}\mapsto \tilde{B}$ be the family moduli spaces 
of $C_0$ and the type $I$ classes 
$e_{i_k}$, $1\leq k\leq p$.
 We define the family moduli space of co-existence of the tuple  
$(C_0, e_{i_1}, e_{i_2}, \cdots, e_{i_p})$
 to be the fiber product of ${\cal M}_{C_0}$ 
 and ${\cal M}_{e_{i_k}}$, $1\leq i\leq p$. 
The above dimension formula follows from viewing the given fiber product
 as the pull-back of the Cartesian product 
${\cal M}_{C_0}\times^{k\leq p} {\cal M}_{e_{i_k}}$ through the diagonal map 
$\Delta_{\tilde{B}^{1+p}}: \tilde{B}\mapsto \tilde{B}^{1+p}$. 

\medskip

As in the $\tilde{B}=pt$ case, 
we use the admissible graphs $\Gamma\in adm(l)$ and the admissible 
stratification $({\cal Y}/\tilde{B})_l=\coprod_{\Gamma\in adm(l)}{\bf Y}_{\Gamma}$
to stratify the relative universal space $({\cal Y}/\tilde{B})_l$. 
For their definitions, please 
consult [Liu5], section 2.

The admissible graphs $\Gamma$ are finite graphs defined by a few 
combinatorial axioms which characterize the combinatorial
 patterns of $l-$consecutive blowing ups on algebraic surfaces.
 The space ${\bf Y}(\Gamma)$ is of complex 
codimension $codim_{\bf C}(\Gamma)$ in $({\cal Y}/\tilde{B})_l$. 
 The $codim_{\bf C}(\Gamma)$ is equal to the number of $1-$edges in
 $\Gamma$.

 By proposition 4 of [Liu5], each
 ${\bf Y}(\Gamma)$ can be viewed as a regular 
complete intersection of smooth hypersurfaces in $({\cal Y}/\tilde{B})_l$,
 which is the locus of co-existence of a finite collection
 of type $I$ exceptional curves $e_1, e_2, \cdots, e_l$.

\medskip

 For a given family ${\cal Y}\mapsto \tilde{B}$, the discussion of the following
 list of topics is completely parallel to the original $\tilde{B}=pt$ version
 in [Liu5].

\medskip

\noindent (1). Given an admissible graph $\Gamma\in adm(l)$ and the
 corresponding closure of the stratum ${\bf Y}_{\Gamma}$, 
${\bf Y}(\Gamma)$, there are
 a collection of type $I$ exceptional classes, $e_i=E_i-\sum_{j_i}
E_{j_i}$ \footnote{The indexes 
$j_i$ are the direct descendent indexes of $i$ in $\Gamma$.}, 
effective over ${\bf Y}(\Gamma)$.
Let $e_{i_k}, k\leq p$ be the sub-collection of 
type $I$ exceptional classes 
which pair negatively with the given class $C_0=C-\sum_{i\leq l} 2E_i$.

\medskip

\noindent (2). Given the class $e_i$, the universal curves form a ${\bf P}^1$
 fibration. The construction of the relative ${\bf P}^1$ fibrations 
$\Xi_i$ and its relatively minimal model, 
the ${\bf P}^1$ bundle $\tilde{\Xi}_i$. See
 section 3.1. of [Liu5] for more details.

\medskip

\noindent (3). The construction of the bundle class $\tau_{\Gamma}$ 
 by using $\tilde{\Xi}_i$. Please consult definition 10 and lemma 17
 of [Liu5] for the details.

\medskip

\noindent (4). The recursive definitions and constructions of the 
modified family invariant ${\cal AFSW}_{({\cal Y}/\tilde{B})_{l+1}\mapsto
 ({\cal Y}/\tilde{B})_l}^{\ast}(\eta_g\cap c_{total}(\tau_{\Gamma}),
 C-2\sum_{i\leq l}E_l)$. 
 Please consult section 5.2. of [Liu5], where we had taken $\eta_g=1$.

\medskip

\noindent (5). The identification of the residual contribution of
the family invariant ${\cal AFSW}_{({\cal Y}/\tilde{B})_{l+1}\mapsto 
({\cal Y}/\tilde{B})_l}(\eta_g, C-2\sum_{i\leq l}E_i)$
above the top stratum ${\bf Y}_{\gamma_l}\subset
 ({\cal Y}/\tilde{B})_l$ with the
 modified family invariant 

${\cal AFSW}_{({\cal Y}/\tilde{B})_{l+1}\mapsto 
({\cal Y}/\tilde{B})_l}^{\ast}(\eta_g, C-2\sum_{i\leq l}E_i)$ based on the
repeatedly application of residual intersection formula (example 14.1.4. of [F])
 and repeatedly blowing ups along 
$Z(s_{canon})\times_{({\cal Y}/\tilde{B})_l}{\bf Y}(\Gamma)$ for various
 $\Gamma\in \Delta(l)-\{\gamma_l\}$. This is essentially the key argument 
 of the universality theorem in [Liu5]. Please consult its section 6 for
 the details.

\medskip

 Unlike the type $I$ exceptional classes, the moduli spaces of the
 type $II$ exceptional classes are usually non-regular.
 To argue that the modified family invariant 
${\cal AFSW}_{({\cal Y}/\tilde{B})_{l+1}\mapsto 
({\cal Y}/\tilde{B})_l}^{\ast}(\eta_g, C-2\sum_{i\leq l}E_i)$ is $l!$ times the
virtual number of $l-$node nodal curves, we adopt the residual intersection
 theory of type $II$ exceptional curves [Liu6].

\medskip

 Let us survey the basic idea of the identification.  
After resolving the $l$-node nodal points, nodal
curves dual to $C$ are resolved into smooth curves dual to
 $C-2\sum_{i\leq l}E_i$.
On the other hand, dimension analysis in [Liu3] indicates that the 
enumeration of curves in $C-2\sum_{i\leq l}E_i$ 
using ${\cal AFSW}_{({\cal Y}/\tilde{B})_{l+1}
\mapsto ({\cal Y}/\tilde{B})_l}(\eta_g, C-2\sum_{i\leq l}E_i)$ 
 counts not only the smooth curves in $C-2\sum_{i\leq l}E_i$ but also
 the other combinations of curves involving type $II$
 exceptional classes.

 In defining the modified family invariant, all the correction terms in
 $$\hskip -.2in {\cal AFSW}_{({\cal Y}/\tilde{B})_{l+1}\mapsto
 ({\cal Y}/\tilde{B})_l}^{\ast}(\eta_g, C-2\sum_{i\leq l}E_i)=
{\cal AFSW}_{({\cal Y}/\tilde{B})_{l+1}\mapsto
 ({\cal Y}/\tilde{B})_l}(\eta_g, C-2\sum_{i\leq l}E_i)$$
$$\hskip -.2in 
-\sum_{\Gamma\in \Delta(l)-\{\gamma_l\}}{\cal AFSW}_{({\cal Y}/\tilde{B})_{l+1}
\times_{({\cal Y}/\tilde{B})_l}{\bf Y}(\Gamma)\mapsto {\bf Y}(\Gamma)}^{\ast}(
\eta_g\cap c_{total}(\tau_{\Gamma}), C-2\sum_{i\leq l}E_i-\sum_{i\leq p}e_{k_i})$$

 are from the various collections of type $I$ exceptional classes
over the various ${\bf Y}(\Gamma)\subset ({\cal Y}/\tilde{B})_l$, 
 $\Gamma\in \Delta(l)-\{\gamma_l\}$. When $C$ is a high multiple, the type
 $II$ exceptional classes do not contribute to the nodal curve counting.
 Without any assumption on $C$, in principle we
 should consider not only the collections of type $I$ exceptional classes but
also the mixtures of type $I/II$ exceptional classes as well as 
the collections of type $II$ exceptional classes. 

  By applying the residual intersection theory to both the type $I$ and type
 $II$ classes, the theory [Liu6] allows us to
 subtract all the excess contributions of the type $I/II$ exceptional
 curves, the resulting residual contribution is the virtual number of
 smooth curves dual to $C-2\sum_{i\leq l}E_i$ in the universal family, 
and is $l!$ times the
 virtual number of $l$-node nodal curves dual to $C$ 
of the $K3$ fiber bundle ${\cal Y}\mapsto \tilde{B}$.

\medskip

Then the key argument is to show
that all these excess contributions involving type $II$ exceptional classes 
\footnote{Including those from the mixtures of type $I/II$ classes and
 those from type $II$ classes.}
 vanish!  Once this is achieved, we can identify 
 ${1\over l!}{\cal AFSW}_{({\cal Y}/\tilde{B})_{l+1}\mapsto
 ({\cal Y}/\tilde{B})_l}^{\ast}(\eta_g, C-2\sum_{i\leq l}E_i)$ with the
 virtual number of $l$-node nodal curves (dual to $C$) within the family
 ${\cal Y}\mapsto \tilde{B}$.

\bigskip

 When ${\cal Y}\mapsto \tilde{B}$ is not a trivial product,
 ${\bf R}^2\pi_{\ast}\bigl({\cal O}_{\cal Y}\bigr)$ is a non-trivial
 line bundle over $\tilde{B}$. So the formal excess base dimension 
$febd(C, {\cal Y}/\tilde{B})=0$.
 
 Consider a collection of type $I$ and
 type $II$ classes, $e_{k_1}, e_{k_2}, \cdots, e_{k_p}$ and
 $e_{II; 1}, \cdots, e_{II; p'}$, $e_{k_i}\cdot e_{k_j}\geq 0$,
 $e_{II; i}\cdot e_{II; j}\geq 0$, $e_{k_i}\cdot e_{II; j}\geq 0$,
 which pair negatively with 
$C-2\sum_{i\leq l}E_i$.
Following the argument\footnote{The fact that
 ${\cal Y}\mapsto \tilde{B}$ is a $K3$ fiber bundle implies automatically that
 $c_1({\bf K}_{{\cal Y}/\tilde{B}})-e_{II; i}$ is non-nef.}
 of the main theorem 
and section 5.3. in [Liu6], schematically the
 mixed invariant attached to the moduli space of
 co-existence ${\cal M}_{e_{k_1}, \cdots, e_{k_p};
 e_{II; 1}, \cdots, e_{II; p'}}=\times_{\tilde{B}}^{i\leq p}
{\cal M}_{e_{k_i}}\times_{\tilde{B}} 
\times_{\tilde{B}}^{i\leq p'}{\cal M}_{II; i}$ is
 of the form

$$\hskip -.4in
{\cal AFSW}_{({\cal Y}/\tilde{B})_{l+1}\mapsto ({\cal Y}/\tilde{B})_l}(
\eta_g\cap 
h_{\ast}[{\cal M}_{e_{k_1}, \cdots, e_{k_p}; e_{II; 1}, 
\cdots, e_{II; p'}}]_{vir}\cap \tau, 
C-2\sum_{i\leq l}E_i-\sum_{i\leq p}e_{k_i}-\sum_{i\leq p'}e_{II; i}).$$

 Here $h:{\cal M}_{e_{k_1}, \cdots, e_{k_p};
 e_{II; 1}, \cdots, e_{II; p'}}\mapsto ({\cal Y}/\tilde{B})_l$ is the
natural projection map to the base space.

 The vanishing of such mixed invariants has nothing to do with the
 detail structure of the class $\tau$ nor 
$[{\cal M}_{e_{k_1}, \cdots, e_{k_p}; e_{II; 1}, 
\cdots, e_{II; p'}}]_{vir}$.  Firstly, we have the following lemma,

\begin{lemm}\label{lemm; vanishing}
Let ${\cal Y}\mapsto \tilde{B}$ be a non-trivial $K3$ fiber bundle.
Then the push-forward image of $[{\cal M}_{e_{II; i}}]_{vir}$ in
 ${\cal A}_{\cdot}(({\cal Y}/\tilde{B})_l)$ is of the form
 $c_1({\cal R}^2\pi_{\ast}\bigl({\cal O}_{\cal Y}\bigr))\cap \beta$ for
 some $\beta \in {\cal A}_{\cdot}(({\cal Y}/\tilde{B})_l)$.
\end{lemm}

\noindent Proof of the lemma: When ${\cal Y}\mapsto \tilde{B}$ is non-trivial,
 ${\cal R}^2\pi_{\ast}{\cal O}_{\cal Y}$ is a non-trivial invertible sheaf
 over $\tilde{B}$. Likewise 
${\cal R}^2{\bf f}_{l\ast}{\cal O}_{({\cal Y}/\tilde{B})_{l+1}}$, isomorphic to
the pull-back of ${\cal R}^2\pi_{\ast}{\cal O}_{\cal Y}$ to 
$({\cal Y}/\tilde{B})_l$, is also non-trivial.
So $febd(e_{II; i}, ({\cal Y}/\tilde{B})_{l+1}\mapsto 
 ({\cal Y}/\tilde{B})_l)=0$ and the expected dimension of the class
 $e_{II; i}$ is $dim_{\bf C}({\cal Y}/\tilde{B})_l+{e_{II; i}^2-
c_1({\bf K}_{{\cal Y}_{l+1}/{\cal Y}_l})\cdot e_{II; i}\over 2}$.

 Then the assertion in the lemma follows from the argument of the
 algebraic family Kuranishi model in the step 1 of proposition \ref{prop; wall}. 
 $\Box$

\medskip

 When ${\cal Y}$ is three dimensional, 
$c_1({\cal R}^2\pi_{\ast}\bigl({\cal O}_{\cal Y}\bigr))\cap [\tilde{B}]\in
 {\cal A}_0(\tilde{B})$ is a zero dimensional cycle class. After some
simple calculation, this implies that the mixed invariant to be 
subtracted from the original family invariant of $C-2\sum_{i\leq l}E_i$,

$$\hskip -.6in 
{\cal AFSW}_{({\cal Y}/\tilde{B})_{l+1}
\mapsto ({\cal Y}/\tilde{B})_l}(\eta_g\cap 
h_{\ast}[{\cal M}_{e_{k_1}, \cdots, e_{k_p}; e_{II; 1}, 
\cdots, e_{II; p'}}]_{vir}\cap \tau, 
C-2\sum_{i\leq l}E_i-\sum_{i\leq p}e_{k_i}-\sum_{i\leq p'}e_{II; i})$$

 involving
 a finite collection of type $II$ exceptional classes $e_{II; i}$, $i\leq p'$,
 is proportional to

$$\hskip -.4in {\cal AFSW}_{({\cal Y}/\tilde{B})_{l+1}
\mapsto ({\cal Y}/\tilde{B})_l}([({\cal Y}_b)_l]\cap\eta, C-2\sum_{i\leq l}
E_i-\sum_{i\leq p}e_{k_i}-\sum_{i\leq p'}e_{II; i})$$ 

for some $\eta\in {\cal A}_{\cdot}(({\cal Y}/\tilde{B})_l)$ and some
 fiber $({\cal Y}_b)_l$ of $({\cal Y}/\tilde{B})_l$, $b\in \tilde{B}$.

 By a similar argument parallel to lemma \ref{lemm; vanishing}, virtual
 fundamental class of ${\cal M}_{C-2\sum_{i\leq l}E_i-
\sum_{i\leq p}e_{k_i}-\sum_{i\leq p'}e_{II; i}}$, and therefore the
 mixed invariant of $C-2\sum_{i\leq l}E_i-
\sum_{i\leq p}e_{k_i}-\sum_{i\leq p'}e_{II; i}$, is proportional to
 $c_1({\cal R}^2{\bf f}_{l\ast}{\cal O}_{({\cal Y}/\tilde{B})_{l+1}})$. 
Since 
${\cal R}^2{\bf f}_{l\ast}{\cal O}_{({\cal Y}/\tilde{B})_{l+1}}$ is isomorphic
 to the pull-back of ${\cal R}^2\pi_{\ast}\bigl({\cal O}_{\cal Y}\bigr)$ 
 to $({\cal Y}/\tilde{B})_l$, its restriction to
 the fiber $({\cal Y}_b)_l$ of the relative universal space
 $({\cal Y}/\tilde{B})_l$ is trivial. So the mixed family
 invariant found above vanishes. \label{vava}

\bigskip

 In the following, we argue that the modified mixed family invariant
${\cal AFSW}_{({\cal Y}/\tilde{B})_{l+1}\mapsto 
({\cal Y}/\tilde{B})_l}^{\ast}(\eta_g, C-2\sum_{i\leq l}E_i)$ can be  
identified to be a mixed algebraic 
family Seiberg-Witten invariant of $C$,
 ${\cal AFSW}_{{\cal Y}\mapsto \tilde{B}}(\vartheta, C)$, 
where the inserted class $\vartheta\in {\cal A}_{\cdot}(\tilde{B})$ is 
 a polynomial expression of $\pi_{\ast}(\hat{C}^a
c_1({\bf T}_{{\cal Y}/\tilde{B}})^a\cap
 c_2({\bf T}_{{\cal Y}/\tilde{B}})^b\cap \hat{C}^c\cap [{\cal Y}])$.

\medskip 

  Firstly, we recall 
that $({\cal Y}/\tilde{B})_{l+1}$ can be constructed
 from ${\cal Y}\times_{\tilde{B}} ({\cal Y}/\tilde{B})_l$ 
by blowing up consecutively along $l$ sections of  
intermediate fiber bundles\footnote{Consult lemma 3.1. and proposition
 3.1. of [Liu1].}.

 Take $({\cal Y}/\tilde{B})_{l,0}={\cal Y}\times_{\tilde{B}}
 ({\cal Y}/\tilde{B})_l$ and inductively let
 $({\cal Y}/\tilde{B})_{l, k}\mapsto ({\cal Y}/\tilde{B})_l$
 denote the fiber bundle constructed by blowing up the section $s_{k-1}:
({\cal Y}/\tilde{B})_l\mapsto ({\cal Y}/\tilde{B})_{l,k-1}$. 
The above cross 
section $s_{k-1}$ is
the pull-back of the relative diagonal $({\cal Y}/\tilde{B})_k\subset 
({\cal Y}/\tilde{B})_k\times_{({\cal Y}/\tilde{B})_{k-1}}
 ({\cal Y}/\tilde{B})_k$ by the composite map

$${\bf f}_{k}\circ \cdots \circ {\bf f}_{l-1}:
({\cal Y}/\tilde{B})_l\mapsto ({\cal Y}/\tilde{B})_k,$$

 after we have identified $({\cal Y}/\tilde{B})_{l, k-1}$ with
 $({\cal Y}/\tilde{B})_k\times_{({\cal Y}/\tilde{B})_{k-1}} 
({\cal Y}/\tilde{B})_l$. 

\medskip

 In enumerating the (type $I$) modified family invariant 

\noindent ${\cal AFSW}_{
({\cal Y}/\tilde{B})_l\mapsto ({\cal Y}/\tilde{B})_l}^{\ast}
(\eta_g, C-\sum_i 2E_i)$, we apply algebraic family blowup formula [Liu2] 
to the individual terms of the form 
 ${\cal AFSW}_{({\cal Y}/\tilde{B})_{l+1}
\times_{({\cal Y}/\tilde{B})_l}{\bf Y}(\Gamma)\mapsto {\bf Y}(\Gamma)}(
\eta_g\cap c_{total}(\tau_{\Gamma}), C-2\sum_{i\leq l}E_i-
\sum_{i\leq p}e_{k_i})$, where $e_{k_i}\cdot (C-2\sum_{i\leq l}E_i)<0$.  

As the fundamental cycle 
$[{\bf Y}(\Gamma))]\in {\cal A}_{\cdot}(({\cal Y}/\tilde{B})_l)$ 
can be expressed as a polynomial expression of the 
 various $E_{i; j}$, the above mixed family invariant
 can be casted into some mixed invariant 

\noindent ${\cal AFSW}_{({\cal Y}/\tilde{B})_{l+1}\mapsto 
({\cal Y}/\tilde{B})_l}(\eta,
 C-\sum 2E_i-\sum_j e_j)$ for some $\eta$ depending on $\Gamma$. Then the 
family blowup formula in [Liu2] relates this mixed invariant to some 
 mixed invariant 
${\cal AFSW}_{{\cal Y}\times_{\tilde{B}} ({\cal Y}/\tilde{B})_l
\mapsto ({\cal Y}/\tilde{B})_l}(\eta', C)$ of the product
 fiber bundle ${\cal Y}\times_{\tilde{B}} ({\cal Y}/\tilde{B})_l\mapsto 
({\cal Y}/\tilde{B})_l$, where
 $\eta'$ is an polynomial expression of ${\bf \pi}_i^{\ast}\hat{C}$, $E_{i; j}$,
 $c_1(s_{k-1}^{\ast}{\bf T}(({\cal Y}/\tilde{B})_{l, k-1}/
({\cal Y}/\tilde{B})_l))$ and 
 $c_2(s_{k-1}^{\ast} {\bf T}(({\cal Y}/\tilde{B})_{l,(k-1)}/
({\cal Y}/\tilde{B})_l))$ for the various
 indexes $i, j$ and $k$. Because both $C$ and the product fiber bundle
 ${\cal Y}\times_{\tilde{B}}({\cal Y}/\tilde{B})_l\mapsto ({\cal Y}/\tilde{B})_l$
 are pulled back from ${\cal Y}\mapsto \tilde{B}$ by $({\cal Y}/\tilde{B})_l
\mapsto \tilde{B}$,
 its family invariant and the algebraic family Kuranishi model
 over $({\cal Y}/\tilde{B})_l$ are pulled back from $\tilde{B}$ by
 $({\cal Y}/\tilde{B})_l\mapsto \tilde{B}$.

\medskip

 Then by pushing forward along $({\cal Y}/\tilde{B})_l\mapsto \tilde{B}$,
 the mixed invariant of $C$ over $({\cal Y}/\tilde{B})_l$ is equal to
 a mixed invariant of $C$ over $\tilde{B}$. By the functorial property of
 the proper push-forward, 
we factorize the push-forward along $({\cal Y}/\tilde{B})_l\mapsto \tilde{B}$ into
 the push-forward along $({\cal Y}/\tilde{B})_l\mapsto 
\times_{\tilde{B}}^l{\cal Y}$ and the push-forward along 
$\times_{\tilde{B}}^l{\cal Y}\mapsto \tilde{B}$. 

After some standard computation of characteristic classes
 similar to proposition 13 of [Liu5], 
the original mixed family invariant over $({\cal Y}/\tilde{B})_l$
 can be re-expressed as a mixed invariant over the $l-$fold fiber product 
$\times_{\tilde{B}}^l{\cal Y}$, where
 the inserted class 
is a polynomial expression of the products of
 $\pi_i^{\ast}\hat{C}$, $\pi_i^{\ast} c_1({\bf T}{\cal Y}/{\bf T}\tilde{B})$ and 
$\pi_i^{\ast}c_2({\bf T}{\cal Y}/{\bf T}\tilde{B})$, $1\leq i\leq l$.

\medskip

  On the other hand, the fiber product 
$\times_{\tilde{B}}^l{\cal Y}$ is
 the pull-back of ${\cal Y}^l$ over $\tilde{B}^l$ by $\Delta_l:\tilde{B}
\mapsto \tilde{B}^l$.
 Thus, the push-forward of the inserted class
 along $\times_{\tilde{B}}^l{\cal Y}\mapsto \tilde{B}$ 
 can be expressed as a polynomial expression of 
 $\pi_{\ast}(c_1({\bf T}{\cal Y}/{\bf T}\tilde{B})^a
\cap c_2({\bf T}{\cal Y}/
{\bf T}\tilde{B})^b\cap \hat{C}^c\cap [{\cal Y}])\in {\cal A}_{\cdot}(\tilde{B})$.

 The numerical constraint $2\leq a+cb+c\leq 2+dim_{\bf C}B$ is necessary
 for the
 intersection pairing to be non-zero \footnote{In the above argument, we have
not used the $dim_{\bf C}\tilde{B}=1$ condition.}.

 After the complicated push-forward operation, 
the resulting family invariant is of the form 

\noindent ${\cal AFSW}_{{\cal Y}\mapsto 
\tilde{B}}(\vartheta, C)$ for some $\vartheta\in {\cal A}_{\cdot}(\tilde{B})$. 
The class $\vartheta$ is a polynomial expression of 
 $\pi_{\ast}(c_1({\bf T}{\cal Y}/{\bf T}\tilde{B})^a
\cap c_2({\bf T}{\cal Y}/{\bf T}\tilde{B})^b\cap \hat{C}\cap [{\cal Y}])$.

\medskip

When $dim_{\bf C}\tilde{B}=1$, a vanishing argument similar to the vanishing
 argument of type $II$ contributions on page \pageref{vava}
 implies that ${\cal AFSW}_{{\cal Y}\mapsto
 \tilde{B}}(\{\vartheta\}_0, C)=0$. So only the degree one term
 $\{\vartheta\}_1$ contributes to the mixed invariant and the above
 family invariant is reduced to ${\cal AFSW}_{{\cal Y}\mapsto 
\tilde{B}}(\{\vartheta\}_1, C)$.

\begin{lemm}\label{lemm; select}
Let $\vartheta\in {\cal A}_{\cdot}(\tilde{B})$ be a polynomial expression of
 the variables $\pi_{\ast}(c_1({\bf T}{\cal Y}/{\bf T}\tilde{B})^a
\cap c_2({\bf T}{\cal Y}/{\bf T}\tilde{B})^b\cap \hat{C}^c)$.
Then when $dim_{\bf C}\tilde{B}=1$, only the variables corresponding to 
$(a, b, c)=(0, 0, 2)$ and $(0, 1, 0)$ tuples contribute to $\{\vartheta\}_1$.
\end{lemm}

\noindent Proof: When $dim_{\bf C}\tilde{B}=1$, $dim_{\bf C}{\cal Y}=3$.
 On the one hand, $a+2b+c\leq 3$ for $c_1({\bf T}{\cal Y}/{\bf T}\tilde{B})^a
\cap c_2({\bf T}{\cal Y}/{\bf T}\tilde{B})^b\cap \hat{C}^c\cap [{\cal Y}]$ to be
 nonzero on ${\cal Y}$. On the other hand, $a+2b+c\leq 2$ in order for the
 push-forward to be of degree one. By the fact that the fibers are $K3$,
 the tuples $(1, 0, 1)$ and $(2, 0, 0)$ can be ruled out. The lemma follows.
  $\Box$

So $\{\vartheta\}_1$ can be reduced to a universal degree $l$ polynomial of
 $\pi_{\ast}(C^2)$ and $c_2(K3)$.

 This ends the proof of the theorem. $\Box$

\bigskip

\begin{rem}\label{rem; same}
 When we take $l={C^2\over 2}+1$, $g=0$ and $\eta_0=1$, it is not hard to see
 from the above discussion that the universal polynomials found in
this theorem are identical to the specific universal polynomials coding
 ``the numbers of
 nodal curves'' of $K3$. The existence of such polynomials has been
 guaranteed by the universality theorem [Liu1], [Liu5]. Given an algebraic surface 
$M=K3$, and a class $C\in H^2(M, {\bf Z})$,
these universal polynomials encode the information of the virtual numbers of
 nodal rational curves on $M$ dual to $C$.
 Through a ${\cal C}^{\infty}$ argument [Liu1] of
 Taubes ``SW=Gr'' [T1], [T2] and [T3], the generating function of such universal
 polynomials can be identified with the well known Yau-Zaslow formula.
\end{rem}

\begin{rem}\label{rem; different}
 When $dim_{\bf C}\tilde{B}>1$, $\{\vartheta\}_k$, $k\not=1$, can also
 contribute to the family invariant. 
Even for $l={C^2\over 2}+1$ and $g=0$
the polynomial expression $\vartheta$ 
found in the above theorem is usually different from the
 one predicted by Yau-Zaslow [YZ] formula.
\end{rem}

\begin{rem}
 The above vanishing argument of the type $II$ class contributions 
is parallel to the similar vanishing argument (of the nodal curves counting) for 
hyperkahler families of $K3$. Please consult [Liu6], section 4.3.1. 

\medskip

Please also refer to section 1.3.1 of [Liu7] for its similarity with the
vanishing of second sheaf cohomology of ample line bundles 
on $K3$, where the nodal curve
 enumeration has been viewed as an enumerative Riemann-Roch theorem.
\end{rem}
  
\bigskip

\subsection{The Extension of the Family Invariant to Non-Monodromy 
Invariant Classes}\label{subsection; extension} 

\bigskip
 In the previous discussion, we have 
paid our attention to monodromy invariant classes
 $C$ along the fiber bundle ${\cal Y}\mapsto \tilde{B}$. In this subsection 
we extend our discussion and drop this restriction. 

  Let $M$ be an algebraic surface and let $\pi:
{\cal Y}\mapsto \tilde{B}$ be
a smooth morphism of connected algebraic varieties such that
 $\pi^{-1}(b)\cong M$ for some $b\in B$. 
 
 If the fundamental group $\pi_1(\tilde{B}, b)\not=\{1\}$, 
then it induces a monodromy representation on the middle cohomology of the
 fiber above $b$, $H^2(\pi^{-1}(b), {\bf Z})$.  

Let ${\bf \rho}:\pi_1(\tilde{B}, b)\mapsto
 Aut(H^2(M, {\bf Z}))$ denote the monodromy representation.

 Given a cohomology class $C\in H^2(M, {\bf Z})$, 
 let ${\bf O}_C=Im({\bf \rho}(
\pi_1(\tilde{B}, b)))\cdot C$ denote the orbit of $C$ under the monodromy group 
action.

 Unless the orbit ${\bf O}_C=\{C\}$, by the spectral sequence argument 
the class $C$ does not define a 
 cohomology class on the total space ${\cal Y}$. 

\medskip

We separate our discussion into two different cases, depending on
whether the cardinality $|{\bf O}_C|$ is finite or infinite.

 If the cardinality of 
${\bf O}_C$ is finite, one takes the subgroup 
$G\subset \pi_1(\tilde{B}, b)$ which fixes $C$. 

 It follows from $|{\bf O}_C|<\infty$ that 
$G$ is of finite index in $\pi_1(\tilde{B}, b)$.  Given the stabilizer $G$,
 we consider a finite covering $g:\hat{B}\mapsto \tilde{B}$ such that
 $\pi_1(\hat{B}, \hat{b})=G$.  
Because $C$ is fixed under $\pi_1(\hat{B}, \hat{b})$,
 the class $C$ determines a monodromy invariant class on
 the fiber product $\hat{B}\times_{\tilde{B}} {\cal Y}$.

\begin{defin}\label{defin; finite}
  Under the assumption $|{\bf O}_C|<\infty$, we define 

${\cal AFSW}_{{\cal Y}\mapsto \tilde{B}}( \eta, C)$ to be 
$${{\cal AFSW}_{{\cal Y}\times_{\tilde{B}}\hat{B}\mapsto 
\hat{B}}(g^{\ast}\eta, C)\over |\pi_1(\tilde{B}, b)/\pi_1(\hat{B},
 \hat{b})|}=
{{AFSW}_{Y\times_{\tilde{B}} \hat{B}\mapsto \hat{B}}
(g^{\ast}\eta, C)\over |{\bf O}_C|},$$
 for an arbitrary $\eta\in {\cal A}_{\cdot}(\tilde{B})$.
\end{defin}

 When $|{\bf O}_C|=1$, the above definition
 is reduced to the usual one.
 It is not manifest
 that ${\cal AFSW}_{{\cal Y}\mapsto \tilde{B}}(\eta, C)
\in {\bf Q}$ is always ${\bf Z}$ valued.

 Following the same idea, we define the orbit family 
 invariant of the orbit ${\bf O}_C$,

\begin{defin}\label{defin; orbit}
 Define ${\cal AFSW}_{{\cal Y}\mapsto \tilde{B}}(\eta, {\bf O}_C)$ to be
 ${\cal AFSW}_{{\cal Y}\times_{\tilde{B}}\hat{B}\mapsto \hat{B}}
(g^{\ast}\eta, C)$.
\end{defin}

 When the orbit ${\bf O}_C$ is not a finite set, the group $G\subset 
\pi_1(\tilde{B}, b)$
 is not of finite index. So the covering space
$\hat{B}$ is not a finite covering of $\tilde{B}$. 
In particular, it is non-compact.

In this case, the following vanishing proposition
 implies that the orbit invariant ${\cal AFSW}_{{\cal Y}\mapsto 
\tilde{B}}(\eta, {\bf O}_C)$ should be defined to be 
 zero.

\begin{prop}\label{prop; inf}
Let ${\cal M}_C\mapsto \tilde{B}$ be the family moduli space of
 algebraic curves along ${\cal Y}\times_{\tilde{B}}\hat{B}\mapsto 
\hat{B}$, dual to the class $C$. 
Under the assumption that $|{\bf O}_C|$ is infinite,
 The space ${\cal M}_C$ is empty.
\end{prop}

\noindent Proof of proposition \ref{prop; inf}: Because
 $C$ is monodromy invariant along $\hat{B}\times_{\tilde{B}} 
 {\cal Y}\mapsto \hat{B}$, it defines a cohomology class
 on all the fibers. For simplicity, we denote it by the
same notation $C$.
 We prove the emptiness of ${\cal M}_C$ by showing that
 the class $C$ can never be of type $(1, 1)$ in the fibers of
 ${\cal Y}\times_{\tilde{B}}\hat{B}\mapsto \hat{B}$.

\medskip

 Suppose that $C$ is of type $(1, 1)$ along the fiber algebraic
surface above some point $w\in \hat{B}$, then we argue that
 it is still of type $(1, 1)$ above all points in the
 pre-image $g^{-1}(g(w))$. It is trivial to prove this if
 $p_g(M)=0$. So we may assume that $p_g(M)>0$.
 Let $\Omega_1, \Omega_2, \cdots, \Omega_{p_g}$ be
 a basis of holomorphic $(2, 0)$ forms on ${\cal Y}\times_{\tilde{B}} \{g(w)\}$. 
 Because $C$ is real, it is of type $(1, 1)$ over $\pi^{-1}(g(w))$ iff 
 $\int_{\pi^{-1}(g(w))} C\cup \Omega_i=0$ for all $1\leq i\leq p_g$.

  Because the fibration $\tilde{Y}\times_{\tilde{B}} 
\hat{B}\mapsto \hat{B}$ is pulled back from $\tilde{B}$, these 
$\Omega_i$ can still be viewed as
 holomorphic $(2, 0)$ forms on the fibers above all $g^{-1}(g(w))$. 
 Once $C$ is of type $(1, 1)$ above $w$, it is of type $(1, 1)$ above 
 all $g^{-1}(g(w))$. 
 So the result follows.
 
\medskip

 Because the orbit ${\bf O}_C$ is infinite, this implies that
 there are an infinite number of elements in 
$H^2({\cal Y}\times_{\tilde{B}}\{g(w)\}, {\bf Z})$, which are all of 
type $(1, 1)$
 and are in a single orbit under the monodromy action. We argue
that this is impossible.

Let $\omega$ be the restriction of the relative ample polarization
 on ${\cal Y}\mapsto \tilde{B}$ to the fiber algebraic surface 
 ${\cal Y}\times_{\tilde{B}} \{g(w)\}$. It is of type $(1, 1)$ and is
 apparently monodromy invariant. Because the cup product
 pairing is preserved under the monodromy action, one can
 show easily that the whole monodromy group orbit of $C$, ${\bf O}_C$, has
 an identical pairing with $\omega$. After replacing 
 $C$ by $C+k\omega$, $k\gg 0$, if necessary, we may always assume 
that the self-intersection number \footnote{It is the same for
the whole orbit.} of the whole orbit $(=C\cdot C)$ is positive.
 
 By Hodge index theorem of algebraic surfaces, the subspace of
 $(1, 1)$ classes is of signature $(1, b_2^-)$. The classes 
 with positive self-intersection pairings and positive pairings with
 $\omega$ are in the forward light-cone.
The above
derivation implies the existence of an infinite number of
 lattice elements in the forward light-cone, which have a
 fixed (bounded) pairing with the class $\omega$. On the other
 hand, for any fixed $K$ 
the set $\{t|t\in H^2( {\cal Y}\times_{\tilde{B}}\{g(w)\}, {\bf R})
 \cap H^{1, 1}({\cal Y}\times_{\tilde{B}}\{g(w)\}, {\bf C}), 
 t\cdot t\geq 0, t\cdot \omega<K\}$  is a compact set. It can
never contain an infinite number of lattice elements. So we get
a contradiction!   $\Box$

\subsection{The Virtual Numbers of Higher Genera Nodal Curves in a 
Three Dimensional $K3$ Fiber Bundle}
\label{subsection; genera}

\bigskip

 In this subsection, we find the appropriated $\eta_g$ and 
define the mixed
 algebraic family Seiberg-Witten 
invariants which enumerate high genera immersed nodal curves along
 the $K3$ fiber bundle. 

As a genus $g$ curve in the class $C$ has to develop 
 $l={C^2\over 2}+1-g$ nodes generically, naively
 one may try ${1\over l!}
{\cal AFSW}_{({\cal Y}/\tilde{B})_{l+1}
\mapsto ({\cal Y}/\tilde{B})_l}^{\ast}(1, 
C-\sum_{i\leq l} 2E_i)$ to resemble the virtual number of genus $g$
 nodal curves along the family ${\cal Y}\mapsto \tilde{B}$. 
However there exists some subtlety about family dimensions that we need to take
care.

\medskip

 The following proposition characterizes the correct mixed invariant 
 enumerating the virtual number of genus $g$ nodal curves.

\begin{prop}\label{prop; virtual}
 Let $C$ be a fiberwise cohomology class in $H^2({\cal Y}, {\bf Z})_f$
 determined by the cycle class $\hat{C}\in {\cal A}_1({\cal Y})$.
 For any given $g\leq {C^2\over 2}+1$, there exists a cycle class $\eta_g$
 in ${\cal A}_{C^2-3g+3}(({\cal Y}/\tilde{B})_{{C^2\over 2}+1-g})$,
  expressible as a universal polynomial of 
${\bf \pi}_i^{\ast}\hat{C}$, $E_{i; j}$, and
 the Chern classes of the relative tangent bundles
 ${\bf T}_{({\cal Y}/\tilde{B})_{i+1}/({\cal Y}/\tilde{B})_i}$, $i, j\leq 
l=(C^2/2-g+1)$,
such that 
 $${1\over (C^2/2-g+1)!}{\cal AFSW}_{({\cal Y}/\tilde{B})_{l+1}
\mapsto ({\cal Y}/\tilde{B})_l}^{\ast}(\eta_g, C-\sum_{i\leq l} 2E_i)$$
represents the virtual number of genus $g$ nodal curves dual to $C$
 along ${\cal Y}\mapsto \tilde{B}$.
\end{prop}

\medskip
\noindent Proof of proposition \ref{prop; virtual}: The grade 
 $C^2-3g+3=2l+dim_{\bf C}\tilde{B}-g$ of $\eta_g$ indicates that it
 represents a cycle of codimension $g$ in $({\cal Y}/\tilde{B})_l$.

  Consider the Poincare dual $[A]\in H_2({\cal Y}, {\bf Z})$ of $C\cup F$. 
 The expected dimension of the Gromov-Witten
 invariant of genus $g$ ($g>1$) maps into ${\cal Y}$
 is given by 
$$c_1({\cal Y})\cap [A] -3(g-1)+3(g-1)=c_1({\cal Y})\cap [A].$$

Besides the genus $g$ immersed curves, Gromov-Witten invariant counts
 multiple coverings of embedded maps as well. 
So the expected dimensions of the genus $g$ immersed curves is equal
the expected dimension of genus $g$ Gromov-Witten invariant and is given
 by the above formula.

 Because of the fiber bundle structure of ${\cal Y}\mapsto \tilde{B}$ 
and because of the special fiberwise nature of the class $A$, 
the above dimension formula can be
 reduced to $c_1({\bf T}_{{\cal Y}/\tilde{B}})\cap [A]=0$, because
the fibers are smooth $K3$ surfaces. 

 On the other hand, the expected dimension of the 
algebraic family Seiberg-Witten invariant of the class $C-\sum_{i\leq l}2E_i$ 
along the relative universal family
 $({\cal Y}/\tilde{B})_{l+1}
\mapsto ({\cal Y}/\tilde{B})_l$, induced from the given
 $\tilde{B}-$family of algebraic surfaces ${\cal Y}$, 
is given by $$dim_{\bf C}\tilde{B}+{C^2-c_1({\bf K}_{{\cal Y}/\tilde{B}})\cdot
 C\over 2}-({C^2\over 2}-g+1)=g$$

 using $dim_{\bf C}B=1$ and $c_1({\bf K}_{{\cal Y}/\tilde{B}})\cup F=0$.

\medskip

  This indicates that
 there is a $g$ dimension difference between the two different expected  
dimensions. It is caused by an additional complex rank $g$ 
 obstruction bundle absent in the family Seiberg-Witten theory.

 In the following, we explain how does the obstruction bundle
 appears and define $\eta_g$ to represent its top Chern class.

  Because of the fiber bundle structure ${\cal Y}\mapsto \tilde{B}$,
 the normal bundle of any given fiber ${\bf \pi}^{-1}(b')
\subset {\cal Y}$ is trivial. Let $f:\Sigma\mapsto {\cal Y}$ be a holomorphic
 map from a genus $g$ curve into some fiber of ${\cal Y}$, $g_{\ast}[\Sigma]=A$.
Given such an $f$, the obstruction space of the infinitesimal 
deformations along the normal direction of the fiber $K3$ is 
$g$ dimensional and is absent in the family Seiberg-Witten theory. 

The obstruction space can be identified with
  $H^1(\Sigma, f^{\ast}{\cal O}_{\pi^{-1}(b')})\cong H^1(\Sigma, 
{\cal O}_{\Sigma})$, and by curve Riemann-Roch
 its dimension is $-\chi({\cal O}_{\Sigma})+1
=g$, exactly what we have expected for.

 When the point 
$[f]$ moves, the $g$ dimensional obstruction space forms a rank $g$ vector 
bundle.  
The rank $g$ vector bundle is nothing but the dual of the hodge bundle
 in Gromov-Witten theory.
Let $\tilde{\pi}:
\tilde{\Sigma}\mapsto \bar{\cal M}_{g, n}$ denote the universal curve over the
 compactified moduli space of genus $g$ curves with $n$ marked points. 
 Then from the domain curve point of view the rank $g$ obstruction bundle
can be identified with the first derived image bundle 
${\bf R}^1\tilde{\pi}_{\ast}({\cal O}_{\tilde{\Sigma}})$ over $\bar{\cal M}_{g, n}$.

\bigskip
 
 In the following, we give a target space (i.e. our universal space) 
interpretation of
this rank $g$ obstruction bundle on $({\cal Y}/\tilde{B})_l$. 

 As before let $A$ be the Poincare dual of $C\cup F$, where we assume
 $C$ to be a monodromy invariant class on ${\cal Y}\mapsto \tilde{B}$.  

Consider the 
class $C-2\sum_{i\leq l} E_i$. For all $l$-node 
nodal curves representing $\hat {C}\in {\cal A}_1({\cal Y})$ (dual to $C\in
 H^2({\cal Y}, {\bf Z})_f$),
 the strict transforms of the resolved curves are
 smooth and are in cycle classes of the form $\hat{C}-2\sum_{i\leq l}\hat{E}_i$. 

 So we may consider all effective curves in the class of the form
 $\hat{C}-2\sum_{i\leq l}\hat{E}_i$ along the family
 ${\bf f}_l: ({\cal Y}/\tilde{B})_{l+1} \mapsto ({\cal Y}/\tilde{B})_l$. They form
 the algebraic family moduli space 
${\cal M}_{C-2\sum_{i\leq l}E_i}$ over $({\cal Y}/\tilde{B})_l$.

 The space ${\cal M}_{C-2\sum_{i\leq l} E_i}$ contains the sub-moduli space
 which is the closure of the sub-space of irreducible smooth curves. The
 particular sub-moduli is usually not open nor dense in ${\cal M}_{C-
\sum_{i\leq l}2E_i}$.

 The technique of [Liu5], and [Liu6] enables us to 
separate the virtual fundamental
 class of the sub-moduli from the whole 
$[{\cal M}_{C-2\sum_{i\leq l} E_i}]_{vir}$. The family invariant attached
 to the sub-moduli is nothing but the modified family invariant 
 ${\cal AFSW}^{\ast}$ used in theorem \ref{theo; main} in 
section \ref{section; proof}. 

 Consider
the push-forward ${\cal R}^1{\bf f}_{l\ast}\bigl(
{\cal O}_{\Sigma_{C-2\sum_{i\leq k}E_i}}\bigr)$ along the
 universal curve \footnote{Notice that we have used an identical notation
 ${\bf f}_l$ to denote the map 
$({\cal Y}/\tilde{B})_{l+1} \mapsto ({\cal Y}/\tilde{B})_l$.
 Knowing that the universal curve $\Sigma_{C-2\sum_{i\leq l}E_i}$ can 
 be embedded into $({\cal Y}/\tilde{B})_{l+1}\times_{({\cal Y}/\tilde{B})_l}
 {\cal M}_{C-\sum_{i\leq l}2E_i}$, this is a minor abuse of notations.}
 ${\bf f}_l:\Sigma_{C-2\sum_{i\leq l}E_i}\mapsto 
{\cal M}_{C-2\sum_{i\leq l}E_i}$.
 By semi-continuity theorem and corollary 12.9, on page 288 of 
[Ha], and curve Riemann-Roch formula,
 this first direct image sheaf is locally free of rank $g$. Under the
 natural projection \footnote{The map is constructed by
 the composition $({\cal Y}/\tilde{B})_{l+1}\mapsto {\cal Y}\times_{\tilde{B}}
({\cal Y}/\tilde{B})_l\mapsto {\cal Y}$.}
 map $({\cal Y}/\tilde{B})_{l+1}\mapsto {\cal Y}$, 
the universal curve 
 is projected onto singular curves in $\hat{C}$ along the fibers of 
${\cal Y}$. At each smooth irreducible genus $g$ fiber $\Sigma$ 
of the universal curve parametrized by the sub-moduli 
$\subset {\cal M}_{C-2\sum_{i\leq 
l}E_i}$, the singular curve in ${\cal Y}$ is a $l-$node nodal curve.
 The projection map onto the target nodal curve can be viewed tautologically
 as a holomorphic map $\Sigma\mapsto {\cal Y}$. Because the map is
 an immersion into ${\cal Y}$, there is no non-trivial automorphism of
 $\Sigma$ which fixes the map. Under such an identification of
 $l-$node nodal curves and holomorphic maps, the restriction of the bundle
${\bf R}^1{\bf f}_{l\ast}\bigl(
{\cal O}_{\Sigma_{C-2\sum_{i\leq k}E_i}}\bigr)$ can be identified
 with the restriction of (the appropriate pull-back of) 
${\bf R}^1\tilde{\pi}_{\ast}({\cal O}_{\tilde{\Sigma}})$.

Such an identification can
 be extended to the closure of the sub-space of smooth irreducible curves in
 $\hat{C}-2\sum_{i\leq l}\hat{E}_i$.

 So we may
 view the derived image bundle ${\bf R}^1{\bf f}_{l\ast}\bigl(
{\cal O}_{\Sigma_{C-2\sum_{i\leq l}E_i}}\bigr)$ as the
 family theory analogue of the dual Hodge bundle in Gromov-Witten theory.
 
 \medskip

\begin{defin}\label{defin; eta}
Define $\eta_g=c_g({\cal R}^1{\bf f}_{l\ast}\bigl(
{\cal O}_{\Sigma_{C-2\sum_{i\leq l}E_i}}\bigr))\cap [({\cal Y}/\tilde{B})_l]$.
\end{defin}

\medskip

 In the following, we explain how to determine the Chern classes
 of ${\bf R}^1{\bf f}_{l\ast}\bigl(
{\cal O}_{\Sigma_{C-2\sum_{i\leq l}E_i}}\bigr)$ inductively.
 
\medskip

\begin{lemm}\label{lemm; express}
The class $\eta_g$ can be expressed as a polynomial of 
$\pi_i^{\ast}\hat{C}$, $E_{i;j}$ \footnote{
The symbol $E_{i;j}$ denotes the exceptional divisor in $({\cal Y}/\tilde{B})_l$
above the $(i, j)$-th partial diagonal of $\times_{\tilde{B}}^l{\cal Y}$. The
 same
 notation notation been used in [Liu5] extensively.}
and $c_{total}({\bf T}_{({\cal Y}/
\tilde{B})_{i+1}/({\cal Y}/\tilde{B})_i})$, etc. for
 $i, j\leq l$. 
\end{lemm}

\noindent Proof of lemma \ref{lemm; express}:
Firstly ${\bf R}^0{\bf f}_{l\ast}\bigl(
{\cal O}_{\Sigma_{C-2\sum_{i\leq l}E_i}}\bigr)$ is isomorphic
 to the structure sheaf of ${\cal M}_{C-2\sum_{i\leq l}E_i}$, so 
$c_g({\bf R}^1{\bf f}_{l\ast}\bigl(
{\cal O}_{\Sigma_{C-2\sum_{i\leq l}E_i}}\bigr))=c_g(-{\bf f}_{l\ast}
\bigl({\cal O}_{\Sigma_{C-2\sum_{i\leq k}E_i}}\bigr)$.

For $1\leq k\leq l$, consider the following short exact 
sequences,

$$0\mapsto {\cal O}_{2E_k}(-\hat{C}+2\sum_{i\leq k}E_i)\mapsto 
{\cal O}_{\Sigma_{C-2\sum_{i\leq k-1}E_i}}\mapsto {\cal O}_{\Sigma_{C-
2\sum_{i\leq k}E_i}}\mapsto 0.$$

By taking their right derived sequences along 
 ${\bf f}_l: ({\cal Y}/\tilde{B})_{l+1}\mapsto 
({\cal Y}/\tilde{B})_l$, we may re-express $\eta_g$ as 
$$c_g(\oplus_{1\leq k\leq l}{\bf f}_{l\ast}{\cal O}_{2E_k}(\hat{C}
-2\sum_{i\leq k}E_i))\cap [({\cal Y}/\tilde{B})_l].$$

 The Chern classes of the direct sums of direct images can be computed by
using Grothendieck-Riemann-Roch theorem along the
 ${\bf P}^1$ fibrations (the exceptional divisors $E_k$). The
computation is completely parallel to the family blowup formula and we
 omit the details. By going through the computation similar to
 proposition 13 and lemma 20 of [Liu5], it
 can be expressed as a universal polynomial
expression of $\pi_i^{\ast}\hat{C}$, $E_{i; j}$ and
 $c_{total}({\bf T}_{({\cal Y}/
\tilde{B})_{i+1}/({\cal Y}/\tilde{B})_i})$, etc. $\Box$

 Finally the modified mixed algebraic family Seiberg-Witten
invariant which enumerates the virtual number of immersed nodal 
curves
should be ${1\over l!}{\cal AFSW}_{({\cal Y}/\tilde{B})_{l+1}\mapsto
 ({\cal Y}/\tilde{B})_l}^{\ast}(\eta_g, C-2\sum_{i\leq l} E_i)$, with 
the given $\eta_g$ defined in definition \ref{defin; eta}. 
$\Box$

\section{The Counting of Nodal Curves on an Algebraic $K3$ fiber Bundle}
\label{section; nodal} 

In this section, we would like to apply the general machineries we  have developed
 in the previous sections to deal with the enumeration of virtual
 numbers of nodal curves in a $K3$ fiber bundle.

  Recall the concept of lattice polarized $K3$ surface introduced By
 Dolgachev [D]. Let ${\bf M}$ be 
 sub-lattice of the $K3$ lattice ${\bf L}=3{\bf H}\oplus -2E_8$ 
(with an even quadratic form) 
of signature $(1, m-1)$. A marked ${\bf M}$ polarized
 algebraic $K3$ surface is by definition a pair $(X, \phi)$ such 
that,

\medskip

 (i). $X$ is an algebraic $K3$ surface.

 (ii). The marking 
$\phi: {\bf L}\stackrel{\cong }{\longrightarrow} H^2(X, {\bf Z})$ 
is an isomorphism of lattices which identifies ${\bf M}$ with
 a sub-lattice of $H^2(X, {\bf Z})$.

 By the construction in [D], the complex moduli space of 
${\bf M}-$polarized algebraic $K3$ surface $\underline{{\cal M}}_{\bf M}$ 
forms a bounded symmetric domain of type $IV$,
 which is a ramified quotient from 
the complex hyperbolic space of dimension
 $20-rank_{\bf Z}{\bf M}$.  The moduli space parametrizes
 the isomorphism classes of marked ${\bf M}-$polarized $K3$ surfaces. Roughly speaking, the moduli space parameterizes the
 $K3$ surfaces whose Picard lattice is at least as large as ${\bf M}$ and the
 marking provides the necessary level structure.

\medskip

 Recall the following 
definition of transcendental lattice ${\bf T}$,

\begin{defin}\label{defin; tran}
Define the transcendental lattice 
${\bf T}\equiv {\bf M}^{\bot}\subset {\bf L}$ to be the sub-lattice consisting
 of all the elements in ${\bf L}$ which are perpendicular to ${\bf M}$. 
\end{defin}

  Only when ${\bf M}$ is unimodular, ${\bf T}$ can be unimodular and then ${\bf L}={\bf M}\oplus {\bf T}$.
 Suppose that through the marking $\phi$ the lattice ${\bf M}$ is identified with the Picard
lattice of the $K3$, 
then ${\bf T}$ can be thought to be the sub-lattice of 
transcendental (non-algebraic) 
classes in $L$ which are perpendicular to all $(1, 1)$ classes. 

\medskip

 The dimension of $\underline{{\cal M}}_{\bf M}$ is given by, 

\begin{lemm}\label{lemm; dim}
 $dim_{\bf C}\underline{{\cal M}}_{\bf M}$ is equal to
 $dim_{\bf Z}{\bf T}-2=20-dim_{\bf Z}{\bf M}$.
\end{lemm}

\noindent Proof: The lemma is proved by identifying the tangent 
space of $\underline{{\cal M}}_{\bf M}$ with the subspace of 
the space of infinitesimal complex 
deformations of $X$, $H^1(X, \Theta_X)\cong H^1(X, \Omega^1_X)\subset H^2(X, {\bf C})$,
 perpendicular to the elements in ${\bf M}$ under the cup product pairing.
 Because the cup product pairing 

$$\cup: H^2(X, {\bf C})\otimes_{\bf Z} H^2(X, {\bf Z})\mapsto H^4(X, {\bf C})\cong
 {\bf C}$$ is non-degenerated, the dimension of
 $T\underline{{\cal M}}_{\bf M}$ is equal to $dim_{\bf C} H^1(X, {\Theta}_X)-
dim_{\bf Z}{\bf M}$.  On the other hand, $dim_{\bf Z}{\bf T}+dim_{\bf Z}{\bf M}=
22$. This implies the equality in the lemma. $\Box$

\medskip

 If ${\bf T}$ can be decomposed further into ${\bf H}\oplus \overline{\bf M}$, then
 $\overline{\bf M}$ is of signature $(1, 20-m)$. The Mirror conjecture
of Dolgachev-Gritsenko-Nikulin [D], [GN2] relates the moduli spaces of 
marked ${\bf M}-$polarized 
$K3$ surfaces and of the marked $\overline{\bf M}-$polarized $K3$ surfaces.

\medskip

 Suppose that $\pi:{\cal Y}\mapsto \tilde{B}$ is a relative
 algebraic $K3$ fiber bundle, $dim_{\bf C}\tilde{B}=1$, 
with a cross section 
$s: \tilde{B}\mapsto {\cal Y}$. 
 Under the additional hypothesis that 
$H^{2,0}({\cal Y}, {\bf C})=0$,
 then all the classes in $H^2({\cal Y}, {\bf Z})$ 
are of type $(1, 1)$ and
 they correspond bijectively to the first Chern classes of
 holomorphic line bundles on ${\cal Y}$. As before, we use $F$ to
denote the cohomology class of the fibers. 

 On the one hand, the triple intersection pairing of $H^2({\cal Y}, 
{\bf R})$ can be restricted to a quadratic pairing of 
 $H^2({\cal Y}, {\bf Z})_{free}$ by the formula 
 $<{\bf a}, {\bf b}>=\int_{\cal Y}{\bf a}\cup {\bf b}\cup F$.
 By the property that $F\cup F=0$, it is apparent that the
 pairing descends to the quotient $H^2({\cal Y}, {\bf Z})_{free}$.
 This motivates us to define,

\begin{defin}\label{defin; M}
Define 
${\bf M}=H^2({\cal Y}, {\bf Z})_{free}/{\bf Z}F$, equipped with
the intersection form $<\cdot, \cdot>$ above.
\end{defin}

 Because the Leray spectral sequence of $\pi:{\cal Y}\mapsto \tilde{B}$
 degenerates at the $E_2$ terms, the ${\bf Z}$ module
 $H^2({\cal Y}, {\bf Z})$ can be identified with the direct
sum of the monodromy invariant part $H^2(\pi^{-1}(b), 
{\bf Z})^{\pi_1(\tilde{B}, b)}$ and $H^2(\tilde{B}, {\bf Z})\cong
 {\bf Z}F$. 
  
Therefore this identification induces an embedding of ${\bf M}$
 into $H^2(\pi^{-1}(b), {\bf Z})\cong {\bf L}$ and ${\bf M}$ is the
 monodromy invariant part of $L$.

 \medskip

Once the lattice ${\bf M}$ has been fixed, 
the relative algebraic $K3$ fiber bundle ${\cal Y}\mapsto \tilde{B}$ determines
 a holomorphic map $\Phi_{\cal Y}:\tilde{B}\mapsto 
\underline{\cal M}_{\bf M}$ from the
base of the fibration $\tilde{B}$ to the moduli space of ${\bf M}$-polarized
 $K3$ surfaces. This map is known as the cosmic string associated
 with this K3 fibration.

 Under the assumption that ${\bf M}$ is unimodular, our goal is to determine the 
virtual number of nodal rational curves in the classes of ${\bf M}$, extending
 the result of theorem \ref{theo; main} to the monodromy non-invariant classes.

 Because $\underline{\cal M}_{\bf M}$ is a quotient of a
 complex hyperbolic space, the domain curve
 $\tilde{B}$ of the holomorphic map $\Phi_{\cal Y}$ cannot be rational. 
Therefore the genus of $\tilde{B}$, $g(\tilde{B})$, has to be positive. 
In particular
 $\pi_1(\tilde{B}, b)\not=\{1\}$ generates a monodromy representation into 
$Aut({\bf L})=Aut(H^2(\pi^{-1}(b), {\bf Z}))$.

 It is well known that the whole arithmetic group $$Aut({\bf L})=
SO_{\bf Z}(3, 19)\equiv 
 SO_{\bf R}(3, 19)
\cap SL_{22}({\bf Z})$$ is generated by the $-2$ reflections
in ${\bf L}$. On the other hand, the sub-lattice ${\bf M}$ is 
monodromy invariant. Thus, the monodromy representation of
 $\pi_1(\tilde{B}, pt)$ into $Aut({\bf L})$ is induced by an arithmetic
 subgroup $\subset Aut({\bf M}^{\bot})$.

The fundamental class of a holomorphic curve in the fiber 
$\pi^{-1}(b')$, $b'\in \tilde{B}$
is in $H_2(\pi^{-1}(b'), {\bf Z})$, and
 therefore in $H_2({\cal Y}, {\bf Z})\cong H^4({\cal Y}, {\bf Z})$. 
Apparently such a class pairs trivially with $F\in H^2({\cal Y}, 
{\bf Z})$.

 Such a class $\in H^4({\cal Y}, {\bf Z})$ 
always lies in the image of $\cup F:H^2({\cal Y}, 
 {\bf Z})\mapsto H^4({\cal Y}, {\bf Z})$, if the $F-$cohomology
$Ker(\cup F)/Im(\cup F)$ is trivial.

 According to remark \ref{rem; cohomology}, the unimodular property of 
${\bf M}$ implies that all classes of holomorphic curves in the fibers 
${\bf \pi}^{-1}(b'), b'\in \tilde{B}$ 
can be viewed as the intersections of divisor classes in ${\bf M}$ with $F$. 

\begin{prop}\label{prop; mod}
 Two homology classes in $H_2({\bf \pi}^{-1}(b'), {\bf Z})$ are identified under $$(i_{b'})_{\ast}: H_2({\bf \pi}^{-1}(b'), {\bf Z})\mapsto H_2({\cal Y}, {\bf Z})$$
 if and only if there difference pairs trivially with 
 the embedded image of ${\bf M}$ into $H^2({\bf \pi}^{-1}(b'), 
{\bf Z})$.
\end{prop}

\noindent Proof: Because $\tilde{B}$ is connected, the image of
 $(i_{b'})_{\ast}$ does not depend on $b'\in B$.
The cokernel of the map $(i_{b'})_{\ast}$ consists of
 the homology classes which are not in the fibers. 
 According to Leray spectral sequence argument on $H_2$, 
it is one dimensional and is generated by the fundamental class of the
 cross section
 $s_{\ast}[\tilde{B}]$. Because $s_{\ast}[\tilde{B}]$ is Poincare dual to the
 fiber class $F$, the pairing between ${\bf M}\cong H^2({\cal Y}, 
{\bf Z})/{\bf Z}F$ and $Im((i_{b'})_{\ast})$ is perfect.
Then the result follows. $\Box$

\medskip

 Assuming that $M$ is unimodular,
 ${\bf L}={\bf M}\oplus {\bf M}^{\bot}$ is an orthogonal decomposition. 
 Under the pairing $H_2(\pi^{-1}(b'), {\bf Z})\otimes H^2(\pi^{-1}(b'), {\bf Z})
\mapsto {\bf Z}$, $(M^{\bot})^{\ast}\subset H_2(\pi^{-1}(b'), {\bf Z})\cong
 {\bf L}^{\ast}$ is exactly the kernel $Ker((i_{b'})_{\ast})$.
Our discussion shows that the fiberwise curve classes
 ($\in Im((i_{b'})_{\ast})$, $b'\in B$) are represented by
 cohomology classes ${\bf x}$ 
in ${\bf M}\cong {\bf L}/{\bf M}^{\bot}$, 
which parameterizes a whole equivalent class of 
elements ${\bf x}+{\bf M}^{\bot}$ in ${\bf L}$. Namely,
 the element ${\bf x}\in {\bf M}$ resembles an equivalence class of 
elements in ${\bf L}$ of the form 
$\{{\bf x}+{\bf y}$, ${\bf y}\in {\bf M}^{\bot}\}$.

 The following simple 
lemma relates the self-intersection pairings of ${\bf x}$ and
 ${\bf y}$.

\begin{lemm}\label{lemm; x}
 Let ${\bf x}+{\bf y}$ be a representative in ${\bf L}$ of the class ${\bf x}$, 
 In order that it contributes to the counting of algebraic nodal rational curves
 in ${\bf x}$, the self-intersection numbers ${\bf x}^2$ and ${\bf y}^2$ have to satisfy the following bound,

$${\bf x}^2+{\bf y}^2\geq -2.$$
\end{lemm}

\noindent proof:
  Given a class $C\in H^2({\cal Y}, {\bf Z})$ which 
restricts non-trivially to the fibers, 
the family expected dimension of the algebraic 
 family Seiberg-Witten invariants is given by 
${C^2-C\cdot c_1({\bf K}_{{\cal Y}/B})\over 2}+dim_{\bf C}dim B=
{C^2\over 2}+1$. In order that the
 class $C={\bf x}+{\bf y}$ has non-negative family dimension,

 $${\bf x}^2+{\bf y}^2={\bf x}^2+2{\bf x}\cdot {\bf y}+
{\bf y}^2=C^2\geq -2.$$ $\Box$ 

  This indicates that when ${\bf y}=0$, the element itself
 ${\bf x}$ has to satisfy ${\bf x}^2\geq -2$.

 If we restrict to a class ${\bf x}$, then
 the above lemma gives a lower bound on ${\bf y}^2$ immediately,

$${\bf y}^2\geq -{\bf x}^2-2.$$

 On the other hand, according to the adjunction formula 
the expected genus of the class ${\bf x}+{\bf y}$ in
 the algebraic $K3$ is given by 
$2g-2=C^2+c_1({\bf K}_{{\cal Y}/\tilde{B}})\cdot C=
{\bf x}^2+{\bf y}^2$. In order the class ${\bf x}+{\bf y}$ to 
 be represented by a rational curve in the fibers of 
${\cal Y}$, the immersed rational curve is expected to develop 
$g={{\bf x}^2+{\bf y}^2\over 2}+1$ nodes in the
 generic situation. 

\medskip

Because ${\bf y}$ is not invariant under the monodromy
 representation\footnote{Otherwise it would have been in 
 ${\bf M}$.}, the class 
${\bf x}+{\bf y}$ does not define a monodromy invariant class of 
the fiber bundle ${\cal Y}\mapsto \tilde{B}$.

 However, one can follow the ideas discussed in section 
\ref{subsection; extension} and extend the definition of the 
family invariant to take into account of the whole orbit of ${\bf x}+
{\bf y}$ instead (by taking an appropriate covering). 

 Thus the primary object we are interested at is

$${1\over [\tilde{B}, B]g!}{\cal AFSW}_{({\cal Y}/\tilde{B})_{g+1}
\mapsto({\cal Y}/\tilde{B})_g}^{\ast}(1, {\bf x}+{\bf y}-
\sum_{i\leq g}2E_i); g={{\bf x}^2+{\bf y}^2\over 2}+1$$

 According to a slight extension of 
theorem \ref{theo; main} to non-monodromic invariant $C$, 
the modified family invariant can be re-expressed to be 
the product of ${1\over [\tilde{B}, B]g!}
{\cal AFSW}_{{\cal Y}\mapsto \tilde{B}}(1, {\bf x}+{\bf y})$, and
 the universal degree $g$ polynomial in terms 
of the variables $C^2={\bf x}^2+{\bf y}^2$ and $c_2=24$. However, from remark 
\ref{rem; same}
we know that the
 universal polynomial is equal to $g!\cdot N_g$ (modulus a 
${\cal C}^{\infty}$ argument, 
where $N_g$ denotes the
 number of rational nodal curves dual to $C$ in an algebraic $K3$
 with the square $C^2=2g-2$. 

 Now we consider all the classes ${\bf x}+{\bf y}$ such that the self-intersection
 number
 ${\bf x}^2+{\bf y}^2=2g-2$ is fixed.
 As ${\bf x}^2+{\bf y}^2=2g-2$, then ${\bf y}^2=2g-2-{\bf x}^2$.
 The following lemma bounds the 
self-intersection number ${\bf y}^2$ from above.

\begin{lemm}\label{lemm; y}
 Let ${\cal Y}\mapsto B$ be a relative algebraic $K3$ fiber bundle and
 let ${\bf y}\in {\bf M}^{\bot}$. Suppose that
${\cal AFSW}_{{\cal Y}\mapsto B}(1, {\bf x}+{\bf y})\not=0$, 
then the class ${\bf y}$ has to satisfy
 ${\bf y}^2<0$.
\end{lemm}

\noindent Proof: By proposition \ref{prop; inf}, 
the non-vanishing of ${\cal AFSW}_{{\cal Y}\mapsto B}(1, 
{\bf x}+{\bf y})$ implies:

\medskip 

\noindent (1). The orbit of ${\bf x}+{\bf y}$, ${\bf O}_{{\bf x}+{\bf y}}$ 
under the monodromy group action of the fiber
bundle is a finite set.

\medskip

\noindent (2). Some element in ${\cal O}_{{\bf x}+{\bf y}}$ has to 
be represented by 
holomorphic curves in some fiber of ${\cal Y}$. If not, the
 algebraic family Seiberg-Witten invariant would have vanished!

\medskip

 The second point implies that 
there exists at least one $b_0\in \tilde{B}$ and an element $h\in 
\pi_1(\tilde{B}, b)$, such that
 $h({\bf x}+{\bf y})={\bf x}+h({\bf y})$ is represented by a holomorphic curve in the fiber ${\pi}^{-1}(b_0)$.
  In particular, ${\bf x}+h({\bf y})$ is of type $(1,1)$ in
 the fiber ${\bf \pi}^{-1}(b_0)$. 
On the other hand, ${\bf x}$ is in ${\bf M}$ and
 it is of type $(1, 1)$ in all the fibers of 
$\pi:{\cal Y}\mapsto \tilde{B}$. 
 So $h({\bf y})={\bf x}+h({\bf y})-{\bf x}$ is of type $(1, 1)$ 
 in the fiber ${\bf \pi}^{-1}(b_0)$. 

 Recall that by Hodge index theorem 
the lattice ${\bf M}$ is of signature $(1, m)$, for some $m\in {\bf N}$.  
 If ${\bf y}^2=(h({\bf y}))^2>0$, 
then the orthogonal direct sum 
${\bf M}\oplus {\bf Z} h({\bf y})$ would have been a
signature $(2, m)$ sub-lattice of $H^{1,1}({\bf \pi}^{-1}(b_0), 
{\bf C})$.  On the other hand, the relative ample polarization
 $\omega_{{\cal Y}/B}$ is monodromy invariant and it 
defines a 'big' element in ${\bf M}$.
By the so-called light cone lemma [LL2], the possibility ${\bf y}^2=0$ along with
 $\omega_{{\cal Y}/B}\cdot {\bf y}=\omega_{{\cal Y}/B}\cdot 
 h({\bf y})=0$ would have implied that $\omega_{{\cal Y}/B}\cdot
\omega_{{\cal Y}/B}=0$ as well. This is impossible as
 $\omega_{{\cal Y}/B}$ is relatively ample along the fibers.  
Therefore, ${\bf y}^2$ must be negative.  $\Box$ 

\medskip

 Lemma \ref{lemm; y} gives a lower bound on ${\bf x}^2$:
 $0>2g-2-{\bf x}^2$.

\bigskip

 In the following, we determine the total contribution of the
 whole equivalence class ${\bf x}+{\bf M}^{\bot}$ to the modified
 family invariant ${\cal AFSW}^{\ast}$, which will be identified (up to $g!$) 
with the virtual number of rational nodal curves of the whole orbit.

\medskip

 If the even lattice 
${\bf M}^{\bot}$ has been negative definite, 
then the enumeration of the family invariant
would have been related to the theta function of 
${\bf M}^{\bot}$ in the following simple way: 

Given a non-positive even number $2r$, $r\leq 0$, 
there are a finite
 number of elements in ${\bf M}^{\bot}$ with square 
$2r$. Let $n_{-r, {\bf M}^{\bot}}$ denote
the cardinality of $\{{\bf y}|{\bf y}\in {\bf M}^{\bot}, {\bf y}^2=2r\}$ for 
a negative definite ${\bf M}^{\bot}$.

 Then 

$${\Theta}_{{\bf M}^{\bot}}(q)\equiv \sum_{r\leq 0} n_{-r, {\bf M}^{\bot}} q^{-r}$$

 is the theta function associated to the lattice.
 Under the un-realistic assumption of a negative definite
 ${\bf M}^{\bot}$, the modified 
family invariant ${1\over [\tilde{B}, B]g!}{\cal AFSW}_{({\cal Y}/\tilde{B})_{g+1}
\mapsto({\cal Y}/\tilde{B})_g}(1, {\bf x}+{\bf y}-
\sum_{i\leq g}2E_i)$ of ${\bf x}+{\bf M}^{\bot}$
would have been
 reduced to 

$${1\over [\tilde{B}, B]}\sum_{{\bf y}\in {\bf M}^{\bot}} {\cal AFSW}_{
{\cal Y}\mapsto \tilde{B}}(1, {\bf x}+{\bf y})
\cdot N_{{{\bf x}^2+{\bf y}^2\over 2}+1}$$

$$=-{\int_{\tilde{B}}{\varpi}_{wp}\over [\tilde{B}, B]}
\sum_r N_{{{\bf x}^2\over 2}+r+1}\cdot n_{-r, {\bf M}^{\bot}}.$$

 We have made use of the fact that ${\cal AFSW}_{{\cal Y}\mapsto
 \tilde{B}}(1, C)=-\int_{\tilde{B}} \varpi_{wp}$ from section \ref{section; 
cosmic}.
 Because the above expression
 depends on ${\bf x}$ only through ${\bf x}^2$, one 
 may cast it into a generation function 

$$-\sum_{k\geq -1} {\int_{\tilde{B}}
{\varpi}_{wp}\over [\tilde{B}, B]}\{\sum_r N_{k+r+1}\cdot
 n_{-r, {\bf M}^{\bot}}\}q^{k+1}$$

$$=-{\int_{\tilde{B}}{\varpi}_{wp}\over [\tilde{B}, B]}\sum_k \{\sum_r N_{k+r+1}
\cdot n_{-r, {\bf M}^{\bot}} \}q^{k+1}.$$

 By a change of variable $k+r+1=\delta$, the above expression can be
 factorized into 

$$=-{\int_{\tilde{B}}{\varpi}_{wp}\over [\tilde{B}, B]}\{\sum_{\delta} 
N_{\delta} q^{\delta}\}\cdot \{\sum_r n_{-r, {\bf M}^{\bot}}q^{-r}\}$$
\label{negam}
$$=-({\int_{\tilde{B}}{\varpi}_{wp}\over [\tilde{B}, B]})
\{\sum_{\delta} N_{\delta}q^{\delta}\}
 \Theta_{{\bf M}^{\bot}}(q).$$ 

\bigskip

 Apparently the ${\bf M}^{\bot}$ is never negative
 definite. Instead, ${\bf M}^{\bot}$ is of signature $(2, 19-m)$.  Thus the
 traditional theta function $\Theta_{{\bf M}^{\bot}}$ blows up term by term 
as $n_{-r}=\infty$ for all $r$.

 In the following, we develop a method from Howe duality 
and intersection theory on $\underline{\cal M}_{\bf M}$ 
to regulate the $\Theta_{{\bf M}^{\bot}}(q)$ into
a power series.

\bigskip

\subsection{The Renormalization of the Theta Function of
 ${\bf M}^{\perp}$}\label{subsection; theta}

\medskip

 Firstly we decompose the indefinite lattice ${\bf M}^{\bot}$ into the disjoint unions of the various orbits
 of the Monodromy group action $Im(\pi_1(\tilde{B}, b))
\subset Aut({\bf M}^{\bot})$,

  $${\bf M}^{\bot}=\coprod_{[{\bf y}]} {\bf O}_{\bf y},$$ 

 where the disjoint union
 is over the equivalence classes $[{\bf y}]$ of the orbits of lattice elements.

 By lemma \ref{lemm; y}, we can discard all the orbits whose self-intersection
 numbers are non-negative.

 Thus we cast the modified family invariant into the following sum,

$$\sum_{[{\bf y}], {\bf y}^2<0} 
{1\over [\tilde{B}, B]}{\cal AFSW}_{{\cal Y}\mapsto 
\tilde{B}}(1, {\bf x}+{\bf O}_{{\bf y}})\cdot
 N_{{{\bf x}^2+{\bf y}^2\over 2}+1}.$$

 The following simple lemma guarantees the finiteness of the sum,

\begin{lemm}\label{lemm; u}
 In the above formal sum of algebraic 
family Seiberg-Witten invariants, all but a finite number of 
terms vanish and the finite sum depends on
 ${\bf x}$ only through ${\bf x}^2$.
\end{lemm}

\noindent Proof: If it is an infinite sum, then there will be an infinite number of orbits
 ${\bf O}_{\bf y}$ with non-zero family invariants.
 Fixing the ${\bf x}$, lemma \ref{lemm; x} and lemma \ref{lemm; y} imply that
 these ${\bf y}$ with non-zero family invariants 
satisfy the constraint $-2\geq {\bf y}^2\geq -2-{\bf x}^2$. Therefore there exists at least one fixed $r$ 
such that there are an infinite number of orbits ${\bf O}_{\bf y}$
 with the fixed ${\bf y}^2=2r$, 

$${\cal AFSW}_{{\cal Y}\mapsto \tilde{B}}(1, {\bf x}+{\bf O}_{\bf y})
\not=0.$$

 This implies that for all these infinite number of
${\bf y}$, all ${\bf x}+{\bf O}_{\bf y}$ from different
 orbits
 are represented by
 algebraic curves within the $\tilde{B}$ family. Therefore for each
 ${\bf y}$, the whole ${\bf O}_{\bf y}$ become of type
 $(1, 1)$ in the same fibers (depending on ${\bf y}$).  
Because $\tilde{B}$ is compact, one can find an accumulation point
 $b_{ac}\in \tilde{B}$ such that the curves accumulate into 
the fiber $\pi^{-1}(b_{ac})$
 above $b_{ac}$. Now we restrict the fiber bundle ${\cal Y}\mapsto \tilde{B}$ 
to a small neighborhood of $b_{ac}$. 

One may adjoin all these infinite numbers of orbits
${\bf O}_{{\bf y}_i}$, ${\bf y}_i^2=2r<0$ (of 
type $(1, 1)$ in the fibers somewhere nearby $\pi^{-1}(b_{ac})$), to
 ${\bf M}$ and form a hyperbolic (i.e. of signature $(1, k)$) sub-lattice of
 ${\bf L}$. 

\begin{lemm}\label{lemm; hyper}
The sub-lattice ${\bf I}$ of ${\bf L}$ formed by adjoining all these ${\bf y}_i$
(which are of type $(1, 1)$ along the fibers nearby $\pi^{-1}(b_{ac})$)
 to ${\bf M}$ is hyperbolic.
\end{lemm}

\noindent Proof: If the positive eigen-space of the sub-lattice ${\bf I}$ 
has been greater than one dimensional, 
then one may find an element $\eta\in {\bf I}$ in
 the positive eigenspace of ${\bf I}$ with
 $\eta^2>0$, $\eta\omega_{{\cal Y}/\tilde{B}}=0$. As the relative polarization 
$\omega_{{\cal Y}/\tilde{B}}\in {\bf M}$
 is of type $(1, 1)$ along all the fibers, such an $\eta$ must be
expressed as linear combinations of $(2, 0)$ and $(0, 2)$ forms along all the
 nearby fibers of $\pi^{-1}(b_{ac})$.  
On the other hand, $\eta\in {\bf I}$ is a linear combination
 of elements from ${\bf M}$ and a finite number of ${\bf y}_i$, each of them
 is of type $(1, 1)$ along some fiber(s). By using these facts,
 $\eta^2=\eta\cdot \eta$ is a combination of pairings between $(2, 0)$ (or
 $(0, 2)$) class and $(1, 1)$ classes. So $\eta^2=0$, violating the
 fact that $\eta$ is in the positive eigen-space of ${\bf I}$.   $\Box$

 Even though all these ${\bf O}_{\bf y}$s do not
form a single orbit under the monodromy group action, the argument
 of proposition \ref{prop; inf} is still applicable and it implies the
existence of infinite number of lattice elements in a compact set of the
forward light-cone. We can still derive a contradiction. $\Box$

\begin{rem}\label{rem; Sachs}
Because all these classes ${\bf x}+{\bf O}_{\bf y}$ have a fixed 
harmonic energy\footnote{It is determined by ${\bf x}$ only
 and has nothing to do with ${\bf O}_{\bf y}\subset {\bf M}^{\perp}$.}, 
it also violates
 Gromov-Sachs-Uhlenback compactness theorem in symplectic topology [MS].
\end{rem}

 Given a finite orbit ${\bf O}_{\bf y}$, 
 one can take the finite covering $\hat B\mapsto \tilde{B}$ and make
 the class ${\bf y}$ monodromy invariant on the pull-back 
fiber bundle ${\cal Y}\times_{\tilde{B}} \hat{B}\mapsto \hat{B}$.
Then the family invariant
 ${\cal AFSW}_{{\cal Y}\mapsto \tilde{B}}(1, {\bf x}+{\bf O}_{\bf y})$ 
is equal to
 $=-{1\over [\hat{B}, \tilde{B}]}\int_{\hat{B}}\Phi_{\cal Y}^{\ast}
\varpi_{wp}$ $=-\int_{\tilde{B}}\Phi_{\cal Y}^{\ast}
\varpi_{wp}$,
 which is independent of ${\bf x}$. 
Thus, the only dependence of the finite sum 

\noindent $\sum_{[{\bf y}], {\bf y}^2<0} 
{1\over [\tilde{B}, B]}{\cal AFSW}_{{\cal Y}\mapsto 
\tilde{B}}(1, {\bf x}+{\bf O}_{{\bf y}})\cdot
 N_{{{\bf x}^2+{\bf y}^2\over 2}+1}$ on ${\bf x}$ is
 through the inequality 
$-2\geq {\bf y}^2\geq -2-{\bf x}^2$ and $N_{{{\bf x}^2+{\bf y}^2\over
 2}+1}$. $\Box$ 

\bigskip

 To enumerate the above finite sum, we cast it into an expression
 with a different flavor.
 Recall that we had introduced in section \ref{section; cosmic} a 
 line bundle $\Phi_{\cal Y}^{\ast}
{\bf U}\mapsto \tilde{B}$ by pulling back the
 universal line bundle (of holomorphic two forms on K3) over 
$\underline{\cal M}_{\bf M}$ using the cosmic string map
 $\Phi_{\cal Y}:\tilde{B}\mapsto \underline{\cal M}_{\bf M}$. 
In the following, we first construct a section
 of $\Phi_{\cal Y}^{\ast}{\bf U}$ on $\tilde{B}$ by using the element ${\bf y}$.

Let $\pi^{-1}(b')$ be a fiber K3 surface.
Given a finite orbit ${\bf O}_{\bf y}$
 and a holomorphic two form ${\bf \Omega}_{b'}$ 
over $b'$, we consider
 the following period integral pairing
$$ {\int_{{\bf \pi}^{-1}(b')}{\bf y}\cup{\bf \Omega}_{b'}}.$$

We extend this into the whole $\tilde{B}$ family by the following recipe.
Consider the complex hyperbolic space ${\bf P}$ covering
 $\underline{\cal M}_{\bf M}$, parametrizing the complex
deformations of an algebraic K3 surface $M$ with Picard lattice
 $\supset {\bf M}$. Because the deformation space ${\bf P}$ is contractible, 
 one can find a global ${\bf P}$-family of relative holomorphic two forms ${\bf \Omega}$ on the K3s with the normalization $\int_M {\bf \Omega}\wedge
 \bar{\bf\Omega}=1$.

\begin{lemm}\label{lemm; des}
 The assignment 
$$b\mapsto {\int_{{\bf \pi}^{-1}(b)}}{\bf y}\cup{\bf \Omega}|_b$$
 descends to a ${\cal C}^{\infty}$ section of 
$\Phi_{\cal Y}^{\ast}
{\bf U}^{\ast}$ over $\tilde{B}$.
\end{lemm} 

\medskip

\noindent Proof:
 Consider the universal covering $\check{B}$ of $\tilde{B}$, the cosmic string
 map $\tilde{B}\mapsto \underline{\cal M}_{\bf M}$ can be lifted to
 a map from $\check{B}$ to ${\bf P}$.
 By pulling back the 
 smooth cross section ${\bf \Omega}$ of ${\bf U}$ 
(which trivializes ${\bf U}$) and the tautological $K3$ family to 
 $\check{B}$, the assignment 
in the above lemma determines a smooth section of 
$\Phi_{\cal Y}^{\ast}{\bf U}^{\ast}$
 over $\check{B}$.
 Since ${\bf O}_{\bf y}$ is a 
finite orbit under $\pi_1(\tilde{B}, b)$, the above section can 
 be descended to one over the quotient $\tilde{B}$ of $\check{B}$.  $\Box$ 

\medskip

  The section defined above is called the period section of ${\bf O}_{\bf y}$.

\begin{prop}\label{prop; euler}
 The algebraic family invariant 
${\cal AFSW}_{{\cal Y}\mapsto \tilde{B}}(1, {\bf O}_{\bf y})$ 
can be identified with the number of zeros (counted with multiplicities) 
of the above "period section".
\end{prop} 

\noindent Proof of prop. \ref{prop; euler}:  
By the previous calculation in section \ref{section; cosmic}, 

\noindent ${\cal AFSW}_{{\cal Y}\mapsto \tilde{B}}(1, {\bf O}_{\bf y})=
-\int_{\tilde{B}}\varpi_{wp}$. On the other hand, the cohomology class of
 the Weil-Peterson form $[\varpi_{wp}]$ is $c_1(\Phi_{\cal Y}^{\ast}{\bf U})$. Thus
 the above family invariant can be identified with the 
Euler number of the underlying real rank two bundle 
of $\Phi_{\cal Y}^{\ast}{\bf U}^{\ast}$ over $\tilde{B}$, 
which is equal to the number of zeros of the 
"period section" constructed in lemma \ref{lemm; des}.  
Because the zeros of the "period section"
 may not be transversal, 
the number has to be counted with multiplicities. $\Box$

\medskip

 The zeros of the "period section" can be 
interpreted alternatively by the following lemma,

\begin{lemm}\label{lemm; period}
 The period section vanishes over a point $b'\in \tilde{B}$ 
if and only if the elements in ${\bf O}_{\bf y}$ are of type $(1, 1)$ in
 the middle cohomology $H^2(\pi^{-1}(b'), {\bf C})$.
\end{lemm} 
\medskip

\noindent Proof of lemma \ref{lemm; period}: Take $M=\pi^{-1}(b')$.
Through the Hodge decomposition

 $$H^2(M, {\bf C})=H^{2, 0}(M, {\bf C})\oplus
 H^{1, 1}(M, {\bf C})\oplus H^{0,2}(M, {\bf C}),$$

the elements in ${\bf O}_{\bf y}$ can be decomposed into different 
components according to their Hodge types.
An element ${\bf y}\in H^2(M, {\bf Z})$ can be 
 decomposed into ${\bf y}={\bf y}^{1,1}\oplus {\bf y}^{2, 0}\oplus 
{\bf y}^{0, 2}$ with $\overline{{\bf y}^{0, 2}}={\bf y}^{2, 0}$.

  It is easy to see that $\int_M{\bf \Omega}|_{b'}\cup
{\bf y}^{1,1}=0$, as the cup product preserves the Hodge type. On the other hand,
 the reality condition on ${\bf y}^{2,0}+{\bf y}^{0,2}$ implies that
 it can be written as $c{\bf \Omega}|_{b'}+\bar c\bar 
{\bf \Omega}|_{b'}$ for some
 $c\in {\bf C}$.  Thus the above period integral becomes
 $\bar c\int_M\bar{\bf \Omega}|_{b'}
\cup{\bf \Omega}|_{b'}=0$.

 This implies that $\bar c=0$ and 
therefore ${\bf y}^{2,0}+ {\bf y}^{0,2}=0$. So the element 
${\bf y}$ is of type $(1, 1)$ in $H^2(M, {\bf Z})$.  The converse of the
 above assertion holds by a simple calculation.

 Finally one applies the above argument to the finite collection of 
elements in ${\bf O}_{\bf y}$ instead of a single monodromy invariant ${\bf y}$.
$\Box$

 Given the element ${\bf z}\in {\bf M}^{\bot}$, ${\bf z}^2<0$, we consider the 
codimension one complex hyperbolic subspace ${\bf P}_{\bf z}\subset
 {\bf P}$ defined by the 
zero locus $\ni b'$ of $\int_{\pi^{-1}(b')}{\bf z}\cup {\bf \Omega}$. According to
 lemma \ref{lemm; period}, it is the locus over which ${\bf z}$ 
becomes of type $(1, 1)$. 
It turns out that this complex hyperbolic 
subspace is stable under a subgroup of finite index of the modular group 
$Aut({\bf M}^{\bot})$.

The quotient of ${\bf P}_{\bf z}$ under the finite index subgroup 
 can be identified with the moduli space of 
${\bf M}\oplus {\bf Z}{\bf z}$-marked
 K3 surfaces, 
$\underline {{\cal M}}_{{\bf M}\oplus {\bf Z}{\bf z}}$, 
 which is mapped into $\underline{{\cal M}}_{\bf M}$ 
as a Weil divisor, denoted 
 by ${\cal D}_{\bf z}$.

By combining
 the above discussion with lemma \ref{lemm; period}, 
we find that
 ${\cal AFSW}_{{\cal Y}\mapsto \tilde{B}}(1, {\bf O}_{\bf y})$ 
can be re-casted into

 $$\sum_{{\bf z}\in {\bf O}_{\bf y}}\sharp 
\{{\cal D}_{\bf z}\cap \Phi_{{\cal Y}\ast}[\tilde{B}]\},$$ 

 where $\sharp\{{\cal D}_{\bf z}\cap \Phi_{{\cal Y}\ast}[\tilde{B}]\}$ 
denotes the signed intersection numbers of the
 divisors and the image of the cosmic string curve 
$\Phi_{\cal Y}(\tilde{B})$!
 So the above family algebraic Seiberg-Witten invariants are identified
 with the sum of the intersection numbers of the
 special Weil divisors ${\cal D}_{\bf z}$ and
the image of the cosmic string map.

\medskip

 Define the regulated theta function ${\Theta}_{reg}(q)$ of ${\bf M}^{\bot}$ to 
be,

\begin{defin}\label{defin; regularize}
$${\Theta}_{{\bf M}^{\bot}}^{reg}(q)
=1-\sum_{r>0}\sum_{{\bf O}_{\bf y}\subset {\bf M}^{\bot};
{\bf y}^2=-2r}{1\over \int_B\Phi_{\cal Y}^{\ast}\varpi_{wp}}
  \sum_{{\bf z}\in {\bf O}_{\bf y}}
\sharp \{{\cal D}_{\bf z}\cap \Phi_{{\cal Y}\ast}[B]\} q^r.$$
\end{defin}
  
Lemma \ref{lemm; u} guarantees the finiteness of the
 individual coefficients of the regularized theta series.

 Then by following a similar calculation as the 
negative definite ${\bf M}^{\bot}$ on page \pageref{negam}, 
the total contribution of ${\bf M}^{\bot}$ to the virtual number of
 family invariant is equal to

 $$-\int_{\tilde{B}}\Phi^{\ast}_{\cal Y}\varpi_{wp}\cdot {\Theta}_{{\bf 
M}^{\bot}}^{reg}(q)\cdot (\sum_{g\in {\bf N}\cup \{0\}} N_g q^g).$$

\begin{rem}
 In the above discussion, we have taken $C^2=2g-2$. If we choose $g$ and $l$
such that $C^2=2(g+l)-2$ and consider the following modified family invariant

$${1\over [\tilde{B}, B]l!}{\cal AFSW}_{({\cal Y}/\tilde{B})_{l+1}\mapsto 
({\cal Y}/\tilde{B})_l}^{\ast}(\eta_g, {\bf x}+{\bf y}-2\sum_{i\leq l}E_i),$$

 then a parallel discussion leads to the parallel formula,

$$-\int_{\tilde{B}}\Phi^{\ast}_{\cal Y}\varpi_{wp}\cdot {\Theta}_{{\bf 
M}^{\bot}}^{reg}(q)\cdot (\sum_{l\in {\bf N}\cup \{0\}} N_l(g) q^l).$$

Here $N_l(g)$ denote the virtual number of genus $g$ $l$-node nodal curves
dual to $C$ ($C^2=2(g+l)-2$) 
 on an algebraic $K3$, with $N_g(0)=N_g$. Notice that $N_l(g)$ is distinct from
 the usual virtual number of $l$-node nodal curves through $g$-generic points in
 the $K3$.
\end{rem}

\subsection{The Howe Duality and Type $II$-Heterotic Duality}
\label{subsection; howe}

\bigskip

  In string theory the counting of rational curves from Calabi-Yau
K3 fibrations are predicted by Harvey and Moore [HM1] using
 type $IIA$-heterotic duality to heterotic string on $K3\times T^2$. In our 
mathematical discussion, the special string duality in string theory
 can be traced to the Howe duality on metaplectic representations.

 By K$\ddot{u}$nneth theorem, the middle cohomology $H^3(K3\times T^2, 
{\bf Z})$, equipped with a symplectic intersection pairing, can
be identified with $H^2(K3, {\bf Z})\otimes H^1(T^2, {\bf Z})$,
 where $H^2(K3, {\bf Z})$ and $H^1(T^2, {\bf Z})$ are equipped
 with an orthogonal and a symplectic intersection pairing, respectively.  It turns
 out that modulo the heavy machineries from Seiberg-Witten
 theory [Liu1], [Liu2], [Liu4], [Liu5], [Liu6], the Howe duality between the meta-plectic 
representations of $SL_2({\bf C})$ and $SO(p, 2)$ are responsible
 for the modularities of the curve counting! 

\medskip

 We are ready to prove the following main theorem,

\medskip

\begin{theo} \label{theo; modular}
 Let ${\cal Y}\mapsto \tilde{B}$ be an algebraic $K3$ 
fiber bundle structure over a smooth
 algebraic curve $\tilde{B}$ and let ${\bf M}$$=H^2({\cal Y}, {\bf Z})_f$
$=H^2({\cal Y}, {\bf Z})/{\bf Z}F$ be an unimodular lattice. 
 Let ${\cal F}(q)$ denote the formal power series of 
normalized virtual
 numbers of modified family invariants, resembling the virtual numbers
 of nodal rational curves along the fibration\footnote{Our earlier
discussion has shown that these numbers depend on $l={{\bf x}^2\over 2}+1$, 
and not the details of
the classes ${\bf x}\in {\bf M}$.}, then
 it can be factorized into the following form,

$$-\int_{\tilde{B}}\Phi_{\cal Y}^{\ast}(\varpi_{wp})\cdot 
\{{1\over \prod_{i\geq 0} (1-q^i)
}\}^{24}\cdot \Theta_{{\bf M}^{\bot}}^{reg}(q), $$
where $\Theta_{{\bf M}^{\bot}}^{reg}(q)$ is a holomorphic 
$SL_2({\bf Z})$-modular form of weight 
${rank_{\bf Z}{\bf M}^{\bot}\over 2}$. 
\end{theo}

\bigskip

 If ${\bf M}$ is isomorphic to either ${\bf H}$ or 
${\bf H}\oplus -{\bf E}_8$,  then ${\bf M}^{\bot}$ is isomorphic to 
$2{\bf H}\oplus -2{\bf E}_8$ or $2{\bf H}\oplus -{\bf E}_8$, respectively.
 Then the weights of the modular forms $\Theta^{reg}_{{\bf M}^{\bot}}(q)$ are
 $10$ and $6$ respectively.   We cite the following simple fact.

\begin{lemm}\label{lemm; 6, 10}
 The vector space of holomorphic $SL_2({\bf Z})$ modular forms 
 of weight $3$ or $5$ is one dimensional and is generated by
 $E_6(q)$ or $E_4\cdot E_6(q)$, respectively.
\end{lemm}

 Here $E_k(q)=1+{-2k\over B_k}
\sum_n \sigma_{k-1}(n)q^n$ is the $k$-th Eisenstein series, where $B_k$ is
 the $k$-th Bernoulli number.
We have $E_4(q)=1+240\sum_n \sigma_3(n)q^n$, $E_6(q)=1-504\sigma_5(n)q^n$.

\noindent Proof: It is well known that $E_4(q)$ and $E_6(q)$ generate
 the ring of $SL_2({\bf Z})$ modular forms. So the lemma follows from
 a simple weight count. $\Box$

  Then one can identify 
$\Theta_{{\bf M}^{\bot}}^{reg}(q)$ uniquely, by 
the condition $\Theta_{{\bf M}^{\bot}}^{reg}(0)=1$. 

\medskip

\noindent Proof of Theorem \ref{theo; modular}: By
 using Taubes argument of "SW=Gr" to the universal families, we
can argue that $\sum_{g\geq 0} 
N_g q^g=\bigl({1\over \prod_{i\geq 1}(1-q^i)}\bigr)^{24}$.

 The factorization of the generating function of
 modified invariants into the product of
${q\over \Delta(q)}=
\bigl({1\over \prod_{i\geq 0} (1-q^i)}\bigr)^{24}$ and the regulated theta function
 $\Theta_{{\bf M}^{\bot}}^{reg}(q)$ has been discussed. 
Our goal is to prove that $\Theta_{{\bf M}^{\bot}}^{reg}(q)$ 
is a holomorphic modular form of
the modular group $SL_2({\bf Z})$ of weight ${rank_{\bf Z}{\bf M}^{\bot}\over 2}$.

 This is the place where type $II$-heterotic duality (and
 its mathematical analogue) appears
 implicitly in our picture. 
 We recast the regulated theta function into a different expression,
 whose modularity is proved by using a beautiful result of Kudla-Millson [KM] on
 Howe duality.

\medskip

 In the following, we recall the result, following the original
notations of Kudla-Millson [KM].

\medskip

  Take $V$ be a real vector space of dimension $m$ and
 $(\cdot, \cdot)$ be a quadratic form on $V$ of signature $(p, q)$.

  Let $\Gamma\backslash D$ be the arithmetic quotient of the symmetric
 space $D$ of ${\bf O}(p, q)$.  Let $\beta$ be a symmetric $n\times n$
 matrix with ${\bf Z}$ coefficients which is positive semi-definite.
 Let $C_{\beta}$ denote the special cycle of dimension
 $(p-t)q$ in $\Gamma\backslash D$ constructed
 in [KM], where $t$ denotes the rank of $\beta$. 
 Let $\eta$ denote a closed rapidly decreasing $(p-n)q$ differential form
 on $\Gamma\backslash D$. 

\medskip

 Following the original notations of Kudla-Millson [KM], 
take ${\bf \tau}\in M_n({\bf C})$ to be the
 variable in the Siegel upper-half space.  Define the power series
 $$P({\bf \tau}, \eta)=\sum_{t=0}^n \sum_{\beta\in {\cal L}(t)}
(\int_{C_{\beta}}\eta\wedge e^{n-t}_q)e_{\ast}(\beta{\bf \tau}),$$
where $e_{\ast}(\beta{\bf \tau})$ is defined to be $e^{{1\over 2}
 tr(\beta{\bf \tau})}$.

 The main theorem of Kudla-Millson, on page 126 of [KM], 
generalizing Shintani's result [Shin] is the following,

\begin{theo}\label{theo; KM}
The function $P({\bf \tau}, \eta)$ is a holomorphic modular form of weighted
$m/2$ for a suitable congruence subgroup of $Sp_n({\bf Z})$.
\end{theo}

 We have specialized the ring of integers ${\cal O}$ to ${\bf Z}$.

  There is a homological version of the above 
theorem, on the same page, page 126 of [KM], which we will use.
 Let\footnote{The class $\tilde{C}$ has nothing to do with
 $C$ used frequently in the paper. We adopt this notation here
 as it has been used in the original paper [KM].}
 $\tilde C\subset \Gamma\backslash D$ be a compact cycle of real dimension $nq$. 
Define $${\bf I}({\bf \tau}, \tilde C)= 
\sum_{t=0}^n \sum_{\beta\in {\cal L}(t)}\tilde C\cdot
 (C_{\beta}\cap e^{n-t}_q)e_{\ast}(\beta{\bf \tau}).$$ 
\medskip

\begin{cor}\label{corr; KM}
 The generating function ${\bf I}({\bf \tau}, \tilde C)$ is a 
holomorphic modular form of weigh$m/2$ for a suitable
 congruent subgroup of $Sp_n({\bf Z})$. 
\end{cor}

\bigskip

 In our application, 
we take $V={\bf M}^{\bot}\otimes_{\bf Z}{\bf C}$. 
Then $p=2m-2$ and $q=2$ with $m=rank_{\bf Z}{\bf M}^{\bot}$. 
Take $n=1$ and
 take $\tilde C$ to be the image of the cosmic string 
$\Phi_{\cal Y}(\tilde{B})$.  In this case, when $n=1$ the Siegel upper-half
space is reduced to the usual upper half plane, 
where ${\bf \tau}\in {\bf C}$,
 $Im({\bf \tau})>0$.

Then ${\bf I}({\bf \tau}, \tilde C)$ is reduced
 to $$\int_{ \Phi_{\cal Y}(\tilde{B})}e+\sum_{\beta\in {\cal L}(1)}
\{\Phi_{{\cal Y}\ast}[\tilde{B}]\cap(C_{\beta})\}e_{\ast}(\beta{\bf \tau}).$$

It is not hard to see from the proof of their
 paper that the congruent subgroup is this special case is
 $Sp_1({\bf Z})=SL_2({\bf Z})$ itself.

\medskip

The rank $1$ symmetric semi-positive definite
 integral matrix $\beta$ are nothing but non-negative even 
integers, which can
 be realized as the negation of the Gram matrix of negative self-intersecting 
lattice elements in ${\bf M}^{\bot}$ or $\{0\}$.  Then the locally finite cycle
 $C_{\beta}=C_{2k}; k\in {\bf N}$, can be re-written as the
 union (in terms of our notations)
$\cup_{{\bf z}, {\bf z}^2=-2k}{\cal D}_{{\bf z}}$.

 On the other hand, the universal $2$-plane bundle is isomorphic to
the underlying real bundle of ${\bf U}^{\ast}$. So 
the Euler class of the universal $2$-plane bundle ($q=2$)
 on $\Gamma\backslash D$ is equal to the negation of the class of 
Weil-Peterson form $-[\varpi_{wp}]$.

 We see without difficulty that ${\bf I}({\bf \tau}, \tilde C)$ match
 perfectly with our
 $-\int_{\Phi_{\cal Y}(\tilde{B})}
\varpi_{wp}\cdot \Theta_{{\bf M}^{\bot}}^{reg}(q)$.

\medskip

 Even though the locally finite cycle 
$C_{\beta}$ is not a 'finite' sum of irreducible cycles,
 lemma \ref{lemm; u} along with proposition \ref{prop; euler} and
 lemma \ref{lemm; period} imply that its intersection number with
 the cosmic string image is still 
 well defined (finite). Therefore the above theorem follows from combining the 
main theorem of Kudla-Millson [KM] with our earlier discussion.
$\Box$

{}

\end{document}